\documentclass[10pt,a4paper]{article}

\usepackage{bbm}
\usepackage{graphicx}
\usepackage{latexsym}
\usepackage{amsmath, amsfonts,amssymb, amsthm, euscript,makeidx, xcolor,mathrsfs}
\usepackage{enumerate}
\usepackage{tabularx}

\usepackage[colorlinks,linkcolor=blue,anchorcolor=green,citecolor=red]{hyperref}

\def\sqr#1#2{{\vcenter{\vbox{\hrule height.#2pt
              \hbox{\vrule width.#2pt height#1pt \kern#1pt \vrule width.#2pt}
              \hrule height.#2pt}}}}
\def\5n{\negthinspace \negthinspace \negthinspace \negthinspace \negthinspace }
\def\4n{\negthinspace \negthinspace \negthinspace \negthinspace }
\def\3n{\negthinspace \negthinspace \negthinspace }
\def\2n{\negthinspace \negthinspace }
\def\1n{\negthinspace }

\def\dbD{\mathbb{D}}

\def\dbF{\mathbb{F}}

\def\dbH{\mathbb{H}}

\def\dbN{\mathbb{N}}

\def\dbP{\mathbb{P}}

\def\dbR{\mathbb{R}}

\def\dbW{\mathbb{W}}

\def\={\buildrel \triangle \over =}

\def\ds{\displaystyle}

\def\ns{\noalign{\ss}}
\def\a{\alpha}

\def\e{\varepsilon}

\def\si{\sigma}
\def\t{\tau}

\def\th{\theta}

\def\G{\Gamma}
\def\D{\Delta}

\def\me{\mathbb{E}}
\def\rd{\,\mathrm d}
\def\bal{\begin{aligned}}
\def\eal{\end{aligned}}
\def\cC{{\cal C}}

\def\cE{{\cal E}}
\def\cF{{\cal F}}
\def\cG{{\cal G}}
\def\cH{{\cal H}}

\def\cN{{\cal N}}

\def\cS{{\cal S}}

\def\Bj{{\bf j}}

\def\no{\noindent}

\def\ss{\smallskip}
\def\ms{\medskip}
\def\bs{\bigskip}
\def\q{\quad}
\def\qq{\qquad}

\def\limsup{\mathop{\overline{\rm lim}}}

\def\lt{\left}
\def\rt{\right}
\def\lan{\langle}
\def\ran{\rangle}

\def\wt{\widetilde}

\def\cd{\cdot}
\def\cds{\cdots}

\def\as{\hbox{\rm a.s.}}

\DeclareMathOperator{\Span}{span}

\def\deq{\triangleq}

\def\({\Big (}
\def\){\Big )}
\def\pr{{\partial}}
\def\bde{\begin{definition}\label}
\def\ede{\end{definition}}
\def\be{\begin{equation}}
\def\bel{\begin{equation}\label}
\def\ee{\end{equation}}
\def\beq{\begin{equation*}\begin{aligned}}
\def\eeq{\end{aligned}\end{equation*}}
\def\bt{\begin{theorem}\label}
\def\et{\end{theorem}}
\def\bc{\begin{corollary}\label}
\def\ec{\end{corollary}}
\def\bl{\begin{lemma}\label}
\def\el{\end{lemma}}
\def\bp{\begin{proposition}\label}
\def\ep{\end{proposition}}
\def\bas{\begin{assumption}\label}
\def\eas{\end{assumption}}
\def\br{\begin{remark}\label}
\def\er{\end{remark}}
\def\bex{\begin{example}\label}
\def\ex{\end{example}}
\def\ba{\begin{array}}
\def\ea{\end{array}}
\def\ed{\end{document}}

\def\square#1{\vbox{\hrule\hbox{\vrule height#1%
     \kern#1\vrule}\hrule}}
\def\rectangle#1#2{\vbox{\hrule\hbox{\vrule height#1%
     \kern#2\vrule}\hrule}}

\font\tenbb=msbm10 \font\sevenbb=msbm7 \font\fivebb=msbm5

\newfam\bbfam
\scriptscriptfont\bbfam=\fivebb \textfont\bbfam=\tenbb
\scriptfont\bbfam=\sevenbb

\allowdisplaybreaks  

\newtheorem{theorem}{\hskip 1.3em Theorem}[section]
\newtheorem{definition}[theorem]{\hskip 1.3em Definition}
\newtheorem{proposition}[theorem]{\hskip 1.3em Proposition}
\newtheorem{corollary}[theorem]{\hskip 1.3em Corollary}
\newtheorem{lemma}[theorem]{\hskip 1.3em Lemma}
\newtheorem{remark}[theorem]{\hskip 1.3em Remark}
\newtheorem{example}[theorem]{\hskip 1.3em Example}
\newtheorem{algorithm}[theorem]{\hskip 1.3em Algorithm}

\newtheorem{assumption}[theorem]{\hskip 1.3em Assumption}

\makeatletter
   
   \@addtoreset{equation}{section}
\makeatother

\begin{document}

\title{Numerical Analysis of the Finite Transposition Method for Backward Stochastic Differential Equations\thanks{
{This work was partially supported by the National Key R\&D Program of China under grant 2024YFA1013400, the New Cornerstone Investigator Program, 
the National Natural Science Foundation of China under Grants 12271298 and 12671156, and the Taishan Scholars Project under Grant TSTP 20250702.}\ms}}

\author{Penghui Wang,\thanks{
School of Mathematics, Shandong University, Jinan 250100, China. {\small\it
e-mail:} {\small\tt phwang@sdu.edu.cn}. \ms}
~~~Yanqing Wang,\thanks{
School of Mathematics and Statistics, Southwest University, Chongqing 400715, China.  {\small\it
e-mail:} {\small\tt yqwang@swu.edu.cn}. \ms}~~~
  and~~~
Xu Zhang\thanks{Corresponding author,
School of Mathematics, Sichuan University, Chengdu 610064, China. {\small\it
e-mail:} {\small\tt zhang$\_$xu@scu.edu.cn}.}}

\date{}
\maketitle

\begin{abstract}
In this paper, we develop a {\em finite transposition method} for the numerical approximation of backward stochastic differential equations (BSDEs). Starting from the transposition formulation of BSDEs, we construct appropriate finite-dimensional transposition spaces and establish the existence and uniqueness of the corresponding finite transposition solutions. We prove the convergence of the proposed method under standard assumptions and derive convergence rates under additional Malliavin regularity assumptions. The proposed approach provides a variational approximation framework for BSDEs and may be viewed as a stochastic counterpart of the classical finite element method for partial differential equations. In addition, the conditional expectations arising in the resulting scheme can be evaluated explicitly, thereby facilitating its numerical implementation.
\end{abstract}

\bs

\no{\bf 2020 Mathematics Subject Classification}.
65C30,  
60H35, 
60H10. 

\bs

\no{\bf Key words}. Backward stochastic differential equation, finite transposition method, convergence rate.

\section{Introduction}
In this paper, for $d, n\in\dbN$, we consider the numerical approximation of the following BSDE:
\bel{bsde0}
\lt\{
\bal
&\mathrm d y(t)=f\big(t,y(t),Y(t)\big)\rd t+Y(t)\rd W(t)\,,\q t\in [0,T]\,,\\
&y(T)=y_T\,,
\eal
\rt.
\ee
where $\{W(t)\}_{t\in[0,T]}$ is an $\dbR^d$-valued Brownian motion, $f(\cdot,\cdot,\cdot)$ is the generator and $y_T\in L_{\cF_T}^2(\Omega;\dbR^n)$ is the terminal condition; see Subsection~\ref{notation-assumption} for the notation and assumptions.

The theory of BSDEs was first studied by Bismut \cite{Bismut76} in the linear setting, in connection with the maximum principle and duality methods in stochastic control. A fundamental well-posedness result for nonlinear BSDEs with Lipschitz generators was subsequently established by Pardoux and Peng \cite{Pardoux-Peng90}. Since then, this theory has developed rapidly, and BSDEs and their generalizations have become important analytical tools in stochastic control, mathematical finance, and the probabilistic representation of nonlinear PDEs; see, e.g., \cite{ElKaroui-Peng-Quenez97, Yong-Zhou99, Lv-Zhang21}.

However, explicit solutions to BSDEs are rarely available, even for linear equations with general random coefficients and terminal conditions.
Therefore, it is important to develop effective numerical methods for BSDEs and, even more challengingly, for backward stochastic partial differential equations (BSPDEs).
Indeed, the solution $\big(y(\cd),Y(\cd)\big)$ to a BSDE or BSPDE consists of two $\mathbf F$-adapted stochastic processes. In particular, the martingale integrand $Y(\cd)$  is indispensable for preserving the adaptedness of the solution, and it constitutes one of the main difficulties in both the analysis and numerical approximation of these equations. This feature fundamentally distinguishes BSDEs and BSPDEs from forward stochastic differential equations (SDEs) and stochastic PDEs.

Since the early numerical studies in the 1990s, including, for example, the four-step-scheme-based approach for a class of forward-backward stochastic differential equations (FBSDEs) \cite{Douglas-Ma-Protter96} and the random-time discretization scheme for BSDEs based on a Poisson grid \cite{Bally97}, a wide variety of numerical methods for BSDEs have been developed. Except for a few classes of methods, such as branching-diffusion methods (e.g., \cite{Henry-Tan-Touzi14}) and four-step-scheme-based methods (e.g., \cite{Douglas-Ma-Protter96, Milstein-Tretyakov06}), most numerical algorithms for BSDEs are constructed in two stages. The first stage consists of a time discretization of the BSDE, based either on one-step schemes (e.g., \cite{Hu-Nualart-Song11, ZhangJF04, Zhao-Chen-Peng06}) or on multistep schemes (e.g., \cite{Chassagneux14, Zhao-Zhang-Ju10}). The second stage concerns the numerical approximation of the resulting conditional expectations, for which several techniques have been proposed, including least-squares regression methods (e.g., \cite{Bender-Denk07, Gobet-Lemor-Warin05}), Malliavin-calculus-based methods (e.g., \cite{Bouchard-Touzi04}), quantization methods (e.g., \cite{Bally-Pages03a, Bally-Pages03b}), tree-based methods (e.g., \cite{Ma-Protter-SanMartin-Torres02}), and Wiener-chaos-decomposition methods (e.g., \cite{Briand-Labart14}).
More recently, deep-learning-based approaches have attracted considerable attention for high-dimensional BSDEs and related PDE problems. Representative approaches include the deep BSDE method, deep Galerkin method, physics-informed neural networks, and deep backward dynamic programming; see, e.g., \cite{E-Han-Jentzen17, Hure-Pham-Warin20, Raissi-Perdikaris-Karniadakis19, Sirignano-Spiliopoulos18}.
For a more comprehensive overview of numerical methods for BSDEs, we refer to the recent survey \cite{Chessari-Kawai-Shinozaki-Yamada23}.

In contrast to the extensive literature on numerical methods for BSDEs, comparatively few studies have addressed the numerical approximation of BSPDEs. We refer, for example, to \cite{Lv-Wang-Wang-Zhang22, Wang16} for spectral Galerkin schemes, to \cite{Dunst-Prohl16, Molla-Qiu23, Prohl-Wang21, Prohl-Wang22, Prohl-Wang26, Yang-Zhao20} for finite element methods, and to \cite{Li-Tang21} for a splitting method. Despite the substantial progress achieved to date, the efficient and reliable numerical approximation of general BSDEs and BSPDEs remains challenging. In particular, the computation of conditional expectations, the approximation of the martingale integrand, and the treatment of non-Markovian terminal conditions continue to pose significant difficulties. These challenges motivate the development of novel numerical approaches for BSDEs.

\ms
In the announcement article \cite{Wang-Zhang11}, the first and third authors of this paper introduced the {\em finite transposition method} for the numerical approximation of BSDE \eqref{bsde0} and laid out its foundational principles, thereby establishing a new framework for designing numerical schemes---an approach later extended to BSPDEs in \cite{Lv-Wang-Wang-Zhang22}. The current paper provides a systematic and complete development of this method. Its main contributions are threefold: first, we construct explicit finite-dimensional transposition spaces based on Wiener chaos expansions and establish the existence and uniqueness of the corresponding finite transposition solutions; second, we prove the convergence of the proposed method under mild assumptions and derive explicit convergence rates under additional regularity conditions, including a separate analysis for Markovian BSDEs; and third, we develop a coefficient-wise Monte Carlo approximation for its numerical implementation and provide corresponding error estimates.

The finite transposition method is built upon the transposition solution theory for BSDEs introduced in \cite{Lv-Zhang13}, which was subsequently extended to BSPDEs in \cite{Lv-Zhang14}. The transposition method was originally developed in the context of PDEs and their numerical approximation. Its underlying idea is to reformulate the original problem in a suitable dual space, thereby facilitating the analysis and computation of its solutions.
The finite transposition method for BSDE \eqref{bsde0} can be summarized as follows.
\begin{enumerate}

\item [(i)]  Construct suitable scalar {\em finite transposition spaces} as finite-dimensional subspaces of $L_{\dbF}^2(0,T;\dbR)$. For a scalar finite transposition space $\dbH$, define its vector-valued counterpart componentwise by
\beq
\dbH^n \deq \big\{u=(u_1,\cds,u_n)^\top\,\big|\, u_\ell\in\dbH,\, 1\leq \ell\leq n\big\} \subseteq L_{\dbF}^2(0,T;\dbR^n)\,.
\eeq
Similarly, define its matrix-valued counterpart by
\beq
\dbH^{n\times d} \deq \big\{\Phi=(\Phi_{k\ell})_{n\times d}\,\big|\, \Phi_{k\ell}\in\dbH,\, 1\leq k\leq n,\, 1\leq\ell\leq d\big\}\,. 
\eeq
These spaces play a role analogous to the finite-dimensional trial spaces used in finite element approximations of PDEs; see Section \ref{fts}.

\ss
\item [(ii)] Introduce the following variational identity:
\begin{align}\label{w629e2}
\me\big[\lan z(T),y_T\ran\big]
=\me\int_{0}^T  \Big[\big \lan z(t), f\big(\nu(t),y(t),Y(t)\big)\big\ran
+ \langle
u(t),y(t)\rangle +\langle v(t),
Y(t)\rangle \Big] \rd t\,,
\end{align}
where $u(\cd)\in\dbH_1^n$ and $v(\cd)\in\dbH_2^{n\times d}$ for suitably chosen scalar finite transposition spaces $\dbH_1$ and $\dbH_2$ (with $\dbH_1^n$ and $\dbH_2^{n\times d}$ defined similarly to the above), and
\beq
z(t)=\int_0^{\mu(t)} u(s)\rd s+\int_0^{\nu(t)} v(s)\rd W(s)\,.
\eeq
Here, $\nu(\cd)$ and $\mu(\cd)$ are piecewise constant functions defined later (in \eqref{w1e1}).
Under suitable assumptions, the variational identity \eqref{w629e2} uniquely determines a pair of adapted stochastic processes $\big(y(\cd),Y(\cd)\big)\in \dbH_1^n\times\dbH_2^{n\times d}$. More precisely, if $\{e_i^1(\cd)\}_{i=1}^{\dim(\dbH_1)}$ and $\{e_i^2(\cd)\}_{i=1}^{\dim(\dbH_2)}$ are suitably chosen bases of $\dbH_1$ and $\dbH_2$, respectively, then
\bel{w629e1}
y(\cd)=\sum_{i=1}^{\dim(\dbH_1)} a_i e_i^1(\cd)\,,
\qquad
Y(\cd)=\sum_{i=1}^{\dim(\dbH_2)} b_i e_i^2(\cd)\,,
\ee
where $a_i\in\dbR^n$ and $b_i\in\dbR^{n\times d}$ are coefficients to be determined. Such a pair is called a {\em finite transposition solution} to BSDE \eqref{bsde0}; see Section \ref{wellposed}.

\ss
\item [(iii)] Determine the coefficients $a_i$ and $b_i$ in the finite transposition solution through the variational identity \eqref{w629e2};
see Section \ref{convergence} for the corresponding convergence analysis.


\end{enumerate}

The above procedure is conceptually analogous to the classical finite element method for solving deterministic PDEs; see, e.g., \cite{Brenner-Scott08}.
Accordingly, the proposed method for BSDEs can be seen as a stochastic counterpart of the classical finite element method.
However, the term ``stochastic finite element method'' has already been used in a different context; see, e.g., \cite{Ghanem-Spanos91, Kleiber-Hien92, Matthies-Keese05, Nouy09} and the references therein. Those works are primarily concerned with the numerical solution of random PDEs.
One common technique in stochastic finite element methods for random PDEs is to employ orthogonal expansions, such as polynomial chaos and Karhunen--Lo\`eve expansions. For example, let $\{\Phi_{m}(\cd)\}_{m=1}^\infty$ be an orthogonal basis of $L^2_{\cF_T}(\Omega;\dbR)$.
Then, for each $t\in[0,T]$ and any $u(t,\cdot)\in L^2_{\cF_T}(\Omega;\dbR)$, there exists an expansion of the form
\begin{align}\label{eq0.3}
u(t,\cdot)=\sum\limits_{m=1}^\infty\alpha_{m}(t) \Phi_{m}(\cdot)\,,
\end{align}
with suitable coefficients $\alpha_{m}(t)$.
However, for the BSDE \eqref{bsde0} with $d=n=1$, the solution $\bigl(y(t,\cdot),Y(t,\cdot)\bigr)$ must satisfy the adaptedness requirement. Consequently, one needs suitable finite-dimensional subspaces of the adapted Hilbert space $L_{\dbF}^2(0,T;\dbR)$, rather than merely an orthogonal basis of $L_{\cF_T}^2(\Omega;\dbR)$. Moreover, if the expansion \eqref{eq0.3} were applied directly to the solution processes, its finite truncations would in general fail to preserve adaptedness. This is precisely why the classical stochastic finite element framework cannot be directly applied to the present problem.
Our proposed approach therefore differs substantially from the stochastic finite element methods developed in the references cited above, and for this reason we refer to it as the {\em finite transposition method}.

\ms
Among the existing approaches for solving BSDEs, least-squares regression methods are conceptually the closest to our finite transposition method, since both involve approximating random variables in suitably chosen finite-dimensional spaces. Their mechanisms for determining the expansion coefficients, however, are fundamentally different. Regression-based methods estimate these coefficients from simulated samples by solving empirical least-squares problems. By contrast, the finite transposition method exploits the orthogonality of Hermite polynomials together with the independence of Brownian increments. Once a random variable has been projected onto a prescribed finite transposition space, its conditional expectation can be evaluated explicitly by retaining only those Hermite terms involving Brownian increments available up to the conditioning time. Consequently, no additional regression problem is required for computing the conditional expectations arising in the scheme.

A further distinction lies in the construction of the finite-dimensional approximation spaces. Various families of basis functions have been employed in least-squares regression methods. For example, \cite{Gobet-Lemor-Warin05} considered indicator-function bases associated with hypercube and Voronoi partitions, local polynomial bases on Voronoi cells, and global polynomial bases; Hermite polynomial bases were used in \cite{Pelsser-Gnameho19}; different truncation strategies for multivariate Hermite expansions were investigated in \cite{Gnameho-Stadje-Pelsser24}; and indicator functions were adopted in \cite{Bender-Denk07}. In many such methods for Markovian BSDEs or FBSDEs, the regression basis functions are constructed in terms of the state process of the associated SDE. Since the exact solution to the associated SDE is generally unavailable, the state process entering the regression basis must itself be replaced by a numerical approximation. As a result, the overall error may include, in addition to time-discretization, regression, and sampling errors, a discretization error inherited from the SDE approximation.

By contrast, the finite transposition method can be applied to general BSDEs. The finite transposition spaces are constructed directly from normalized Brownian increments and the associated Hermite polynomials. Their basis functions are therefore explicitly known and need not be generated from a numerically approximated forward state process.
Unlike standard pathwise schemes, the finite transposition method does not directly compute the values of the solution along individual sample paths. Instead, the numerical problem is reduced to determining the deterministic coefficients in the finite-dimensional expansions \eqref{w629e1}. The basis functions are piecewise-constant stochastic processes whose realized values depend on the Brownian increments. Once the deterministic coefficients have been obtained, the approximate solution along any realization can be recovered by evaluating the prescribed basis functions and forming the corresponding finite sums. In this sense, the finite transposition method avoids the additional forward-state approximation error associated with state-dependent regression bases, while providing an explicit finite-dimensional representation of both components of the BSDE solution.

In \cite{Briand-Labart14}, a Wiener-chaos-decomposition method was proposed. This method exploits the structure of Wiener chaos to represent the solution in terms of orthogonal polynomial expansions and obtains numerical approximations by truncating the resulting series. Although this method and the finite transposition method both employ tools from Wiener chaos analysis, their underlying principles are radically different. The finite transposition method is based on the variational identity \eqref{w629e2}, whereas the method in \cite{Briand-Labart14} relies on Malliavin calculus and the representation $Y(t)=\big(D^1_t y(t),\cds,D^d_ty(t)\big)$.
Moreover, the finite transposition method is well-defined and convergent under the relatively mild assumption $y_T\in L^2_{\cF_T}(\Omega;\dbR^n)$ together with Assumption {\bf(B)} on the generator $f(\cd,\cd,\cd)$; see Section \ref{notation-assumption}. These assumptions are also sufficient to establish the convergence result in Theorem \ref{ftm-convergence} of Section \ref{sec-ftm-convergence}.
By contrast, the method in \cite{Briand-Labart14} requires Malliavin differentiability of the terminal condition $y_T$ and differentiability of the generator $f(\cd,\cd,\cd)$; see \cite[Section 4]{Briand-Labart14}.


\ms
The paper is organized as follows. In Section \ref{preliminary}, we briefly review the Wiener chaos expansion and construct a class of finite transposition spaces.
Section \ref{wellposed} establishes the well-posedness of finite transposition solutions.
In Section \ref{convergence}, we establish convergence results and derive convergence rates for the finite transposition method, with the detailed analysis presented in Sections \ref{ftm1}--\ref{sec-picard}.
Section \ref{Markov-BSDE} is devoted to the convergence rate for Markovian BSDEs. Since the expectations appearing in the variational identity \eqref{w629e2} enter the coefficients of the expansions \eqref{w629e1}, Section \ref{MC} develops a Monte Carlo approximation of these coefficients.
Finally, Section \ref{examples} presents numerical experiments illustrating the performance of the finite transposition method.

\section{Preliminaries}\label{preliminary}

\subsection{Notation and assumptions}\label{notation-assumption}

Denote by $\dbN$ the set of positive integers, and set $\dbN_0=\dbN\cup \{0\}$.
Let $T>0$ be a fixed terminal time. Let $(\Omega,\cF,\mathbf{F},\dbP)$ (with $\mathbf{F}=\{\cF_t\}_{t\in[0,T]}$) be a complete filtered probability space, on which an $\dbR^d$-valued standard Brownian motion $\{W(t)\}_{t\in[0,T]}=\big\{\big(W^1(t), W^2(t),\cds, W^d(t)\big)^\top\big\}_{t\in[0,T]}$ is defined such that $\mathbf{F}$ is the natural filtration generated by $W(\cd)$ (augmented by all $\dbP$-null sets).
For $1\leq p\leq \infty$ and $t\in[0,T]$, denote by $L_{\cF_t}^p(\Omega;\dbR^n)$ the Banach space consisting of all $\cF_t$-measurable $\dbR^n$-valued random variables with finite $L^p$-norm, by $\dbF$ the progressive $\si$-field with respect to $\mathbf{F}$,
and by $L^2_{\dbF}\big(\Omega;L^p(0,T;\dbR^n)\big)$ the Banach space consisting of all $\dbR^n$-valued, $\dbF$-measurable processes $X(\cdot)$ such that $\me\big[|X(\cdot)|_{L^p(0,T;\dbR^n)}^2\big]<\infty$.
For simplicity, we will write $L^2_\dbF(0,T;\dbR^n)$ instead of $L^2_\dbF\big(\Omega;L^2(0,T;\dbR^n)\big)$.

Throughout the paper, $|\cdot|$ denotes the Euclidean norm for
vector-valued quantities and the Frobenius norm for matrix-valued
quantities, according to the context. The notation
$\lan\cdot,\cdot\ran$ denotes the corresponding Euclidean or Frobenius
inner product. In particular, for matrices $A$ and $B$ of the same
size, $\lan A,B\ran=\operatorname{tr}(AB^\top)$.
All stochastic integrals with respect to $W(\cd)$ are understood componentwise.

\ss
To propose our numerical method, we impose the following assumptions:
\begin{enumerate}
\item [{\bf (P)}]
Let $I_\t:\,0=t_0<t_1<\cds<t_N=T$ be a uniform partition of $[0,T]$, where $N\in \dbN$,
$\t=T/N$, and $t_k=k \tau $ for $k=0,1,\cds,N$.

\item [{\bf (B)}]
The generator $ f: [0,T]\times\dbR^n\times\dbR^{n\times d}\longrightarrow\dbR^n$ is a deterministic function, and there exists a constant $L>0$ such that, for all $s_1,s_2,t\in[0,T]$, $y_1,y_2\in\dbR^n$, and $Y_1,Y_2\in\dbR^{n\times d}$,
\beq
 |f(s_1,y_1,Y_1)-f(s_2, y_2,Y_2)| \leq
     L\big(|s_1-s_2|^{1/2}+|y_1-y_2|+|Y_1-Y_2|\big)\,,
\eeq
and
\beq
 |f(t,0, 0)|\leq L\,.
\eeq
\end{enumerate}

%

For a given time partition $I_\t$ and each $k=0,1,\cds,N-1$, set
\bel{w9a2}
\D_{k+1}W=\big(\D_{k+1}W^1,\D_{k+1}W^2,\cds,\D_{k+1}W^d\big)^\top\,,\q
\ee
where $\D_{k+1}W^\ell=W^{\ell}(t_{k+1})-W^{\ell}(t_k)$ for $\ell=1,\cds,d$, and define
\bel{w1e1}
\bal
&\nu(t)= t_k\,,\q \pi(t)= k\,, &&\q \forall\, t\in [t_k,t_{k+1})\,,\\
&\mu(t)= t_{k+1}\,,  &&\q \forall\, t\in (t_k,t_{k+1}]\,,\\
\eal
\ee
with the conventions
\beq
\mu(0)=0\,,\q \nu(T) = T\,,\q\pi(T) = N\,.
\eeq

Throughout this paper, let $\cC$ be a generic positive constant that may depend on the problem data $T$, $n$, $d$, $y_T$ and $f(\cd,\cd,\cd)$, but is independent of the numerical parameters.

 \br{w0420r0}
{\rm (1)} In Assumption {\rm(P)}, we adopt a uniform partition for simplicity; the method can also be extended to general time partitions.

{\rm (2)} Under Assumption {\rm(B)}, there exists a unique adapted solution $\big(y(\cdot),Y(\cdot)\big)$ to BSDE~\eqref{bsde0}; see, e.g., \cite[Theorem 2.1]{ElKaroui-Peng-Quenez97}.


 \er

\subsection{Wiener chaos expansion and finite transposition spaces}\label{fts}

We next construct {\em finite transposition spaces} for the $d$-dimensional Brownian motion $W(\cd)$ by using Wiener chaos expansions.
We first recall some notation and basic results concerning Wiener chaos. Further details can be found in \cite{Nualart06}.
Define the Wiener integral operator $\dbW:L^2(0,T;\dbR^d)\longrightarrow L^2_{\cF_{T}}(\Omega;\dbR)$ by
\beq
\dbW(g)=\sum_{\ell=1}^d\int_0^T g_\ell(t)\rd W^\ell(t)\,,
\qquad g=(g_1,\cds,g_d)^\top\in L^2(0,T;\dbR^d)\,.
\eeq
For $m\in \dbN$, denote by $C_p^\infty(\dbR^m)$ the set of all smooth functions $\psi:\dbR^m\longrightarrow \dbR$ such that $\psi(\cd)$ and all of its partial derivatives have polynomial growth.
Let $\cS$ be the set of smooth random variables of the form
\begin{eqnarray}\label{eq5.1}
F=\psi \big(\dbW(g_1),\dbW(g_2),\cds,\dbW(g_m)\big)\,,
\end{eqnarray}
where $\psi(\cd)\in C_p^\infty(\dbR^m)$. Here $g_i(\cd)\in L^2(0,T;\dbR^d)$ for $i=1,\cds,m$, and we write $g_i=(g_i^1,\cds,g_i^d)^\top$.
For any $F\in\cS$ of the form \eqref{eq5.1} and each $\ell=1,2,\cds,d$, define its Malliavin derivative with respect to the $\ell$-th component of $W(\cd)$ by
\beq
D_t^\ell F
\deq
\sum_{i=1}^m
\frac{\partial\psi}{\partial x_i}
\big(\dbW(g_1),\dbW(g_2),\cds,\dbW(g_m)\big)
g_i^\ell(t)\,,\q t\in[0,T]\,.
\eeq
The first Malliavin derivative of $F$ is then defined by
\beq
D^{(1)}F \deq \big(D^1F,D^2F,\cds,D^dF\big)^\top.
\eeq
More generally, for $k\geq2$, the $k$-th Malliavin derivative $D^{(k)}F$ is defined iteratively. For $\theta_1,\cds,\theta_k\in[0,T]$ and $\ell_1,\cds,\ell_k\in\{1,\cds,d\}$, we write $D_{\theta_1,\cds,\theta_k}^{\ell_1,\cds,\ell_k}F \deq D_{\theta_k}^{\ell_k}\cds D_{\theta_1}^{\ell_1}F$.
For any $j\in\dbN$, we use the convention
\beq
\big|D^{(j)}F\big|_{L^2((0,T)^j)}^2
\deq \sum_{\ell_1,\cds,\ell_j=1}^d
\int_{(0,T)^j} \big| D_{\theta_1,\cds,\theta_j}^{\ell_1,\cds,\ell_j}F \big|^2 \rd\theta_1\cds\rd\theta_j\,.
\eeq
For $k\in\dbN$ and $p\in[1,\infty)$, set
\beq
\|F\|_{k,p}
\deq
\Big( \me\big[|F|^p\big]
+ \sum_{j=1}^k \me\big[ \big|D^{(j)}F\big|_{L^2((0,T)^j)}^p\big]
\Big)^{1/p}.
\eeq
The completion of $\cS$ under the norm $\|\cdot\|_{k,p}$ is denoted by $\dbD^{k,p}$.

For $\dbR^r$-valued random variables, where $r\in\dbN$, we use the componentwise definition and denote the corresponding space by $\dbD^{k,p}(\dbR^r)$. In this case, we write $|D^{(j)}F|_{L^2((0,T)^j;\dbR^r)}^2$ in place of $|D^{(j)}F|_{L^2((0,T)^j)}^2$.
We further set
\beq
\dbD^{k,\infty}(\dbR^r)
\deq \bigcap_{1\leq p<\infty}\dbD^{k,p}(\dbR^r).
\eeq

For $x\in\dbR$, let
\beq
H_m(x)=
\left\{\!\!\!
\begin{array}{ll}
\ds \frac{(-1)^m} {m!}e^{\frac{x^2} {2}}\frac{\rd^m}{\rd x^m}(e^{-{\frac{x^2} {2}}})\,,\qq & m>0\,,\\
\ns\ds 1\,, & m=0
\end{array}
\right.
\eeq
be the $m$-th Hermite polynomial, and set $H_{-1}(x)=0$.
It is easy to see that
$H_1(x)=x\,, H_2(x)=\frac 1 2(x^2-1)$, and $H_3(x)=\frac 1 6(x^3-3x)$.

The following lemma gives a fundamental relation between Hermite polynomials and jointly Gaussian random variables.

\begin{lemma}[{\cite[Lemma 1.1.1]{Nualart06}}]\label{w512l1}
Let $(\xi,\eta)$ be a jointly Gaussian random vector such that
$\me[\xi]=\me[\eta]=0$, and $\me[\xi^2]=\me[\eta^2]=1$.
Then, for any $m,r\in\dbN_0$,
\beq
\me\big[H_m(\xi)H_r(\eta)\big]
=
\lt\{\!\!\!
\begin{array}{ll}
0\,,&\q  m\neq r\,,\\
\ns\ds\frac 1{m!}\big[\me(\xi\eta)\big]^m\,,  &\q m=r\,.
\end{array}
\rt.
\eeq
\end{lemma}

Let $\Lambda_d$ be the set of all sequences $\a=(\a_j^\ell)_{ j\geq 1 , 1\leq\ell\leq d}$ with $\a_j^\ell\in\dbN_0$ such that only finitely many components are nonzero.
For each $k\in\dbN$, set $\Lambda_d(k)\deq\dbN_0^{k\times d}$.
For consistency, we also set $\Lambda_d(0)=\{\emptyset\}$, where $\emptyset$ denotes the unique empty multi-index.
For any $\alpha\in \Lambda_d \, (\mbox{resp. }\a\in\Lambda_d(k))$, define
\beq
\alpha!= \prod\limits_{j=1}^\infty\prod\limits_{\ell=1}^d \alpha_j^\ell! \, \q\Big(\mbox{resp.}\prod\limits_{j=1}^k\prod\limits_{\ell=1}^d \alpha_j^\ell!\Big)\,,\quad
\mbox{and}
\quad |\alpha|= \sum\limits_{j=1}^\infty  \sum\limits_{\ell=1}^d \alpha_j^\ell \, \q\Big(\mbox{resp.}\sum\limits_{j=1}^k  \sum\limits_{\ell=1}^d \alpha_j^\ell\Big)\,.
\eeq
For $k=0$, the unique empty multi-index has order zero and factorial one.

Under Assumption {\rm\bf(P)}, we define
$g_j(\cd)= \frac{\chi_{[t_{j-1},t_j)}(\cd)}{\sqrt{\tau }}$ for $j=1,2,\cds,N$.
Let $\{{\bf e}_\ell\}_{\ell=1}^d$ be the canonical basis of $\dbR^d$ and put
$g_{j,\ell}(t)=g_j(t){\bf e}_\ell\in L^2(0,T;\dbR^d).$
Then
\bel{w907e1}
\dbW(g_{j,\ell})=\frac{\D_jW^\ell}{\sqrt\tau}\,, \q 1\leq j\leq N\,, 1\leq\ell\leq d\,.
\ee
For $m\in\dbN_0$ and $k=0,1,\cds,N$, define the {\em Wiener chaos of order $m$} up to the time node $t_k$ by
\beq
\cH_m(k)=\Span\Big\{  \prod_{j=1}^k \prod_{\ell=1}^d H_{\a_j^\ell}\big(\dbW(g_{j,\ell})\big) \;\big|\; \a\in \Lambda_d(k),\,|\a|=m\Big\}\,,
\eeq
where, for $k=0$, the empty product is understood as $1$.

To construct the space containing the finite transposition solutions, set $\mathcal{G}_0=\cF_0$ and, for $k=1,2,\cds,N$, define
\bel{w9a1}
\mathcal{G}_{t_k}= \si\big(\big\{\D_jW^{\ell}\;\big|\;1\leq j\leq k,\,1\leq \ell\leq d\big\}\cup \cN\big)\,,
\ee
where $\D_kW$ is defined in \eqref{w9a2} and $\cN$ is the collection of $\dbP$-null sets. For $r\in\dbN$, denote by $L^2_{\mathcal{G}_{t_k}}(\Omega;\dbR^r)$ the $\dbR^r$-valued, $\mathcal{G}_{t_k}$-measurable square-integrable random variables.
The Wiener chaos decomposition yields the following theorem.

\begin{theorem}[{\cite[Theorem 1.1.1]{Nualart06}}]\label{w9l2}
For any $k=0,1,\cds,N$, it holds that
\beq
L^2_{\cG_{t_k}}(\Omega;\dbR)=\bigoplus\limits_{m=0}^\infty \cH_m(k)\,.
\eeq
\et

Motivated by Theorem \ref{w9l2}, we define
\bel{w9d1}
\bal
\cH^M(k)\deq \bigoplus\limits_{m=0}^M \cH_m(k)\,,
\eal
\ee
which is a finite-dimensional subspace of $L^2_{\cG_{t_k}}(\Omega;\dbR)$. We then introduce the {\em finite transposition space}
\bel{w9d2}
\begin{array}{c}
\dbH_{N,M}\deq \bigoplus\limits_{k=0}^{N-1}\big\{\chi_{[t_k,t_{k+1})}(\cd)\xi\,|\,\xi\in\cH^M(k)\big\}
\subseteq L^2_\dbF(0,T;\dbR)\,,
\end{array}
\ee
equipped with the inherited norm $|\cd|_{\dbH_{N,M}}\deq|\cd|_{L^2_\dbF(0,T;\dbR)}$.
For the BSDE variables in \eqref{bsde0}, we use the componentwise extensions
\beq
&\dbH_{N,M}^{n}
\deq\big\{(u_1,\cds,u_n)^\top\,\big|\,u_r\in\dbH_{N,M},\,1\leq r\leq n\big\}\,,\\
\mbox{and }\q&
\dbH_{N,M}^{n\times d}
\deq\big\{(v_{r\ell})_{n\times d}\,\big|\,v_{r\ell}\in\dbH_{N,M},\,1\leq r\leq n,\, 1\leq\ell\leq d\big\}\,.
\eeq
These spaces are endowed with the natural product norms, which coincide with the Euclidean and Frobenius norms.
For any $t\in[0,T)$, we further use the notation
$\dbH^n_{N,M}(\nu(t),T) \deq \big\{\chi_{[\nu(t),T)}(\cd)u\,\big|\,u\in\dbH^n_{N,M}\big\}$,
equipped with the norm $|\cdot|_{\dbH^n_{N,M}(\nu(t),T)} \deq |\cdot|_{L^2_\dbF(0,T;\dbR^n)}$.
The space $\dbH^{n\times d}_{N,M}(\nu(t),T)$ and its norm are defined analogously.

\br{w531r1}
Similar to Theorem \ref{w9l2}, the following orthogonal decomposition theorem for $L^2_{\cF_T}(\Omega;\dbR)$ also holds:
\beq
L^2_{\cF_T}(\Omega;\dbR)=\bigoplus\limits_{m=0}^\infty \cH_m\,,
\eeq
where $\cH_m$ is the Wiener chaos of order $m$.
More precisely, let
$\{\widehat g_j(\cd)\}_{j=1}^\infty$
be an orthonormal basis of $L^2(0,T;\dbR)$, and let
$\{{\bf e}_\ell\}_{\ell=1}^d$ be the canonical basis of $\dbR^d$.
Set
\beq
\widehat g_{j,\ell}(\cd)
=\widehat g_j(\cd){\bf e}_\ell \in L^2(0,T;\dbR^d)\,,
\qquad j\in\dbN,\quad \ell=1,2,\cds,d\,.
\eeq
Then $\{\widehat g_{j,\ell}\}_{j\geq1,\,1\leq\ell\leq d}$ is an orthonormal basis of $L^2(0,T;\dbR^d)$, and
\beq
\cH_m
=\Span\Big\{\prod_{j=1}^\infty\prod_{\ell=1}^dH_{\a_j^\ell}\big(\dbW(\widehat g_{j,\ell})\big)\;\Big|\;\a\in\Lambda_d,\ |\a|=m\Big\}\,.
\eeq
Consequently, for any
$y_T\in L^2_{\cF_T}(\Omega;\dbR^n)$, there exists a unique sequence of deterministic kernels $h_m=\big\{h_m^{\ell_1,\cds,\ell_m}\big\}_{1\leq\ell_1,\cds,\ell_m\leq d}\in L^2_{S,d}(T^m;\dbR^n)$
such that $y_T$ admits the following Wiener chaos expansion
\bel{ex-y}
\bal
y_T
&=: \me[y_T]
+\sum_{m=1}^{+\infty}J_m(h_m)\,,
\eal
\ee
where $J_m(h_m)$ denotes the sum of
$m$-fold iterated It\^o integrals over all $\ell_1,\cds,\ell_m\in\{1,\cds,d\}$:
\beq
J_m(h_m)=\sum_{\ell_1,\cds,\ell_m=1}^d
\!\int_0^T\!\!\int_0^{s_m}\!\!\!\cds\! \int_0^{s_2} h_m^{\ell_1,\cds,\ell_m}(s_1,\cds\!,s_m) \rd W^{\ell_1}(s_1) \cds \rd W^{\ell_m}(s_m)\,.
\eeq
Here $L^2_{S,d}(T^m;\dbR^n)$ is the space of families $h=\{h^{\ell_1,\cds,\ell_m}\}_{1\leq\ell_1,\cds,\ell_m\leq d}$ such that $h^{\ell_1,\cds,\ell_m}\in L^2\big((0,T)^m;\dbR^n\big)$ and
\beq
h^{\ell_1,\cds,\ell_m}(t_1,\cds,t_m)=h^{\ell_{\varrho(1)},\cds,\ell_{\varrho(m)}}(t_{\varrho(1)},\cds,t_{\varrho(m)})
\eeq
for any $\ell_1,\cds,\ell_m\in\{1,\cds,d\}$, any permutation $\varrho$ of $\{1,2,\cds,m\}$, and a.e. $(t_1,\cds,t_m)\in(0,T)^m$. We equip this space with the norm
\beq
|h|_{L^2_{S,d}(T^m;\dbR^n)}^2
\deq \sum_{\ell_1,\cds,\ell_m=1}^d \big|h^{\ell_1,\cds,\ell_m}\big|^2_{L^2((0,T)^m;\dbR^n)}\,.
\eeq

For any $M,N\in\dbN$, define the componentwise $L^2$-orthogonal projections
\bel{w417e1}
\G_M:
L^2_{\cF_T}(\Omega;\dbR^n)
\longrightarrow
\bigoplus_{m=0}^M\big(\cH_m\big)^n\,,
\qquad
\G_M^N:
L^2_{\cF_T}(\Omega;\dbR^n)
\longrightarrow
\big[\cH^M(N)\big]^n\,.
\ee
For any $r\in\dbN$, we use the same notation $\G_M$ for the componentwise $L^2$-orthogonal projection onto $\bigoplus_{m=0}^M(\cH_m)^r$.
For any $\xi\in L^2_{\cG_T}(\Omega;\dbR^n)$, it holds that
\beq
\G_M\xi=\G_M^N\xi\,.
\eeq
\er

The following lemma presents an explicit representation of $\G_M^Ny_T$ for $y_T\in L^2_{\cF_T}(\Omega;\dbR^n)$.

\bl{w524l1}
Suppose that $y_T\in L^2_{\cF_T}(\Omega;\dbR^n)$ has the Wiener chaos expansion \eqref{ex-y}. Then, for any $N,M\in\dbN$,
\begin{align}
\G_M^N y_T
&= \me[y_T]+\sum_{m=1}^M J_m\big(\bar h_m\big)\,, \label{w518e3}
\end{align}
where
$\bar h_m=\big\{ \bar h_m^{\ell_1,\cds,\ell_m}\big\}_{1\leq\ell_1,\cds,\ell_m\leq d}\in L^2_{S,d}(T^m;\dbR^n)$ is given by
\bel{w524e7}
\bar h_m^{\ell_1,\cds,\ell_m}(s_1,\cds,s_m)
= \frac{1}{|I_{\Bj_m}|} \int_{I_{\Bj_m}} h_m^{\ell_1,\cds,\ell_m}(t_1,\cds,t_m) \rd t_1\cds\rd t_m
\ee
for a.e. $(s_1,\cds,s_m)\in I_{\Bj_m}$ and all $\ell_1,\cds,\ell_m\in\{1,\cds,d\}$, with
\bel{w524e3}
\bal
&\Bj_m\deq(j_1,j_2,\cds,j_m)\,,\q
0\leq j_1,j_2,\cds,j_m\leq N-1\,,\\
&I_{\Bj_m}
\deq
[j_1\t,(j_1+1)\t)\times\cds\times
[j_m\t,(j_m+1)\t)\,.
\eal
\ee
\el

\begin{proof}
We divide the proof into three steps.

\ss\no
{\bf Step 1.}
We first show that $J_m\big(\bar h_m\big)\in\big[\cH_m(N)\big]^n$ for $m=1,2,\cds,M$.
Relying on \eqref{w907e1},
$\big\{\dbW(g_{j,\ell})\,\big|\,1\leq j\leq N,\ 1\leq\ell\leq d \big\}$
is a family of independent standard normal random variables.

Since each $\bar h_m^{\ell_1,\cds,\ell_m}$ is constant on every rectangle $I_{\Bj_m}$, $J_m(\bar h_m)$ is a finite sum of iterated It\^o integrals whose deterministic integrands are constant on products of the partition intervals.
By the relation between multiple Wiener--It\^o integrals of step functions and Hermite polynomials (see, e.g., \cite[Proposition 1.1.4]{Nualart06}), each component of such an integral is a finite linear combination of random variables of the form
\beq
\prod_{j=1}^N\prod_{\ell=1}^d
H_{\a_j^\ell}\big(\dbW(g_{j,\ell})\big)\,,
\qquad
\a_j^\ell\in\dbN_0\,,
\qquad
\sum_{j=1}^N\sum_{\ell=1}^d \a_j^\ell=m\,.
\eeq
By the definition of the $m$-th Wiener chaos up to the time node $t_N=T$, 
$\prod_{j=1}^N\prod_{\ell=1}^d H_{\a_j^\ell}\big(\dbW(g_{j,\ell})\big) \in\cH_m(N)$.
Consequently,
\beq
J_m\big(\bar h_m\big)
\in
\big[\cH_m(N)\big]^n
\subseteq
\big[\cH^M(N)\big]^n\,.
\eeq

\ss\no
{\bf Step 2.}
We claim that, for any $\xi\in[\cH^M(N)]^n$, there exist deterministic kernels $\{\widehat h_m\}_{m=1}^M\subseteq L^2_{S,d}(T^m;\dbR^n)$, whose components are constant on every rectangle $I_{\Bj_m}$, such that
\bel{w524e5}
\xi=\me[\xi]
+\sum_{m=1}^M J_m\big(\widehat h_m\big)\,.
\ee

Indeed, by the construction of $\cH^M(N)$, it suffices to consider, componentwise, basis elements of the form
$\prod_{j=1}^N\prod_{\ell=1}^d H_{\alpha_j^\ell}\big(\dbW(g_{j,\ell})\big)$ and $\sum_{j=1}^N\sum_{\ell=1}^d\alpha_j^\ell\leq M$ for $\alpha_j^\ell\in\dbN_0$.
By the relation between Hermite polynomials and multiple Wiener--It\^o integrals, each such basis element can be represented by a multiple Wiener--It\^o integral whose deterministic kernel is constant on the rectangles generated by the partition $I_\tau$.
After symmetrization with respect to the simultaneous permutations of the time variables and the corresponding Brownian component indices, the resulting kernel belongs to $L^2_{S,d}(T^m;\dbR)$.
Applying this argument to every component of $\xi$, we obtain $\widehat h_m\in L^2_{S,d}(T^m;\dbR^n)$ and hence \eqref{w524e5}.

\ss\no
{\bf Step 3.}
We now identify the $L^2$-orthogonal projection.
Since $\G_M^N$ is the $L^2$-orthogonal projection onto $[\cH^M(N)]^n$, it suffices to determine the unique minimizer of
\beq
\xi \longmapsto \me\big[\big|\G_My_T-\xi\big|^2\big]\,,
\qquad
\xi\in[\cH^M(N)]^n\,.
\eeq
For any $\xi\in[\cH^M(N)]^n$ with the representation \eqref{w524e5}, we have
\bel{w709e1}
\G_My_T-\xi
=\me[y_T-\xi]
+\sum_{m=1}^{M}J_m\big(h_m-\bar h_m\big)
+\sum_{m=1}^{M}J_m\big(\bar h_m-\widehat h_m\big)\,.
\ee

By \eqref{w524e7}, for every $m=1,2,\cds,M$ and every $\ell_1,\cds,\ell_m\in\{1,\cds,d\}$, $\bar h_m^{\ell_1,\cds,\ell_m}$ is the $L^2((0,T)^m;\dbR^n)$-orthogonal projection of $h_m^{\ell_1,\cds,\ell_m}$ onto the space of functions that are constant on the rectangles $I_{\Bj_m}$.
Moreover, since $h_m$ is symmetric with respect to simultaneous permutations of the time variables and the corresponding Brownian component indices, the kernel $\bar h_m$ defined by \eqref{w524e7} also belongs to $L^2_{S,d}(T^m;\dbR^n)$.
Therefore,
$\sum_{\ell_1,\cds,\ell_m=1}^d
\int_{(0,T)^m}\lan h_m^{\ell_1,\cds,\ell_m}
-\bar h_m^{\ell_1,\cds,\ell_m},
\bar h_m^{\ell_1,\cds,\ell_m}
-\widehat h_m^{\ell_1,\cds,\ell_m}
\ran
\rd s_1\cds\rd s_m
=0$.
The multiple It\^o isometry then yields
$\me\big[ \big\lan J_m(h_m-\bar h_m),
J_m(\bar h_m-\widehat h_m) \big\ran \big]
=0$.
Moreover, Wiener chaoses of different orders are mutually orthogonal, and all multiple It\^o integrals of positive order are centered. Hence,
\beq
\me\big[\big|\G_My_T-\xi\big|^2\big]
&= \me\Big[ \Big| \sum_{m=1}^M J_m(h_m-\bar h_m) \Big|^2 \Big]
\!+\!\me\Big[ \Big| \me[y_T-\xi]
\!+\!\sum_{m=1}^M J_m(\bar h_m-\widehat h_m) \Big|^2 \Big]\\
&\geq \me\Big[ \Big| \sum_{m=1}^M J_m(h_m-\bar h_m) \Big|^2 \Big]\,.
\eeq
The equality holds if and only if
\beq
\me[\xi]=\me[y_T]\,, \q\mbox{and}\q
\widehat h_m=\bar h_m\,,
\quad m=1,2,\cds,M\,.
\eeq
Therefore, the unique minimizer is $\me[y_T]+\sum_{m=1}^{M}J_m\big(\bar h_m\big)$.

Finally, $y_T-\G_My_T$ is orthogonal to $\bigoplus_{m=0}^M(\cH_m)^n$ and therefore to its subspace $[\cH^M(N)]^n$.
Hence, the same minimizer is the $L^2$-orthogonal projection of $y_T$ onto $[\cH^M(N)]^n$,
which proves \eqref{w518e3}.
\end{proof}

\br{w9r1}
{\rm (1)}
For each $k=0,1,\cds,N$, let $\big\{\a^{k,i}\big\}_{i=1}^{\dim(\cH^M(k))}$ be a fixed graded enumeration of the set $\big\{\a\in\Lambda_d(k)\,\big|\,|\a|\leq M\big\}$. Then
\beq
\cE_0(k)=\Big\{h_{k,i}\deq\sqrt{\a^{k,i}!}
\prod_{j=1}^k\prod_{\ell=1}^d
H_{\a_j^{k,i,\ell}}\big(\dbW(g_{j,\ell})\big)\,\Big|\,1\leq i\leq\dim(\cH^M(k))\Big\}
\eeq
is an orthonormal basis of $\cH^M(k)$.
Thus, for each $0\leq m\leq M$, the first $\dim(\cH^m(k))$ elements correspond precisely to the multi-indices with total degree not exceeding $m$.
Accordingly, an orthonormal basis of $\dbH_{N,M}$ is given by
\beq
\cE = \Big\{e_{k,i}\deq \chi_{[t_k,t_{k+1})}(\cd)\frac{1}{\sqrt{\tau }}h_{k,i}\,\Big|\, 1\leq i\leq \dim\big(\cH^M(k)\big), 0\leq k\leq N-1\Big\}\,.
\eeq

{\rm (2)}
We next determine the dimensions of the finite-dimensional spaces introduced above.
For $k=0$, we have $\cH^M(0)=\dbR$, and hence
$
\dim\big(\cH^M(0)\big)=1=C_M^M.
$
Here
\beq
C_p^q=\frac{p!}{q!(p-q)!}
\eeq
denotes the binomial coefficient for $p,q\in\dbN_0$ with $p\geq q$.

For $k\geq1$, the space $\cH_m(k)$ is spanned by the Hermite polynomials indexed by $\a\in\Lambda_d(k)$ satisfying $|\a|=m$.
Thus, $\dim\big(\cH_m(k)\big)$ equals the number of nonnegative integer solutions of $\sum_{j=1}^k\sum_{\ell=1}^d\a_j^\ell=m$.
Since there are $kd$ indices $\a_j^\ell$, the standard stars-and-bars formula gives $\dim\big(\cH_m(k)\big) = C_{kd+m-1}^{m}$.
Consequently, for $k\geq1$, $\dim\big(\cH^M(k)\big)=\sum_{m=0}^{M}C_{kd+m-1}^{m} =C_{kd+M}^{M}$.
Together with $\dim\big(\cH^M(0)\big)=1=C_M^M$, we obtain
\beq
\dim\big(\cH^M(k)\big)
= C_{kd+M}^{M}\,, \qquad k=0,1,\cds,N\,.
\eeq

By the definition of $\dbH_{N,M}$ in \eqref{w9d2} and the disjointness of the time intervals $[t_k,t_{k+1})$, we have
\beq
\dim\big(\dbH_{N,M}\big)
=\sum_{k=0}^{N-1}\dim\big(\cH^M(k)\big)
=\sum_{k=0}^{N-1}C_{kd+M}^{M}\,.
\eeq
Accordingly, for the componentwise extensions used for the BSDE variables,
\beq
\dim\big(\dbH_{N,M}^{n}\big)
= n\sum_{k=0}^{N-1}C_{kd+M}^{M}\,,
\q \mbox{and}\q
\dim\big(\dbH_{N,M}^{n\times d}\big)
= nd\sum_{k=0}^{N-1}C_{kd+M}^{M}\,.
\eeq

For fixed $N$ and $M$, the dimension of the finite transposition space grows polynomially with respect to the dimension $d$ of the Brownian motion. More precisely,
\beq
\sum_{k=0}^{N-1}C_{kd+M}^{M}
= \frac{d^M}{M!}\sum_{k=1}^{N-1}k^M +\mathcal O(d^{M-1})
= \mathcal O(d^M)\,.
\eeq
Moreover, since
$\sum_{k=1}^{N-1}k^M=\mathcal O(N^{M+1})$, for fixed $M$ we have
\beq \dim\big(\dbH_{N,M}\big)
=\mathcal O\big(d^M N^{M+1}\big)\,.
\eeq
\er

\br{w408r2}
The computation or approximation of conditional expectations is one of the central difficulties in the numerical solution of BSDEs; see Theorem \ref{ftm-rate} below and, e.g., \cite[Section 5]{ZhangJF04}, \cite[Section 6]{Bouchard-Touzi04}, and \cite[Section 2]{Bender-Denk07}.
For any $\xi\in L^2_{\cF_T}(\Omega;\dbR^n)$, the conditional expectation $\me\big(\xi\,\big|\,\cF_{t_k}\big)$ may be approximated by $\me\big(\G_M^N\xi\,\big|\,\cF_{t_k}\big)$, where $\G_M^N$ is the componentwise $L^2$-orthogonal projection onto the finite-dimensional Wiener chaos space $[\cH^M(N)]^n$ defined in \eqref{w417e1}.

More precisely, write
\beq
\G_M^N\xi
=\sum_{m=0}^M \sum_{|\a|=m, \a\in\Lambda_d(N)} d^\a\sqrt{\a!} \prod_{j=1}^N\prod_{\ell=1}^d H_{\a_j^\ell}\big(\dbW(g_{j,\ell})\big)\,,
\eeq
where $d^\a\in\dbR^n$.
Since the Brownian increments $\{\D_iW^\ell\,|\,1\leq i\leq N,\ 1\leq\ell\leq d\}$ are mutually independent, and the increments after $t_k$ are independent of $\cF_{t_k}$, we have
\begin{align}
\me\big(\G_M^N\xi\,\big|\,\cF_{t_k}\big)
&= \sum_{m=0}^M \sum_{\substack{|\a|=m\\ \a\in\Lambda_d(N)}} d^\a\sqrt{\a!} \prod_{j=1}^k\prod_{\ell=1}^d H_{\a_j^\ell}\big(\dbW(g_{j,\ell})\big)
\prod_{j=k+1}^N\prod_{\ell=1}^d \me\Big[ H_{\a_j^\ell}\big(\dbW(g_{j,\ell})\big)\Big] \notag\\
&= \sum_{m=0}^M \sum_{\substack{ |\a|=m,\ \a\in\Lambda_d(N)\\
\a_j^\ell=0,\ k+1\leq j\leq N,\ 1\leq\ell\leq d }}
d^\a\sqrt{\a!} \prod_{j=1}^k\prod_{\ell=1}^d H_{\a_j^\ell}\big(\dbW(g_{j,\ell})\big)\,.
\label{w408e5}
\end{align}
Equivalently, \eqref{w408e5} is obtained from the Wiener chaos expansion of $\G_M^N\xi$ by retaining only those terms involving the Brownian increments up to the time node $t_k$.

Finally, we have
\beq
\lim_{M\to\infty}\lim_{N\to\infty}
\me\big[\big|\xi-\G_M^N\xi\big|^2\big]=0\,.
\eeq
By the $L^2$-contraction property of conditional expectation,
\beq
\me\big[\big|\me\big(\xi\,\big|\,\cF_{t_k}\big)
-\me\big(\G_M^N\xi\,\big|\,\cF_{t_k}\big)\big|^2\big]
\leq \me\big[\big|\xi-\G_M^N\xi\big|^2\big]\,.
\eeq
Thus, $\me\big(\G_M^N\xi\,\big|\,\cF_{t_k}\big)$ given by \eqref{w408e5} provides an $L^2$-approximation of $\me\big(\xi\,\big|\,\cF_{t_k}\big)$ in the above iterated-limit sense.
\er

\section{Well-posedness of finite transposition solutions}\label{wellposed}

In this section, we begin by introducing the definition of the finite transposition solution and proposing the finite transposition method for BSDEs. We then establish the existence and uniqueness of finite transposition solutions.
\begin{definition}\label{defi-of-ftm}
For any $N,M\in\dbN$, a pair $\big(y(\cdot), Y(\cdot)\big)\in \dbH_{N,M}^n \times \dbH_{N,M-1}^{n\times d}$ is called a finite transposition solution to BSDE \eqref{bsde0} with terminal approximation $y_T^\pi\in [\cH^M(N)]^n$ if, for every $u(\cd)\in \dbH_{N,M}^{n},\,v(\cd)\in \dbH_{N,M-1}^{n\times d}$, the following variational identity holds
\bel{w9e3}
\bal
\me\big[\lan z(T),y_T^\pi\ran\big]
=\me\int_{0}^T  \Big[ \big\langle z(t), f\big(\nu(t),y(t),Y(t)\big)\big\rangle
+ \langle u(t),y(t)\rangle +\langle v(t), Y(t)\rangle \Big] \rd t\,,
\eal
\ee
where the test process $z(\cd)$ solves the following SDE:
\bel{w9e4}
z(s)=\int_{0}^{\mu(s)} u(t)\rd t+\int_{0}^{\nu(s)} v(t) \rd W(t)\,,
\ee
and $\mu(\cd)$, $\nu(\cd)$ are given by \eqref{w1e1}.

The numerical scheme defined by the variational identity \eqref{w9e3} is called the finite transposition method for BSDE \eqref{bsde0}.
\end{definition}


For a general BSDE, one may take $y^\pi_T=\G_M^Ny_T$; see Theorem \ref{ftm-rate} in Section \ref{ftm1}. For a Markovian BSDE, one may instead take $y^\pi_T=\G_M\varphi\big(x^\pi(T)\big)$; see Remark \ref{w531r2} in Section \ref{Markov-BSDE}.

We recall the following standard result for later use.
\bl{w10l1}
Suppose that Assumption {\rm\bf(P)} holds. For any $N,M\in\dbN$, let $u(\cd)\in \dbH_{N,M}^{n},\,v(\cd)\in \dbH_{N,M-1}^{n\times d}$.
Then the solution $z(\cd)$ to SDE \eqref{w9e4} satisfies $z(\cd)\in \dbH_{N,M}^n$ and $z(t)\in \big[\cH^M\big(\pi(t)\big)\big]^n$ for any $t\in [0,T]$.
Furthermore, it holds that
\bel{w10e1}
|z(\cd)|_{L^2_{\dbF}(\Omega;L^\infty(0,T;\dbR^n))}
\leq \cC\big[|u(\cd)|_{\dbH_{N,M}^{n}}+|v(\cd)|_{\dbH_{N,M-1}^{n\times d}}\big]\,,
\ee
where $\cC$ is a constant depending only on $T$.
\el

The following theorem establishes the existence and uniqueness of finite transposition solutions to BSDEs.

\begin{theorem}\label{w9t1}
Under Assumptions {\rm\bf(P)} and {\rm\bf(B)}, there exists $\t_0>0$ such that, for any $N,M\in\dbN$ satisfying $\tau=T/N\leq \t_0$, BSDE \eqref{bsde0} with $y(T)=y^\pi_T\in [\cH^M(N)]^n$ admits a unique finite transposition solution $\big(y(\cd),Y(\cd)\big)\in \dbH_{N,M}^{n}\times\dbH_{N,M-1}^{n\times d}$.
\end{theorem}

\begin{proof}
We follow the strategy used in the proofs of \cite[Theorems 3.1 and 4.1]{Lv-Zhang13} for linear and semilinear BSDEs, respectively. We first consider the linear case and then proceed to the semilinear case.

\ss\no
{\bf Step 1. Linear BSDEs.}
We first consider the case where the generator is independent of $y$ and $Y$, that is, $f\big(t,y(t),Y(t)\big)=f(t)$ for $t\in[0,T]$, where $f(\cd)$ is allowed to be stochastic. More precisely, we assume that $f(\cd)\in L_{\dbF}^2\big(\Omega;L^1(0,T;\dbR^n)\big)$ and $f(\cd)$ is piecewise constant with respect to the partition $I_\tau$, i.e., $f(t)=f\big(\nu(t)\big)$ a.e.
We prove that the corresponding BSDE admits a unique finite transposition solution $ \big(y(\cd),Y(\cd)\big) \in\dbH_{N,M}^{n}\times\dbH_{N,M-1}^{n\times d}$ satisfying the following strengthened variational identity:
\bel{w10e2}
\bal
\me\big[\lan z^s(T),y_T^\pi\ran\big] -\me\big[\big\lan \eta, y\big(\nu(s)\big)\big\ran\big]
=\me\int_{\nu(s)}^T  \Big[ \big\langle z^s(t), f\big(\nu(t)\big)\big\rangle
+ \langle u(t),y(t)\rangle +\langle v(t), Y(t)\rangle \Big] \rd t\,,
\eal
\ee
for any $s\in [0,T)$, where $z^s(\cd)$ solves the following test SDE:
\bel{w10e4}
z^s(t)=\eta+\int_{\nu(s)}^{\mu(t)}u(\th)\rd\th+\int_{\nu(s)}^{\nu(t)}v(\th)\rd W(\th)\,,\q 0\leq s\leq t\leq T\,,
\ee
with $\eta\in [\cH^M\big(\pi(s)\big)]^n,\,u(\cd)\in \dbH_{N,M}^{n},\,v(\cd)\in \dbH_{N,M-1}^{n\times d}$.

The variational identity \eqref{w10e2} is a strengthened version of \eqref{w9e3}. In particular, taking $s=0$ and $\eta=0$ in \eqref{w10e2} recovers \eqref{w9e3}.

We establish the assertion in four steps.

\ss\no
{\bf Step 1(a).} For any $t\in [0,T]$, set $\mathfrak{S}_t=\dbH_{N,M}^{n}\big(\nu(t),T\big)\times\dbH_{N,M-1}^{n\times d}\big(\nu(t),T\big)\times \big[\cH^M\big(\pi(t)\big)\big]^n$, and define a linear functional $\mathcal{V}_t$ on $\mathfrak{S}_t$ by
\beq
\mathcal{V}_t \big(\chi_{[\nu(t),T]}u,\chi_{[\nu(t),T]}v,\eta\big)=\me\big[\lan z^t(T),y_T^\pi \ran\big]-\me\int_{\nu(t)}^T\big\lan z^t(s),f\big(\nu(s)\big) \big\ran \rd s\,,
\eeq
where $z^t(\cd)$ solves \eqref{w10e4}.
Applying Lemma \ref{w10l1} to $z^t(\cd)-\eta$ and then using the triangle inequality, we have
$|z^t(\cd)|_{L^2_{\dbF}(\Omega;L^\infty(\nu(t),T;\dbR^n))}
\leq \cC |(\chi_{[\nu(t),T]}u,\chi_{[\nu(t),T]}v,\eta)|_{\mathfrak{S}_t}$.
Then, the Cauchy--Schwarz inequality yields
\begin{align*}
&|\mathcal{V}_t(\chi_{[\nu(t),T]}u,\chi_{[\nu(t),T]}v,\eta)|\\
&\leq  |z^t(T)|_{L^2_{\cF_T}\!(\Omega;\dbR^n\!)} |y_T^\pi|_{L^2_{\cF_T}\!(\Omega;\dbR^n\!)}
        \!+\!|z^t(\cd)|_{L^2_{\dbF}(\Omega;L^{\!\infty\!}(\nu(t),T;\dbR^n))}\big|f\!\big(\nu(\cd)\big)\big|_{L^2_\dbF(\Omega;L^1\!(\nu(t),T;\dbR^n))}\\
&\leq   \cC\big[|f(\cd)|_{L^2_\dbF(\Omega; L^1(\nu(t),T;\dbR^n))}+|y^\pi_T|_{L^2_{\cF_T}(\Omega;\dbR^n)}\big]
              |(\chi_{[\nu(t),T]}u,\chi_{[\nu(t),T]}v,\eta)|_{\mathfrak{S}_t}\,.
\end{align*}
Hence, by Riesz's representation theorem, there is a unique $\big(y^{t}(\cd),Y^{t}(\cd),\varsigma^{t}\big)\in \mathfrak{S}_t$ such that
\bel{wang8}
\bal
\me\big[\lan z^t(T),y^\pi_T \ran\big]-\me\int_{\nu(t)}^T\big\lan z^t(s), f\big(\nu(s)\big)\big\ran \rd s
=\me\int_{\nu(t)}^T \Big[\lan u(s),y^{t}(s) \ran+\lan v(s),Y^{t}(s)\ran\Big] \rd s+\me\big[\lan \eta, \varsigma^{t}\ran\big]\,.
\eal
\ee
Moreover, there exists a positive constant $\cC$, independent of $t$, such that
\bel{wang9}
\bal
|(y^{t}(\cd),Y^{t}(\cd),\varsigma^{t})|_{\mathfrak{S}_t}
\leq \cC\big[|f(\cd)|_{L^2_\dbF(\Omega; L^1(\nu(t),T;\dbR^n)) }+|y^\pi_T|_{L^2_{\cF_T}(\Omega;\dbR^n)}\big]\,.
\eal
\ee

Taking $t=T$ in \eqref{wang8}, we have $\me\big[\lan y^\pi_T-\varsigma^T, \eta\ran\big]=0$, for any $\eta \in [\cH^M(N)]^n$. Since $y^\pi_T,\,\varsigma^T\in [\cH^M(N)]^n$, it follows that $\varsigma^T=y^\pi_T\,\as$

\ss\no
{\bf Step 1(b).} We claim that, for any $t_a,t_b$ satisfying $0\leq t_a <t_b < T$,
\bel{w10-consistency}
\big(y^{t_a}(\cd),Y^{t_a}(\cd)\big)=\big(y^{t_b}(\cd),Y^{t_b}(\cd)\big)\q\mbox{in }\dbH_{N,M}^{n}\big(\nu(t_b),T\big)\times \dbH_{N,M-1}^{n\times d}\big(\nu(t_b),T\big)\,.
\ee


Fix any $\bar u(\cd)\in\dbH_{N,M}^n$
and choose $\eta=0,\,u(\cd)=\chi_{[\nu(t_b),T)}\bar u(\cd),\,v(\cd)=0$.
Applying \eqref{wang8} with $t=t_b$ gives
\bel{wang10}
\me\big[\lan z^{t_b}(T),y^\pi_T\ran\big]-\me\int_{\nu(t_b)}^T\big\lan z^{t_b}(s),f\big(\nu(s)\big)\big\ran \rd s
  =\me\int_{\nu(t_b)}^T\lan \bar u(s),y^{t_b}(s)\ran \rd s\,.
\ee
Applying \eqref{wang8} with $t=t_a$ and the same choice of $\big(\eta,u(\cd),v(\cd)\big)$ gives
\bel{wang11}
\me\big[\lan z^{t_a}(T),y^\pi_T\ran\big]-\me\int_{\nu(t_a)}^T\big\lan z^{t_a}(s),f\big(\nu(s)\big)\big\ran \rd s
=\me\int_{\nu(t_b)}^T\lan \bar u(s),y^{t_a}(s)\ran \rd s\,.
\ee
For this choice of $u(\cd)$ and $v(\cd)$,
\beq
z^{t_a}(s)=
\left\{\!\!\!
\begin{array}{ll}
\ds z^{t_b}(s)\,, \q& s\in[\nu(t_b),T]\,,\\
\ns\ds 0\,,    &s\in [\nu(t_a),\nu(t_b))\,.
\end{array}
\right.
\eeq
Hence, comparing \eqref{wang10} with \eqref{wang11}, we deduce that
\beq
\me\int_{\nu(t_b)}^T\big\lan \bar u(s),y^{t_a}(s)-y^{t_b}(s)\big\ran \rd s
=\me\int_{\nu(t_b)}^T\big\lan
\bar u(s),\chi_{[\nu(t_b),T]}y^{t_a}(s)-y^{t_b}(s)\big\ran \rd s=0\,.
\eeq
Taking
\beq
\bar u(s)=
\left\{\!\!\!
\begin{array}{ll}
\ds \chi_{[\nu(t_b),T)}y^{t_a}(s)-y^{t_b}(s)\,,\q & s\in [\nu(t_b),T)\,,\\
\ns \ds 0\,, & s\in [0,\nu(t_b))\,,
\end{array}
\right.
\eeq
we conclude that $y^{t_a}(\cd)=y^{t_b}(\cd)$ in $\dbH_{N,M}^n\big(\nu(t_b),T\big)$.

Similarly, by taking $\eta=0,\,u(\cd)=0,\,v(\cd)=\chi_{[\nu(t_b),T)}\bar v(\cd)$ with arbitrary $\bar v\in\dbH_{N,M-1}^{n\times d}$, we obtain $Y^{t_a}(\cd)=Y^{t_b}(\cd)$ in $\dbH_{N,M-1}^{n\times d}\big(\nu(t_b),T\big)$. This proves the claim \eqref{w10-consistency}.

%

Set
\beq
y(\cd)=y^0(\cd)\,,\q Y(\cd)=Y^0(\cd)\,.
\eeq
Then, for each $t\in [0,T)$,
\beq
\big(y^{t}(\cd),Y^{t}(\cd)\big)=\big(y(\cd),Y(\cd)\big)\q \mbox{in } \dbH_{N,M}^n\big(\nu(t),T\big)\times \dbH_{N,M-1}^{n\times d}\big(\nu(t),T\big)\,.
\eeq

\ss\no
{\bf Step 1(c).} We claim that, for any $s\in [0,T)$,
\beq
\varsigma^{\nu(s)}=y\big(\nu(s)\big) \q \as
\eeq

Let $s\in[t_k,t_{k+1})$ for some $k=0,1,\cds,N-1$. Then
$\nu(s)=t_k$. Fix $\sigma\in[\cH^M(k)]^n$.
First, in \eqref{wang8}, take $t=t_k$, $\eta=0,\,u(\cd)=\chi_{[t_k,t_{k+1})}(\cd)\sigma $, $v(\cd)=0$.
For this choice,
\beq
z^{t_k}(s)=
\left\{\!\!\!
\begin{array}{ll}
\ds \int_{t_k}^{\mu(s)}u(t)\rd t=(t_{k+1}-t_k)\sigma\,, \q&s\in (t_k,T]\,,\\
\ns \ds 0\,, & s\in [0,t_k]\,.
\end{array}
\right.
\eeq
Consequently,
\bel{wang21}
\bal
&\me\big[\lan \t \sigma,y^\pi_T \ran\big]-\me\int_{t_k}^T\big\lan \t \sigma, f\big(\nu(s)\big)\big\ran \rd s
=\me\int_{t_k}^{t_{k+1}} \lan \sigma,y(s) \ran \rd s=\t \me\big[\lan\sigma,y(t_k)\ran\big]\,,
\eal
\ee
where the last equality follows from $y(\cd)\in\dbH_{N,M}^n$.

Next, in \eqref{wang8}, take $t=t_k$, $u(\cd)=v(\cd)=0$, and $\eta=\t \sigma$.
Then
\bel{wang22}
\bal
\me\big[\lan \t \sigma,y_T^\pi \ran\big]-\me\int_{t_k}^T\big\lan \t \sigma, f\big(\nu(s)\big)\big\ran \rd s
=\me\big[\lan \t \sigma,\varsigma^{t_k}\ran\big]\,.
\eal
\ee
Comparing \eqref{wang21} and \eqref{wang22}, we obtain $\varsigma^{t_k}=y(t_k)\,\,\as$
This proves the claim.

Combining \eqref{wang8} with the claims in Steps 1(b) and 1(c), we conclude that, for each $t\in[0,T)$,
\bel{w10e5}
\bal
&\me\big[\lan z^t(T),y_T^\pi \ran\big]-\me\big[\lan \eta, y(\nu(t))\ran\big]
=\me\int_{\nu(t)}^T \Big[\big\lan z^{t}(s), f\big(\nu(s)\big)\big\ran+\lan u(s),y(s) \ran+\lan v(s),Y(s)\ran\Big] \rd s\,.
\eal
\ee

\ss\no
{\bf Step 1(d).} Setting $t=0$ and $\eta=0$ in \eqref{w10e5} shows that $\big(y(\cd),Y(\cd)\big)$ is a finite transposition solution to BSDE \eqref{bsde0} (satisfying the variational identity \eqref{w9e3}).

It remains to prove uniqueness.
Let $\big(y_1(\cd),Y_1(\cd)\big), \big(y_2(\cd),Y_2(\cd)\big)\in \dbH_{N,M}^{n}\times \dbH_{N,M-1}^{n\times d}$ be two finite transposition solutions. Subtracting their variational identities yields
\beq
0
=\me\int_{0}^T  \Big[ \langle u(t),y_1(t)-y_2(t)\rangle
+\langle v(t), Y_1(t)-Y_2(t)\rangle \Big] \rd t\,.
\eeq
Taking $u(\cd)=y_1(\cd)-y_2(\cd),\,v(\cd)=Y_1(\cd)-Y_2(\cd)$, we conclude that
\beq
\big(y_1(\cd),Y_1(\cd)\big)=\big(y_2(\cd),Y_2(\cd)\big)\q\mbox{in } \dbH_{N,M}^{n}\times \dbH_{N,M-1}^{n\times d}\,.
\eeq
This proves the well-posedness result for the linear BSDE.

\ss\no
{\bf Step 2. Semilinear BSDEs.} We now establish the well-posedness result for the semilinear BSDE.

\ss\no
{\bf Step 2(a).}
Let $\big(p(\cd),q(\cd)\big)\in \dbH_{N,M}^{n}\times \dbH_{N,M-1}^{n\times d}$. By Assumption {\bf(B)}, it follows that
\beq
f_N\big(\cd,p(\cd),q(\cd)\big) \deq f\big(\nu(\cd), \chi_{[t_{N-1},T)}p(\cd),\chi_{[t_{N-1},T)}q(\cd)\big) \in L^2_\dbF\big(\Omega;L^1(0,T;\dbR^n)\big)\,.
\eeq
Therefore, by the well-posedness result established in Step 1, the following linear BSDE admits a unique finite transposition solution $\big(y_N(\cd),Y_N(\cd)\big) \in \dbH_{N,M}^{n}\times \dbH_{N,M-1}^{n\times d}$
\bel{eq2.3.3}
\lt\{
\bal
&\mathrm d y(t)=f_N\big(t,p(t),q(t)\big)\rd t+Y(t)\rd W(t)\,,\q t\in [0,T]\,,\\
&y(T)=y^\pi_T\,.
\eal
\rt.
\ee

Define ${\mathcal T}: \dbH_{N,M}^{n}\times \dbH_{N,M-1}^{n\times d} \to
\dbH_{N,M}^{n}\times \dbH_{N,M-1}^{n\times d}$ by
$
{\mathcal T}\big(p(\cd),q(\cd)\big)=\big(y_N(\cd),Y_N(\cd)\big)$.
Let $\big(\wt p(\cd),\wt q(\cd)\big)\in \dbH_{N,M}^{n}\times \dbH_{N,M-1}^{n\times d}$ and write
$\big(\wt y_N(\cd),\wt Y_N(\cd)\big)={\mathcal T}\big(\wt{p}(\cd),\wt{q}(\cd)\big)$.
Set
\begin{align*}
&\hat{y}_N(\cd)=y_N(\cd)-\wt{y}_N(\cd)\,,\quad \hat{Y}_N(\cd)=Y_N(\cd)-\wt{Y}_N(\cd)\,,\\
\mbox{and}\q& \hat{f}_N(\cdot)=f_N\big(\cd,p(\cd),q(\cd)\big)-f_N\big(\cd,\wt{p}(\cd),\wt{q}(\cd)\big) \,.
 \end{align*}
It is clear that $\big(\hat y_N(\cd),\hat Y_N(\cd)\big)$ is the finite transposition solution in $\dbH_{N,M}^{n}\times \dbH_{N,M-1}^{n\times d}$ to the following equation:
\beq
\left\{
\bal
&\mathrm d\hat{y}_N(t)=\hat{f}_N(t)\rd t+\hat{Y}_N(t)\rd W(t)\,,\q t\in [0,T]\,,\\
&\hat{y}_N(T)=0\,.
\eal
\right.
\eeq
The solution $(\hat{y}_N(\cd),\hat{Y}_N(\cd))$ therefore satisfies the following variational identity:
\begin{eqnarray}\label{eq2.3.5}
0=
\me\int_{0}^T \Big[  \langle z(t ),\hat{f}_N(t )\rangle+\langle u(t ),\hat{y}_N(t )\rangle
+\langle v(t ),\hat{Y}_N(t)\rangle \Big] \rd t \,,
\end{eqnarray}
where $\big(u(\cd),\,v(\cd)\big)\in \dbH_{N,M}^{n}\times \dbH_{N,M-1}^{n\times d}$, and
$z(\cd)$ satisfies \eqref{w9e4}.
Then, by \eqref{eq2.3.5} and Lemma \ref{w10l1}, we arrive at
\begin{align*}
&\bigg|\me\int_{0}^T \Big[ \langle u(t ),\hat{y}_N(t )\rangle
+\langle v(t ),\hat{Y}_N(t)\rangle \Big] \rd t\bigg|
\leq |\hat{f}_N(\cd)|_{L_{\dbF}^2(\Omega;L^1(0,T;\dbR^n))}
   |z(\cd)|_{L^2_\dbF(\Omega;L^\infty(t_{N-1},T;\dbR^n))}\nonumber\\
&\qq\qq\leq  \cC |\hat{f}_N(\cd)|_{L_{\dbF}^2(\Omega;L^1(0,T;\dbR^n))}
\big|\big(u(\cd),v(\cd)\big)\big|_{\dbH_{N,M}^{n}\times \dbH_{N,M-1}^{n\times d}}\,.
\end{align*}
Hence,
\beq
\big|\big(\hat{y}_N(\cd),\hat{Y}_N(\cd)\big)\big|_{ \dbH_{N,M}^{n}\times \dbH_{N,M-1}^{n\times d}}
\leq \cC |\hat{f}_N(\cd)|_{L_{\dbF}^2(\Omega;L^1(0,T;\dbR^n))}\,.
\eeq
%
On the other hand, Assumption {\bf(B)} implies that
\begin{align*}
&|\hat{f}_N(\cd)|^2_{L_{\dbF}^2(\Omega;L^1(0,T;\dbR^n))}\\
&=\me\Big[\Big(\int_0^T \big| f\big(\nu(t),\chi_{[t_{N-1},T)}p(t),\chi_{[t_{N-1},T)}q(t)\big)       -f\big(\nu(t),\chi_{[t_{N-1},T)}\wt{p}(t),\chi_{[t_{N-1},T)}\wt{q}(t)\big)\big| \rd t\Big)^2\Big]\\
&\leq    2L^2\tau  \big|\big(p(\cd)-\wt p(\cd), q(\cd)-\wt q(\cd)\big)\big|^2_{\dbH_{N,M}^{n}\times \dbH_{N,M-1}^{n\times d}}\,.
\end{align*}
The preceding two estimates show that ${\mathcal T}$ is a contraction whenever $\tau\leq\tau_0$ for some sufficiently small $\tau_0>0$. Consequently, by the Banach fixed-point theorem, ${\mathcal T}$ has a unique fixed point $\big(y_N(\cd),Y_N(\cd)\big)\in\dbH_{N,M}^{n}\times \dbH_{N,M-1}^{n\times d}$, which is the unique finite transposition solution to the following semilinear BSDE
\beq
\lt\{
\bal
&\mathrm d y(t)=f_N\big(t,y(t),Y(t)\big)\rd t+Y(t)\rd W(t)\,,\q t\in [0,T]\,,\\
&y(T)=y^\pi_T\,.
\eal
\rt.
\eeq

In the same vein, we conclude that, for any $k=N-1, N-2,\cds,1$, the following BSDE admits a unique finite transposition solution $\big(y_{k}(\cd),Y_{k}(\cd)\big)\in \dbH_{N,M}^{n}\times \dbH_{N,M-1}^{n\times d}$:
\beq
\left\{
\bal
&\mathrm d y(t)=f_{k}\big(t,y(t),Y(t)\big)\rd t+Y(t)\rd W(t)\,, \q t\in [0,T]\,,\\
&y(T)=y^\pi_T\,,
\eal
\right.
\eeq
where
\beq
&f_k\big(\cd,y(\cd),Y(\cd)\big)\deq f\big(\nu(\cd),\bar y_k(\cd),\bar Y_k(\cd)\big)\,,\\
 &
\bar y_k(\cd)=\chi_{[t_{k-1},t_{k})}y(\cd)+\chi_{[t_{k},t_{k+1})}y_{k+1}(\cd)+\cds+\chi_{[t_{N-1},T)}y_{N}(\cd)\,,\\
\mbox{and}\q&\bar Y_k(\cd)=\chi_{[t_{k-1},t_{k})}Y(\cd)+\chi_{[t_{k},t_{k+1})}Y_{k+1}(\cd)+\cds+\chi_{[t_{N-1},T)}Y_{N}(\cd)\,.
\eeq

\ss\no
{\bf Step 2(b).}  We claim that, for any $1\leq k<l\leq N$,
\bel{w10e10}
\big(y_k(\cd),Y_k(\cd)\big)=\big(y_l(\cd),Y_l(\cd)\big)\q\mbox{in } \dbH_{N,M}^n(t_{l-1},T)\times \dbH_{N,M-1}^{n\times d}(t_{l-1},T)\,.
\ee

By the definition of the finite transposition solution (see Definition \ref{defi-of-ftm}), we have
\bel{w11e1}
\me\big[\lan z(T),y_T^\pi\ran\big]
\!=\!\me\!\int_{0}^T\!  \Big[\big \langle z(t), f_k\big(\nu(t),y_k(t),Y_k(t)\big)\big\rangle
\!+\! \langle u(t),y_k(t)\rangle \!+\!\langle v(t), Y_k(t)\rangle \Big] \rd t\,,
\ee
and
\bel{w11e2}
\me\big[\lan z(T),y_T^\pi\ran\big]
\!=\!\me\!\int_{0}^T\!  \Big[ \big\langle z(t), f_l\big(\nu(t),y_l(t),Y_l(t)\big)\big\rangle
\!+\! \langle u(t),y_l(t)\rangle \!+\!\langle v(t), Y_l(t)\rangle \Big] \rd t\,.
\ee
Choose $u(\cd)$ and $v(\cd)$ such that $\chi_{[0,t_{l-1})}u(\cd)=0,\,\chi_{[0,t_{l-1})}v(\cd)=0$. Then $z(\cd)=0$ on $[0,t_{l-1})$, and subtracting \eqref{w11e2} from \eqref{w11e1} gives
\bel{w11e3}
\bal
&\me\int_{t_{l-1}}^T  \Big[ \big\langle z(t), f_k\big(\nu(t),y_k(t),Y_k(t)\big)\big\rangle
+ \langle u(t),y_k(t)\rangle +\langle v(t), Y_k(t)\rangle \Big] \rd t\\
&=\me\int_{t_{l-1}}^T  \Big[ \big\langle z(t), f_l\big(\nu(t),y_l(t),Y_l(t)\big)\big\rangle
+ \langle u(t),y_l(t)\rangle +\langle v(t), Y_l(t)\rangle \Big] \rd t\,.
\eal
\ee
The definitions of $f_k(\cd)$ and $f_l(\cd)$ imply that, for $t\in[t_{l-1},T)$,
\beq
f_k\big(\nu(t),y_k(t),Y_k(t)\big)=f_l\big(\nu(t),y_l(t),Y_l(t)\big)\,.
\eeq
Therefore, \eqref{w11e3} leads to
\beq
\me\int_{t_{l-1}}^T  \Big[  \langle u(t),y_k(t)-y_l(t)\rangle
 +\langle v(t), Y_k(t)-Y_l(t)\rangle \Big] \rd t
=0\,,
\eeq
which, together with the arbitrariness of $u(\cd)$ and $v(\cd)$, gives \eqref{w10e10}.

Set
\beq
y(\cd)=y_1(\cd)\,,\q Y(\cd)=Y_1(\cd)\,.
\eeq
It follows immediately that $f\big(\nu(t),y(t),Y(t)\big)=f_1\big(\nu(t),y_1(t),Y_1(t)\big)$, and the pair $\big(y(\cd),Y(\cd)\big)$ satisfies
\eqref{w9e3}.

It remains to verify uniqueness for the original semilinear problem.
Let the pair $\big(\wt y(\cd),\wt Y(\cd)\big)$ be another finite transposition solution to \eqref{bsde0}. On the last interval $[t_{N-1},T]$, both $\big(y(\cd),Y(\cd)\big)$ and $\big(\wt y(\cd),\wt Y(\cd)\big)$ solve the finite transposition problem corresponding to $f_N(\cd,\cd,\cd)$. By the uniqueness established in Step 2(a), they coincide on $[t_{N-1},T)$.
Proceeding backward interval by interval and using the uniqueness of the problems corresponding to $f_{N-1}(\cd,\cd,\cd),\cds,f_1(\cd,\cd,\cd)$, we obtain $\big(\wt y(\cd),\wt Y(\cd)\big)=\big(y(\cd),Y(\cd)\big)$.
Hence the finite transposition solution to \eqref{bsde0} is unique.
This completes the proof.
\end{proof}

\br{w808r1}
In what follows, let $\tau_0>0$ be the constant given in Theorem \ref{w9t1}, and set
\bel{w-Nstar}
N_*
\deq
\max\big\{\lceil T/\tau_0\rceil,\,\lceil T\rceil,\,\lceil 72L^2T\rceil\big\}\,.
\ee
Then, for every $N\geq N_*$, $\tau=\frac{T}{N}\leq \tau_0\wedge 1\wedge\frac{1}{72L^2}$.
Henceforth, whenever a finite transposition solution is considered, we assume that $N\geq N_*$ unless otherwise stated.

\er

By Lemma \ref{w10l1}, for any $t\in[0,T]$, $z(t)\in \big[\cH^M\big(\pi(t)\big)\big]^n\subseteq [\cH^M(N)]^n$. Moreover, for any $\big(y(\cd),Y(\cd)\big)\in \dbH_{N,M}^n\times\dbH_{N,M-1}^{n\times d}$, we have $\big(y(t),Y(t)\big)\in \big[\cH^M\big(\pi(t)\big)\big]^n\times\big[\cH^{M-1}\big(\pi(t)\big)\big]^{n\times d}$ for each $t\in [0,T)$. Since $f(\cd,\cd,\cd)$ is deterministic, it follows that $ f\big(\nu(t),y(t),Y(t)\big) \in L_{\cG_{\nu(t)}}^2(\Omega;\dbR^n) \subseteq L_{\cF_T}^2(\Omega;\dbR^n)$. Since $\G_M$, defined in \eqref{w417e1}, is the $L^2$-orthogonal projection onto $\bigoplus_{m=0}^M(\cH_m)^n$, we obtain
\beq
\me\int_{0}^T  \big\lan z(t), f\big(\nu(t),y(t),Y(t)\big) \big\ran \rd t
=\me\int_{0}^T \big \lan z(t), \G_M f\big(\nu(t),y(t),Y(t)\big) \big\ran \rd t\,.
\eeq
We therefore obtain the following corollary.

\begin{corollary}\label{w19c1}
Under Assumptions {\rm\bf(P)} and {\rm\bf(B)}, for any $N\geq N_*$ and $M\in\dbN$, the following two BSDEs have the same finite transposition solution in $\dbH_{N,M}^{n}\times\dbH_{N,M-1}^{n\times d}$:
\bel{w19a3}
\left\{
\bal
&\mathrm d y(t)=f\big(t,y(t),Y(t)\big)\rd t+Y(t)\rd W(t)\,,\quad t\in [0,T]\,,\\
&y(T)= y^\pi_T\in [\cH^M(N)]^n\,,
\eal
\right.
\ee

\bel{w19a4}
\left\{
\bal
&\mathrm d y(t)=\G_M f\big(t,y(t),Y(t)\big)\rd t+Y(t)\rd W(t)\,,\quad t\in [0,T]\,,\\
&y(T)= y^\pi_T\in [\cH^M(N)]^n\,.
\eal
\right.
\ee
\end{corollary}

\br{w21r1} {\rm (1)}
By Theorem \ref{w9t1} and Corollary \ref{w19c1}, the finite transposition solution associated with the partition parameter $N$ and the Wiener chaos truncation order $M$ depends on both $N$ and $M$. Henceforth, we denote it by $\big(y_{N,M}(\cd),Y_{N,M-1}(\cd)\big)$, and write the adapted solution to \eqref{bsde0} as $\big(y(\cd),Y(\cd)\big)$.
Whenever pointwise values at the terminal time are involved, we use the convention $y_{N,M}(T)=y_T^\pi$.
%


{\rm (2)}
Although $\G_M$ acts on random variables rather than pointwise in $\omega$, BSDE \eqref{w19a4} is still well-posed.
Indeed, $\G_M$ preserves $\cF_t$-measurability of $\cF_t$-measurable random variables and is an $L^2$-contraction.
Consequently, Assumption {\rm\bf(B)} gives
$\me\big[ \big|\G_M f(t,y_1,Y_1)-\G_M f(t,y_2,Y_2)\big|^2\big]
\leq \cC\big(\me\big[|y_1-y_2|^2\big]+\me\big[|Y_1-Y_2|^2\big]\big)$.
A standard Picard contraction argument in $L^2_{\dbF}(0,T;\dbR^n)\times L^2_{\dbF}(0,T;\dbR^{n\times d})$ therefore yields a unique adapted solution $\big(y_M(\cd),Y_M(\cd)\big)$ to \eqref{w19a4}.
In general, this adapted solution differs from the adapted solution of \eqref{w19a3}, and both adapted solutions generally differ from the finite transposition solution $\big(y_{N,M}(\cd),Y_{N,M-1}(\cd)\big)$.

{\rm (3)}
The proof of Theorem \ref{w9t1} also applies when different Wiener chaos truncation orders are used for the two components. More precisely, let $N,M_1,M_2\in\dbN$, and take the test processes $\big(u(\cd),\,v(\cd)\big)\in \dbH_{N,M_1}^n\times\dbH_{N,M_2}^{n\times d}$.
Then the corresponding test process belongs to $\dbH_{N,M_1\vee(M_2+1)}^n$.
Therefore, if $y(T)=y^\pi_T\in \big[\cH^{M_1\vee(M_2+1)}(N)\big]^n$, then, for sufficiently small $\tau$, the BSDE with terminal condition $y(T)=y_T^\pi$ admits a unique finite transposition solution $\big(y_{N,M_1}(\cd),Y_{N,M_2}(\cd)\big)\in \dbH_{N,M_1}^n\times\dbH_{N,M_2}^{n\times d}$.

\er

\section{Error analysis of the finite transposition method}\label{convergence}

This section is devoted to the error analysis of the finite transposition method. In Sections \ref{ftm1} and \ref{sec-picard}, we establish its relationship with the Euler method and derive convergence rates under appropriate assumptions; see Theorems \ref{ftm-rate} and \ref{w602t1}. As a consequence of Theorem \ref{ftm-rate}, we also obtain convergence under weaker assumptions; see Theorem \ref{ftm-convergence} in Section \ref{sec-ftm-convergence}.
Section \ref{Markov-BSDE} presents the convergence analysis for Markovian BSDEs, whereas Section \ref{MC} addresses the implementation of the finite transposition method by means of Monte Carlo approximation.

\subsection{Convergence rate of the finite transposition method for general BSDEs}\label{ftm1}

Remark \ref{w9r1} provides the following orthonormal basis of $\dbH_{N,M}$: $\big\{e_{k,i}\,|\,1\leq i\leq \dim\big(\cH^M(k)\big),\,0\leq k\leq N-1\big\}$.
By Theorem \ref{w9t1}, the unique finite transposition solution $\big(y_{N,M}(\cd),Y_{N,M-1}(\cd)\big)\in\dbH_{N,M}^{n}\times\dbH_{N,M-1}^{n\times d}$ admits the expansion
\bel{w23e1}
\bal
y_{N,M}(\cd)&=\sum_{k=0}^{N-1}\sum_{i=1}^{\dim(\cH^M(k))}a_{k,i}e_{k,i}(\cd)\,,\q a_{k,i}\in\dbR^n\,,\\
\mbox{and}\q Y_{N,M-1}(\cd)&=\sum_{k=0}^{N-1}\sum_{i=1}^{\dim(\cH^{M-1}(k))}b_{k,i}e_{k,i}(\cd)\,,\q b_{k,i}\in\dbR^{n\times d}\,.
\eal
\ee
Thus, the finite transposition solution is completely determined by the vector coefficients $\{a_{k,i}\}_{k,i}$ and the matrix coefficients $\{b_{k,i}\}_{k,i}$.
By Definition \ref{defi-of-ftm}, these coefficients are characterized by the following variational identity
\bel{w23e2}
\bal
\me\big[\lan z(T), y^\pi_T \ran\big]
=\me\int_{0}^T  \Big[ \big\langle z(t), f\big(\nu(t),y_{N,M}(t),Y_{N,M-1}(t)\big)\big\rangle
+ \langle u(t),y_{N,M}(t)\rangle
+\langle v(t), Y_{N,M-1}(t)\rangle \Big] \rd t\,,
\eal
\ee
where $\big(u(\cd), v(\cd)\big)\in \dbH_{N,M}^n\times \dbH_{N,M-1}^{n\times d}$ and $z(\cd)$ solves \eqref{w9e4}.

To be precise, let $c\in\dbR^n$ be arbitrary and take $u(\cd)=e_{k,i}(\cd)c$ and $v(\cd)=0$, where $1\leq i\leq\dim(\cH^M(k))$. Then
\beq
z(t)=\int_0^{\mu(t)} e_{k,i}(s)c \rd s=
\left\{\!\!\!
\begin{array}{ll}
0\,, \q & t\in[0,t_k]\,,\\
\sqrt{\tau } h_{k,i}c\,, \q & t\in(t_k,T]\,.\\
\end{array}
\right.
\eeq
Substituting this choice into \eqref{w23e2}, we obtain
\begin{align*}
&\me\big[\lan \sqrt{\tau } h_{k,i}c, y^\pi_T\ran\big]\\
&=\me\int_{t_k}^T\big\lan\sqrt{\tau } h_{k,i}c,  f\big(\nu(t),y_{N,M}(t),Y_{N,M-1}(t)\big)\big\ran \rd t
+\frac{1}{\tau }\me \int_{t_k}^{t_{k+1}}\!\!\Big\lan h_{k,i}c, \sum_{l=0}^{N-1}\sum_{j=1}^{\dim(\cH^M(l))}a_{l,j}\chi_{[t_l,t_{l+1})}(t)h_{l,j}\Big\ran \rd t\\
&=\me \int_{t_k}^T\big\lan \sqrt{\tau } h_{k,i}c,  f\big(\nu(t),y_{N,M}(t),Y_{N,M-1}(t)\big)\big\ran \rd t
+\frac{1}{\tau }\me \int_{t_k}^{t_{k+1}}\big\lan h_{k,i}c, a_{k,i}h_{k,i}\big\ran \rd t \\
&=\me \int_{t_k}^T\big\lan c,\sqrt{\tau } h_{k,i}  f\big(\nu(t),y_{N,M}(t),Y_{N,M-1}(t)\big)\big\ran \rd t
+\lan c,a_{k,i}\ran\,.
\end{align*}
Here we use the fact that, for each fixed $k$, $\cE_0(k)=\{h_{k,i}\}_{i=1}^{\dim(\cH^M(k))}$ is an orthonormal basis of $\cH^M(k)$; see Lemma \ref{w512l1} and Remark \ref{w9r1}.
Notice that the union of the families $\cE_0(k)$ over different values of $k$ need not be orthonormal. The reduction above is instead justified by the disjoint time supports of the functions $e_{l,j}(\cd)$.
Since $c\in\dbR^n$ is arbitrary, we conclude that
\bel{w23e6}
a_{k,i}
=\me\big[\sqrt{\tau}h_{k,i}y^\pi_T\big]
-\me\int_{t_k}^T \sqrt{\tau}h_{k,i}
f\big(\nu(t),y_{N,M}(t),Y_{N,M-1}(t)\big) \rd t
\in\dbR^n\,.
\ee

\ss
Similarly, let $C\in\dbR^{n\times d}$ be arbitrary and take $u(\cd)=0$ and $v(\cd)=e_{k,i}(\cd)C$, where $1\leq i\leq\dim(\cH^{M-1}(k))$. Then
\beq
z(t)=\int_0^{\nu(t)} e_{k,i}(s)C\rd W(s)=
\left\{\!\!\!
\begin{array}{ll}
0\,,\q & t\in[0,t_{k+1})\,,\\
\ns\ds\frac{h_{k,i}}{\sqrt{\tau}}C\D_{k+1}W\,,\q & t\in[t_{k+1},T]\,,\\
\end{array}
\right.
\eeq
and
\begin{align*}
\me\Big[\Big\lan \frac{h_{k,i}}{\sqrt{\tau}}C\D_{k+1}W, y^\pi_T\Big\ran\Big]
&=\me\int_{t_{k+1}}^T\Big\lan \frac{h_{k,i}}{\sqrt{\tau}}C\D_{k+1}W,  f\big(\nu(t),y_{N,M}(t),Y_{N,M-1}(t)\big)\Big\ran \rd t\\
&\q\qq+\frac{1}{\tau }\me\int_{t_k}^{t_{k+1}}\Big\lan h_{k,i}C, \sum_{l=0}^{N-1}\sum_{j=1}^{\dim(\cH^{M-1}(l))}b_{l,j}\chi_{[t_l,t_{l+1})}(t)h_{l,j}\Big\ran \rd t\\
&=\me\int_{t_{k+1}}^T\Big\lan C\D_{k+1}W, \frac{h_{k,i}}{\sqrt{\tau}} f\big(\nu(t),y_{N,M}(t),Y_{N,M-1}(t)\big)\Big\ran \rd t
+\lan C,b_{k,i}\ran\,,
\end{align*}
where the inner product between matrices is the Frobenius inner product.

For any $x\in\dbR^n$, $C\in\dbR^{n\times d}$, and $q\in\dbR^d$, we have
$\lan Cq,x\ran_{\dbR^n}=\lan C,xq^\top\ran_{\dbR^{n\times d}}$.
Therefore,
\begin{align*}
\lan C,b_{k,i}\ran
&=
\me\Big[ \Big\lan C, \frac{h_{k,i}}{\sqrt{\tau}} y^\pi_T(\D_{k+1}W)^\top \Big\ran \Big]\\
&\q- \me\int_{t_{k+1}}^T
\Big\lan C, \frac{h_{k,i}}{\sqrt{\tau}} f\big(\nu(t),y_{N,M}(t),Y_{N,M-1}(t)\big) (\D_{k+1}W)^\top \Big\ran \rd t\,.
\end{align*}
Since $C\in\dbR^{n\times d}$ is arbitrary, we arrive at
\bel{w23e9}
b_{k,i}\!=\!\me\Big[ \frac{h_{k,i}}{\sqrt{\tau }}  y^\pi_T(\D_{k+1}W)^\top\!\Big]
  \!-\!\me\!\int_{t_{k+1}}^T\!\frac{h_{k,i}}{\sqrt{\tau }} f\big(\nu(t),y_{N,M}(t),Y_{N,M-1}(t)\big)(\D_{k+1}W)^\top\! \rd t\,.
\ee

Since the generator $f(\cd,\cd,\cd)$ contains $y_{N,M}(\cd)$ and $Y_{N,M-1}(\cd)$, the relations \eqref{w23e6} and \eqref{w23e9} determine the coefficients only implicitly. This nonlinear algebraic system will be solved by Picard iteration in Section \ref{sec-picard}.

To derive convergence rates for the finite transposition method, we impose the following additional regularity assumption on $y_T$ and $f(\cd,\cd,\cd)$.
\begin{enumerate}
\item [{\bf (H)$_m$}]
There exist $q>4$ and a positive constant $L$, both independent of $m$, such that $y_T\in\dbD^{2,q}(\dbR^n)\cap\dbD^{m,2}(\dbR^n)$ and
\bel{w1228e1}
\bal
&\max_{1\leq\ell\leq d}\me\big[|D_{\th_1}^{\ell}y_T-D_{\th_2}^{\ell}y_T|^2\big]
\leq L|\th_1-\th_2|\,,\q \forall\, \th_1,\th_2\in [0,T]\,,\\
&\sup_{1\leq\ell\leq d}\sup_{\th\in[0,T]}\me\big[|D_{\th}^{\ell}y_T|^q\big]
+\sup_{1\leq\ell_1,\ell_2 \leq d}\sup_{\th_1,\th_2\in[0,T]}
\me\big[|D_{\th_1,\th_2}^{\ell_1,\ell_2}y_T|^q\big]\leq L\,.
\eal
\ee
Moreover, for any $k=1,2,\cds,m$, any
$(\ell_1,\cds,\ell_k)\in\{1,2,\cds,d\}^k$, and all
$(\th_1,\cds,\th_k)$, $(\bar \th_1,\cds,\bar \th_k)\in[0,T]^k$,
\bel{lip-h-2}
\bal
&\big|\me\big[D_{\th_1,\cds,\th_k}^{\ell_1,\cds,\ell_k}y_T\big]
-\me\big[D_{\bar \th_1,\cds,\bar \th_k}^{\ell_1,\cds,\ell_k}y_T\big]\big|
\leq L^{1/2}\sum_{j=1}^k|\th_j-\bar \th_j|^{1/2}\,.
\eal
\ee
The deterministic function $f:[0,T]\times\dbR^n\times\dbR^{n\times d}\longrightarrow\dbR^n$ is $\frac 1 2$-H\"older continuous with respect to the first variable 
and is twice continuously differentiable with respect to the second and third variables, with uniformly bounded first- and second-order derivatives.

\end{enumerate}

\br{w1102r1}
Under Assumption {\rm\bf (H)$_1$}, the argument of \cite[Theorem 2.6]{Hu-Nualart-Song11}, extended componentwise to the present multidimensional setting, yields
\bel{w1102e1}
\me\big[|Y(t)-Y(s)|^2\big]\leq  \cC|t-s|\,,\q \forall\,t,s\in[0,T]\,.
\ee
To derive the convergence rate of the Euler scheme associated with a given partition $I_\t$, however, it suffices to establish the following weaker $L^2$-regularity estimate
\bel{w709e2}
\max_{0\leq k\leq N-1}
\sup_{t\in[t_k,t_{k+1}]}
\me\big[|y(t)-y(t_k)|^2\big]
+
\sum_{k=0}^{N-1}\me\int_{t_k}^{t_{k+1}}|Y(t)-\bar Y(t_k)|^2\rd t \leq \cC\tau \,,
\ee
where
\bel{w911bar}
\bar Y(t_k)
= \frac{1}{\tau}\me\Big(\int_{t_k}^{t_{k+1}}Y(s)\rd s\,\Big|\,\cF_{t_k}\Big)\,,
\qq k=0,1,\cds,N-1\,;
\ee
see \cite{ZhangJF01,Bouchard-Touzi04}.
To obtain this estimate, the first condition in \eqref{w1228e1} may be replaced by the following weaker condition:
\beq
\max_{1\leq \ell\leq d}\me\big[|D_{\th_1}^{\ell}y_T-D_{\th_2}^{\ell}y_T|^2\big]\leq \cC\t\,,\q \forall\, \th_1,\th_2\in [t_k,t_{k+1}), k=0,1,\cds,N-1\,;
\eeq
see \cite[Theorems 2.3 and 2.6]{Hu-Nualart-Song11}.
For Markovian BSDEs, the pointwise estimate \eqref{w1102e1} may not be available under the standard assumptions. Nevertheless, under the milder Assumptions {\rm\bf (B)} and {\rm\bf (F)}, to be introduced in Section \ref{Markov-BSDE}, estimate \eqref{w709e2} follows from \cite[Theorem 3.4.3]{ZhangJF01}; see also \cite[Lemma 3.2]{Bouchard-Touzi04}.
\er


We now combine \eqref{w23e1}, \eqref{w23e6}, and \eqref{w23e9} to obtain an equivalent conditional-expectation representation of the finite transposition solution and to establish its convergence rate.

\bt{ftm-rate}
Let $M\in\dbN$ and $N\geq N_*$, where $N_*$ is given in \eqref{w-Nstar}, and assume that Assumptions {\rm\bf(P)} and {\rm\bf(B)} hold.
Then, for any $t\in[0,T)$, the finite transposition solution $\big(y_{N,M}(\cd),Y_{N,M-1}(\cd)\big)\in \dbH_{N,M}^n\times\dbH_{N,M-1}^{n\times d}$ admits the representation
\bel{w23e10a}
\lt\{
\bal
& y_{N,M}(t)=\me\Big(y^\pi_T-\int_{\nu(t)}^T
\G_M f\big(\nu(s),y_{N,M}(s),Y_{N,M-1}(s)\big)\rd s
\,\Big|\,\cF_{\nu(t)}\Big)\,,\\
& Y_{N,M-1}(t)=\me\Big(
\frac{1}{\tau}\Big[y^\pi_T
-\int_{\mu(t)}^T
\G_M f\big(\nu(s),y_{N,M}(s),Y_{N,M-1}(s)\big)\rd s\Big]
(\D_{\pi(t)+1}W)^\top
\,\Big|\,\cF_{\nu(t)}\Big)\,.
\eal
\rt.
\ee
Furthermore, if Assumption {\rm\bf(H)$_M$} holds and $y_T^\pi=\G_M^Ny_T$, then
\bel{w416e1}
\bal
\sup_{0\leq t\leq T} \me\big[\big|y_{N,M}(t)\!-\!y\big(\nu(t)\big)\big|^2\big]
\!+\!\me\!\int_0^T\! |Y_{N,M-1}(t)\!-\!Y(t)|^2\rd t
\leq \cC\big(\frac 1 N \!+\!\frac 1 {M+1}\big)\,,
\eal
\ee
where $\big(y(\cd),Y(\cd)\big)$ is the adapted solution to BSDE \eqref{bsde0},
and $\cC$ is independent of $N$ and $M$.
\et

\begin{proof}
Since the proof is long, we divide it into three steps.

\ss
\no
{\bf Step 1.} In this step, we derive the representation \eqref{w23e10a} of the finite transposition solution.

\ss\no
{\bf Step 1(a).}
By \eqref{w23e1} and \eqref{w23e6}, we know that
\begin{align*}
&y_{N,M}(\cd)=\sum_{k=0}^{N-1}\sum_{i=1}^{\dim(\cH^M(k))}a_{k,i}e_{k,i}(\cd)\\
&=\sum_{k=0}^{N\!-\!1}\sum_{i=1}^{\dim(\cH^{M\!}(k))}\!\!\Big(\me \big[ \sqrt{\tau } h_{k,i} y^\pi_T\big]
\!-\!\me\!\int_{t_k}^T\!\! \sqrt{\tau } h_{k,i}  f\big(\nu(t),y_{N,M}(t),Y_{N,M-1}(t)\big) \rd t\Big)e_{k,i}(\cd)\,.
\end{align*}
For the first term on the right-hand side, since
$e_{k,i}(\cd)= \chi_{[t_k,t_{k+1})}(\cd) \frac{1}{\sqrt{\tau}}h_{k,i}$,
we have
\bel{w23e13}
\bal
\me\big[ \sqrt{\tau } h_{k,i}\, y^\pi_T\big] e_{k,i}(\cd)
&=\chi_{[t_k,t_{k+1})}(\cd)\me\big[ h_{k,i}\, y^\pi_T\big] h_{k,i}\,.
\eal
\ee
We next identify the sum over $i$ in \eqref{w23e13}.
Since $y^\pi_T\in[\cH^M(N)]^n$, Remark \ref{w408r2} implies that
\bel{w23e14}
\me\big( y^\pi_T\,\big|\,\cF_{t_k} \big) \in [\cH^M(k)]^n\,.
\ee
Indeed, $y^\pi_T$ is measurable with respect to the Brownian increments on the grid, and taking the conditional expectation with respect to $\cF_{t_k}$ removes all Hermite terms involving increments after $t_k$.
Since $\cE_0(k)=\{h_{k,i}\}_{i=1}^{\dim(\cH^M(k))}$ is an orthonormal basis of $\cH^M(k)$ and $y^\pi_T\in [\cH^M(N)]^n$, it follows from \eqref{w23e14} that
\bel{w23e16}
\bal
\sum_{i=1}^{\dim(\cH^M(k))} \me\big[ h_{k,i} y^\pi_T\big]h_{k,i}
= \sum_{i=1}^{\dim(\cH^M(k))} \me\big[ h_{k,i}\me(y^\pi_T\,|\,\cF_{t_k})\big]h_{k,i}
= \me(y^\pi_T\,|\,\cF_{t_k})\,.
\eal
\ee

We next consider the generator term. Set $F_\t(\cd) = f\big( \nu(\cd),y_{N,M}(\cd),Y_{N,M-1}(\cd)\big)$. Since $h_{k,i}\in\cH^M(k)\subseteq\cH^M(N)$ and $\G_M$ is the componentwise $L^2$-orthogonal projection onto $\bigoplus_{m=0}^M(\cH_m)^n$, we have $\me\big[ h_{k,i}F_\tau(\cd)\big]=\me\big[h_{k,i}\G_MF_\tau(\cd)\big]$.
Under the present setting, for $t\in[0,T)$, $F_\tau(t)$ is measurable with respect to $\cG_{\nu(t)}$. Hence $\G_MF_\tau(t)=\G_M^NF_\tau(t)$ and $\G_MF_\tau(t)\in[\cH^M(N)]^n$.
Consequently,
\beq
\int_{t_k}^T\G_MF_\tau(t)\rd t
\in[\cH^M(N)]^n\,.
\eeq
Applying again Remark \ref{w408r2} and the componentwise orthonormal expansion with respect to $\{h_{k,i}\}$ gives
\begin{align}
&\sum_{k=0}^{N-1}\sum_{i=1}^{\dim(\cH^M(k))}\me\Big[\int_{t_k}^T\sqrt{\tau } h_{k,i} F_\tau(t) \rd t \Big]e_{k,i}(\cd)
=\sum_{k=0}^{N-1}\chi_{[t_k,t_{k+1})}(\cd)\sum_{i=1}^{\dim(\cH^M(k))}\me\Big[h_{k,i} \int_{t_k}^T \G_M F_\tau(t) \rd t \Big] h_{k,i}\notag \\
&\qq\qq\q=\sum_{k=0}^{N-1}\chi_{[t_k,t_{k+1})}(\cd)\me\Big(\int_{t_k}^T  \G_M F_\tau(t) \rd t \,\Big|\,\cF_{t_k}\Big)\,.\label{w23e17}
\end{align}

Combining \eqref{w23e13}, \eqref{w23e16} and \eqref{w23e17}, we have
\beq
\bal
y_{N,M}(\cd)=\sum_{k=0}^{N-1}\chi_{[t_k,t_{k+1})}(\cd)\Big[\me(y^\pi_T\,|\,\cF_{t_k})-\me\Big(\int_{t_k}^T
\G_M F_\tau(t) \rd t \,\Big|\,\cF_{t_k}\Big)\Big]\,,
\eal
\eeq
which proves the first assertion of \eqref{w23e10a}.

\ss\no
{\bf Step 1(b).}
%
We next derive the representation of $Y_{N,M-1}$. Since $y^\pi_T\in[\cH^M(N)]^n$, by Remark \ref{w9r1} it admits the finite Wiener chaos expansion
\beq
y^\pi_T=\sum_{m=0}^M\sum_{\substack{|\a|=m\\ \a\in\Lambda_d(N)}}d^\a \sqrt{\a!} \prod_{j=1}^{N}\prod_{\ell=1}^{d} H_{\a_j^\ell}\big(\dbW(g_{j,\ell})\big)\,,
\eeq
where the coefficients $d^\a=\me\big[\sqrt{\a!} \prod_{j=1}^{N}\prod_{\ell=1}^{d} H_{\a_j^\ell}\big(\dbW(g_{j,\ell})\big) y^\pi_T\big]$.
Using the following recurrence relation for Hermite polynomials (see, e.g., \cite[Page 5]{Nualart06}),
\beq
xH_m(x)=(m+1)H_{m+1}(x)+H_{m-1}(x)\qq \forall\,m\in\dbN_0\,,
\eeq
we have, for $1\leq r\leq d$,
\begin{align*}
\frac{\D_{k+1}W^r}{\sqrt{\t}} y^\pi_T
&=\sum_{m=0}^M\sum_{\substack{|\a|=m\\ \a\in\Lambda_d(N)}}d^\a \sqrt{\a!} \prod_{j=1}^{N}\prod_{\ell=1}^{d} H_{\a_j^\ell}\big(\dbW(g_{j,\ell})\big)H_1\big(\dbW(g_{k+1,r})\big)\\
&=\sum_{m=0}^M\sum_{\substack{|\a|=m\\ \a\in\Lambda_d(N)}}d^\a \sqrt{\a!} \prod_{\substack{1\leq j\leq N,\ 1\leq\ell\leq d\\
(j,\ell)\neq(k+1,r)}} H_{\a_j^\ell}\big(\dbW(g_{j,\ell})\big) \\
&\q\times
\Big[(\a_{k+1}^r+1)H_{\a_{k+1}^r+1}\big(\dbW(g_{k+1,r})\big) + H_{\a_{k+1}^r-1}\big(\dbW(g_{k+1,r})\big)\Big]\,.
\end{align*}
Taking the conditional expectation with respect to $\cF_{t_k}$, all terms containing a Hermite polynomial of positive degree in any increment after $t_k$ vanish. In particular, the term containing $H_{\a_{k+1}^r+1}\big(\dbW(g_{k+1,r})\big)$ always has zero conditional expectation. The second term can survive only if
$\a_{k+1}^r=1$ and $\a_j^\ell=0$ for all $j\geq k+1$ as well as $(j,\ell)\neq(k+1,r)$.
Therefore,
\begin{align*}
\me\Big(\frac{\D_{k+1}W^r}{\sqrt{\t}} y^\pi_T\,\big|\,\cF_{t_k}\Big) 
= \sum_{m=1}^M\sum_{\substack{
|\a|=m,\, \a\in\Lambda_d(N),\,\a_{k+1}^r=1\\
\a_j^\ell=0,\, j\geq k+1,\, (j,\ell)\neq(k+1,r)
}} d^\alpha \sqrt{\alpha!}  \prod_{j=1}^k\prod_{\ell=1}^d H_{\alpha_j^\ell}\big(\dbW(g_{j,\ell})\big)
\in [\cH^{M-1}(k)]^n\,.
\end{align*}

Since $r\in\{1,\cds,d\}$ is arbitrary, collecting the $d$ columns gives
\bel{w23el3c}
\me\Big(
\frac{1}{\sqrt{\tau}}
y^\pi_T(\D_{k+1}W)^\top
\,\Big|\,
\cF_{t_k}
\Big)
\in
[\cH^{M-1}(k)]^{n\times d}\,.
\ee
Applying the orthonormal expansion componentwise to the matrix-valued random variable in \eqref{w23el3c}, we have
\begin{align}
\sum_{i=1}^{\dim(\cH^{M-1}(k))}\me\Big[ \frac{h_{k,i}}{\sqrt{\tau }}  y^\pi_T(\D_{k+1}W)^\top\Big] h_{k,i}
&=\sum_{i=1}^{\dim(\cH^{M-1}(k))}\me\Big[ h_{k,i}\me\Big( \frac{1}{\sqrt{\tau }} y^\pi_T(\D_{k+1}W)^\top\,\Big|\,\cF_{t_k}\Big)\Big] h_{k,i}\notag\\
&=\me\Big( \frac{1}{\sqrt{\tau }} y^\pi_T(\D_{k+1}W)^\top\,\Big|\,\cF_{t_k}\Big)\,.\label{w23el1}
\end{align}

A similar argument, together with $\D_{k+1}W^r h_{k,i}\in\cH^M(k+1)\subseteq\cH^M(N)$ for $1\leq r\leq d$ and $i=1,\cds,\dim(\cH^{M-1}(k))$, yields
\begin{align}
&\sum_{k=0}^{N-1}\sum_{i=1}^{\dim(\cH^{M-1}(k))} \me\Big[\int_{t_{k+1}}^T \frac{h_{k, i}}{\sqrt{\tau }}  F_\tau(t)(\D_{k+1}W)^\top \rd t\Big]\, e_{k,i}(\cd)\notag\\
&=\sum_{k=0}^{N-1}\chi_{[t_k,t_{k+1})}(\cd)\sum_{i=1}^{\dim(\cH^{M-1}(k))} \me \Big[h_{k,i} \int_{t_{k+1}}^T \frac{1}{\tau }  \G_M F_\tau(t) \rd t (\D_{k+1}W)^\top \Big]\, h_{k,i} \notag\\
&=\sum_{k=0}^{N-1}\chi_{[t_k,t_{k+1})}(\cd) \me \Big(\int_{t_{k+1}}^T \frac{1}{\tau }  \G_M F_\tau(t) \rd t\, (\D_{k+1}W)^\top\,\Big|\,\cF_{t_k}\Big)\,.\label{w23el50}
\end{align}
%
Combining \eqref{w23el1}, \eqref{w23el50}, \eqref{w23e1} and \eqref{w23e9}, we obtain the representation of $Y_{N,M-1}(\cd)$
\begin{align*}
&Y_{N,M-1}(\cd)
=\sum_{k=0}^{N-1}\sum_{i=1}^{\dim(\cH^{M-1}(k))}b_{k,i}e_{k,i}(\cd)\\
&\q=\sum_{k=0}^{N-1}\chi_{[t_k,t_{k+1})}(\cd)\me\Big( \frac{1}{\tau }\Big[ y^\pi_T- \int_{t_{k+1}}^T \G_M F_\tau(t)\rd t\Big](\D_{k+1}W)^\top \Big|\,\cF_{t_k}\Big)\,,
\end{align*}
which is the second assertion of \eqref{w23e10a}.

\ss
\no
{\bf Step 2.} In this step, we establish the following error estimate
\begin{align}
&\max_{0\leq k\leq N} \me\big[|y_{N,M}(t_k)-y(t_k)|^2\big]+\me\int_0^T|Y_{N,M-1}(t)-Y(t)|^2\rd t \notag\\
&\leq \cC\Big\{\tau \!+\!\me\big[|y_T\!-\!y^\pi_T|^2\big]\!+\!\me \!\int_0^T\!\! \big|f\big(s,y(s),Y(s)\big)\!-\!\G_M f\big(s,y(s),Y(s)\big)\big|^2 \rd s\Big\}\,. \label{w1102e2}
\end{align}
To prove \eqref{w1102e2}, we adapt the standard Euler error argument, using the conditional-average approximation \eqref{w911bar} for the martingale integrand.

\ss\no
{\bf Step 2(a).} We introduce an auxiliary pair $\big(y_0(\cd),Y_0(\cd)\big)$ as the adapted solution of the following BSDE:
\bel{w0310e1}
\lt\{
\bal
&y_0(t_{k+1})-y_0(t)=\int_t^{t_{k+1}} \G_M f\big(t_k,y_0(t_k), \bar Y_0(t_k)\big)\rd s\\
&\qq\qq\q +\int_t^{t_{k+1}} Y_0(s) \rd W(s)\,, \q
    t\in [t_k,t_{k+1}),\,k=0,1,\cds,N-1\,,\\
&y_0(T)=y^\pi_T\,,
\eal
\rt.
\ee
where
\bel{w0311e2}
\bar Y_0(t_k)=\frac {1}{\tau } \me\Big(\int_{t_k}^{t_{k+1}}Y_0(s)\rd s\,\Big|\, \cF_{t_k}\Big)\,.
\ee
This auxiliary process serves as a bridge between the adapted solution $\big(y(\cdot),Y(\cdot)\big)$ and the finite transposition solution $\big(y_{N,M}(\cdot),Y_{N,M-1}(\cdot)\big)$.
A standard fixed-point argument shows that, for sufficiently small $\tau$, equation \eqref{w0310e1} admits a unique adapted solution $\big(y_0(\cd),Y_0(\cd)\big)$.
%

Now we estimate the difference between $\big(y_0(\cd),Y_0(\cd)\big)$ and $\big(y(\cd),Y(\cd)\big)$. By setting
$
F(\cd)=f\big(\cd,y(\cd),Y(\cd)\big)
$
and applying BSDEs \eqref{bsde0} and \eqref{w0310e1}, we have
\begin{align*}
&y(t_k)-y_0(t_k)+\int_{t_k}^{t_{k+1}} \big[Y(s)-Y_0(s)\big]\rd W(s)\\
&= y(t_{k+1})-y_0(t_{k+1}) -\int_{t_k}^{t_{k+1}}\big[F(s)-
     \G_M f\big(t_k,y_0(t_k), \bar Y_0(t_k)\big)\big]\rd s\,.
\end{align*}
Using the orthogonality of the stochastic integral to $\cF_{t_k}$-measurable random variables, we obtain
\begin{align}
&\me\big[|y(t_k)-y_0(t_k)|^2\big]+\me\int_{t_k}^{t_{k+1}}|Y(s)-Y_0(s)|^2\rd s \notag\\
&=\me\bigg\{\Big| y(t_{k+1})-y_0(t_{k+1})
    -\int_{t_k}^{t_{k+1}}\![F(s)-\G_M F(s)]\rd s 
  -\int_{t_k}^{t_{k+1}}\!\big[\G_M F(s)-
     \G_M f\big(t_k,y_0(t_k), \bar Y_0(t_k)\big) \big]\rd s\Big|^2\bigg\} \notag \\
&\leq   \big(1+6\frac{\tau }{\e}\big)\me\big[|y(t_{k+1})-y_0(t_{k+1})|^2\big]
   +2\big(1+\frac{\e}{6\tau }\big)\tau  \me\Big[\int_{t_k}^{t_{k+1}}|F(s)-\G_M F(s)|^2\rd s\Big] \notag \\
&\qq\qq+2\big(1+\frac{\e}{6\tau }\big)\tau  \me\Big[\int_{t_k}^{t_{k+1}}\big|F(s)-
      f\big(t_k,y_0(t_k), \bar Y_0(t_k)\big) \big|^2\rd s\Big]\,.\label{w0310e3}
\end{align}
Here we have used the Cauchy--Schwarz inequality, the contraction property of the orthogonal projection $\G_M$, and Young's inequality.
Since
\begin{align*}
&\big|F(s)-
f\big(t_k,y_0(t_k), \bar Y_0(t_k)\big)\big|\\
&\qq\leq  L\Big\{|s-t_k|^{1/2}+|y(s)-y(t_k)|+|y(t_k)-y_0(t_k)|+|Y(s)-\bar Y(t_k)|+|\bar Y(t_k)-\bar Y_0(t_k)|\big\} \,,
\end{align*}
where $\bar Y(t_k)$ is defined in \eqref{w911bar},
and
\beq
|\bar Y(t_k)-\bar Y_0(t_k)|
\leq \frac{1}{\tau} \me\Big( \int_{t_k}^{t_{k+1}}|Y(s)-Y_0(s)|\rd s \,\Big|\,\cF_{t_k}\Big)\,,
\eeq
it follows that
\begin{align*}
&\me\big[|y(t_k)-y_0(t_k)|^2\big]+\me\int_{t_k}^{t_{k+1}}|Y(s)-Y_0(s)|^2\rd s\\
&\leq\big(1+6\frac{\tau }{\e}\big)\me\big[|y(t_{k+1})-y_0(t_{k+1})|^2\big]
    +2\big(1+\frac{\e}{6\tau }\big)\tau  \me \int_{t_k}^{t_{k+1}}|F(s)-\G_M F(s)|^2\rd s\\
  &\q+2\big(1+\frac{\e}{6\tau }\big)\tau  6L^2 \Big\{ \frac{\tau ^2}{2}
    +\me\int_{t_k}^{t_{k+1}}\big[ |y(s)-y(t_k)|^2
      + |Y(s)-\bar Y(t_k)|^2\big] \rd s\\
  &\q\qq +\tau \me\big[|y(t_k)-y_0(t_k)|^2\big]+\me\int_{t_k}^{t_{k+1}}  |Y(s)-Y_0(s)|^2 \rd s \Big\} \,.
\end{align*}
Consequently,
\begin{align}
&\Big[1-2\big(1+\frac{\e}{6\tau }\big) 6L^2\tau ^2 \Big]\me\big[|y(t_k)-y_0(t_k)|^2\big]
+\Big[1-2\big(1+\frac{\e}{6\tau }\big) 6L^2\tau  \Big]\me \int_{t_k}^{t_{k+1}}|Y(s)-Y_0(s)|^2\rd s \notag \\
&\leq\big(1\!+\!6\frac{\tau }{\e}\big)\me\big[|y(t_{k+1})\!-\!y_0(t_{k+1})|^2\big]
    \!+\!2\big(1\!+\!\frac{\e}{6\tau }\big)\tau  \me \!\int_{t_k}^{t_{k+1}}\!\!|F(s)\!-\!\G_M F(s)|^2\rd s \notag \\
    &\qq+\!2\big(1\!+\!\frac{\e}{6\tau }\big)\tau  6L^2 \Big\{ \frac{\tau ^2}{2}
    \!+\!\me\!\int_{t_k}^{t_{k+1}}\!\big[ |y(s)\!-\!y(t_k)|^2
      \!+\! |Y(s)\!-\!\bar Y(t_k)|^2\big] \rd s\Big\} \,.\label{w0310e6}
\end{align}
Set $\e=\frac{1-36L^2\tau }{6L^2}$.
Since $N\geq N_*$, we have $ \tau \leq 1\wedge\frac{1}{72L^2}$. Then
\beq
 1-2\big(1+\frac{\e}{6\tau }\big) 6L^2\tau  =\frac 23\,, \qq
 1-2\big(1+\frac{\e}{6\tau }\big) 6L^2\tau ^2=1-\frac 13\tau \,,\\
\eeq
and
\begin{align*}
&\frac{1+6\frac{\tau }{\e}}{1-2\big(1+\frac{\e}{6\tau }\big) 6L^2\tau ^2}
    \leq 1+2\big(72L^2+\frac 13 \big)\tau\,,\\
&\frac{1-2\big(1+\frac{\e}{6\tau }\big) 6L^2\tau }{1-2\big(1+\frac{\e}{6\tau }\big) 6L^2\tau ^2}
    \in \big[\frac 12, 1\big]\,, \qq \frac{2\big(1+\frac{\e}{6\tau }\big)\tau  6L^2}{1-2\big(1+\frac{\e}{6\tau }\big) 6L^2\tau ^2}
=\frac{1}{3-\tau }\leq \frac 23\,.
\end{align*}
Therefore, \eqref{w0310e6} leads to
\begin{align}
&\me\big[|y(t_k)-y_0(t_k)|^2\big]+\frac 12 \me\int_{t_k}^{t_{k+1}}|Y(s)-Y_0(s)|^2\rd s \notag\\
&\leq (1+\cC \tau )\me\big[|y(t_{k+1})-y_0(t_{k+1})|^2\big]
 + \frac 1 {12L^2}
  \me \int_{t_k}^{t_{k+1}}|F(s)-\G_M F(s)|^2\rd s \notag \\
 &\q+\frac 23 \Big\{ \frac{\tau ^2}{2}
    +\me\int_{t_k}^{t_{k+1}}\big[ |y(s)-y(t_k)|^2
      + |Y(s)-\bar Y(t_k)|^2\big] \rd s \Big\}\,. \label{w0310e9}
\end{align}
Applying the discrete Gronwall inequality and Assumption {\rm\bf(H)$_M$} (in fact, only {\rm\bf(H)$_1$} is used here; see Remark \ref{w1102r1}), we obtain
\begin{align}
\max_{0\leq k\leq N-1}\me\big[|y(t_k)-y_0(t_k)|^2\big] 
&\leq e^{\cC T}\bigg\{ \me\big[|y(T)-y_0(T)|^2\big]
  + \frac 1 {12L^2}
   \sum_{k=0}^{N-1}\me\int_{t_k}^{t_{k+1}}|F(s)-\G_M F(s)|^2\rd s \notag \\
 &\qq\q+\frac T3\tau +\frac 23 \sum_{k=0}^{N-1}\me \int_{t_k}^{t_{k+1}}\big[|y(s)-y(t_k)|^2+|Y(s)-\bar Y(t_k)|^2\big]\rd s
 \bigg\} \notag \\
&\leq  e^{\cC T}\Big\{ \me\big[|y(T)-y_0(T)|^2\big]+\cC \me \int_0^T |F(s)-\G_M F(s)|^2 \rd s+\cC \tau \Big\}\,.\label{w0310e10}
\end{align}

Summing \eqref{w0310e9} over $k=0,\ldots,N-1$, we obtain
\begin{align}
&\me\big[|y(t_0)-y_0(t_0)|^2\big]+\frac 1 2 \me \int_{0}^{T}|Y(s)-Y_0(s)|^2\rd s \notag \\
&\leq \cC \tau  \sum_{k=0}^{N-1}\me\big[|y(t_{k+1})-y_0(t_{k+1})|^2\big]+ (1+\cC \tau )\me\big[|y(t_{N})-y_0(t_N)|^2\big]
+\cC
  \me \int_{0}^{T}|F(s)-\G_M F(s)|^2\rd s +\cC \tau \notag \\
&\leq  \cC \Big\{\tau +\me\big[|y(T)-y_0(T)|^2\big] +\me\int_0^T |F(s)-\G_M F(s)|^2 \rd s\Big\}\,.\label{w0311e1}
\end{align}

\ss
\no
{\bf Step 2(b).}
Comparing the conditional-expectation representation of \eqref{w0310e1} with \eqref{w23e10a}, we obtain, for $k=0,1,\cds,N-1$,
\bel{w0311e3}
 y_0(t_k)=y_{N,M}(t_k)\,,\q \bar Y_0(t_k)=Y_{N,M-1}(t_k)\,.
\ee

We now estimate the error between $Y(\cd)$ and $ Y_{N,M-1}(\cd)$.
By the triangle inequality, the $L^2$-best approximation property of $\bar Y_0(t_k)$, \eqref{w709e2}, and \eqref{w0311e1}, we obtain
\begin{align}
&\sum_{k=0}^{N-1}\me \int_{t_k}^{t_{k+1}} |Y(s)-Y_{N,M-1}(t_k)|^2\rd s
\leq 2 \sum_{k=0}^{N-1}\me\int_{t_k}^{t_{k+1}} \big[|Y(s)-Y_0(s)|^2+|Y_0(s)-\bar Y_0(t_k)|^2\big]\rd s\notag \\
&\qq\qq\leq 2  \sum_{k=0}^{N-1}\me\int_{t_k}^{t_{k+1}} \big[3|Y(s)-Y_0(s)|^2+2|Y(s)-\bar Y(t_k)|^2\big]\rd s\notag \\
&\qq\qq\leq  \cC \Big\{\tau +\me\big[|y(T)-y_0(T)|^2\big] +\me\int_0^T |F(s)-\G_M F(s)|^2 \rd s\Big\}\,.\label{w0311e4}
\end{align}
Combining \eqref{w0310e10}--\eqref{w0311e4}, we arrive at \eqref{w1102e2}.

\ss\no
{\bf Step 3.}
Assumption {\rm\bf(H)$_M$} implies, in particular, that $y_T\in\dbD^{1,2}(\dbR^n)$ and that $f(\cd,\cd,\cd)$ is continuously differentiable with respect to its second and third variables, with uniformly bounded derivatives.
By \cite[Proposition 5.3]{ElKaroui-Peng-Quenez97}, $y(t)\in\dbD^{1,2}(\dbR^n)$, $Y(t)\in\dbD^{1,2}(\dbR^{n\times d})$ for almost every $t\in[0,T]$, and $\sum_{\ell=1}^d \me\int_0^T\int_0^T\big[|D_\theta^\ell y(t)|^2+|D_\theta^\ell Y(t)|^2\big]\rd\theta\rd t<\infty$.
Consequently,
\begin{align}
\me\int_0^T\big|F(s)-\G_MF(s)\big|^2\rd s
\leq\frac{1}{M+1} \sum_{\ell=1}^d\int_0^T\int_0^T\me\big[|D_\th^\ell F(s)|^2\big]\rd\th\rd s\,.
\label{w1227e6a}
\end{align}

For each $\ell=1,\cds,d$, the Malliavin chain rule yields
\beq
\bal
D_\th^\ell F(s)
&=\pr_yf\big(s,y(s),Y(s)\big)D_\th^\ell y(s)
+\pr_Yf\big(s,y(s),Y(s)\big)\big[D_\th^\ell Y(s)\big]\,,
\eal
\eeq
where $\pr_yf\big(s,y(s),Y(s)\big)\in\dbR^{n\times n}$, while $\pr_Yf\big(s,y(s),Y(s)\big): \dbR^{n\times d}\longrightarrow\dbR^n$ denotes the Fr\'echet derivative of $f(\cd,\cd,\cd)$ with respect to its third variable. Since these derivatives are uniformly bounded, we have
\beq
|D_\th^\ell F(s)|^2
\leq \cC \big(|D_\th^\ell y(s)|^2+|D_\th^\ell Y(s)|^2\big)\,.
\eeq
Consequently, \eqref{w1227e6a} gives
\begin{align}
\me\int_0^T\big|F(s)-\G_MF(s)\big|^2\rd s
\leq \frac{\cC}{M+1} \sum_{\ell=1}^d\int_0^T\int_0^T\me\big[|D_\th^\ell y(s)|^2+|D_\th^\ell Y(s)|^2\big]\rd\th\rd s
\leq \frac{\cC}{M+1}\,.
\label{w1227e6}
\end{align}

We next estimate $\me\big[|y_T-y^\pi_T|^2\big]$.
Since $y^\pi_T=\G_M^Ny_T$, we have
\bel{w1227e2}
\me\big[|y_T-\G_M^Ny_T|^2\big]\leq 2\me\big[|y_T-\G_M y_T|^2\big]+2\me\big[|\G_My_T-\G_M^Ny_T|^2\big]\,.
\ee
We first estimate the Wiener chaos truncation error.
Since $y_T\in \dbD^{M,2}(\dbR^n)$, an argument similar to that used in \eqref{w1227e6} yields
\begin{align}
\me\big[|y_T-\G_M y_T|^2\big]
&\leq \frac 1{(M+1)M\cds(M-r+2)}
\sum_{\ell_1,\cds,\ell_r=1}^d \me\int_{(0,T)^r}
\big|D_{\th_1,\cds,\th_r}^{\ell_1,\cds,\ell_r}y_T\big|^2\rd\th_1\cds\rd\th_r\notag\\
&\leq \frac {\|y_T\|^2_{r,2}}{(M+1)M\cds(M-r+2)}\,,\qq \forall\, r=1,2,\cds,M\,.\label{w1227e3}
\end{align}
We now estimate the error caused by the discretization of the Wiener chaos kernels. By the Wiener chaos expansion \eqref{ex-y},
$\G_My_T=\me[y_T]+\sum_{m=1}^M J_m(h_m)$,
where $h_m=\big\{h_m^{\ell_1,\cds,\ell_m}\big\}_{1\leq\ell_1,\cds,\ell_m\leq d}\in L^2_{S,d}(T^m;\dbR^n)$.
On the other hand, Lemma \ref{w524l1} gives
$\G_M^Ny_T=\me[y_T]+\sum_{m=1}^M J_m(\bar h_m)$,
where $\bar h_m=\big\{\bar h_m^{\ell_1,\cds,\ell_m}\big\}_{1\leq\ell_1,\cds,\ell_m\leq d}\in L^2_{S,d}(T^m;\dbR^n)$
is defined by \eqref{w524e7}.
Since Wiener chaoses of different orders are mutually orthogonal, the It\^o isometry for iterated stochastic integrals yields
\begin{align}
&\me\big[|\G_My_T-\G_M^Ny_T|^2\big]
=\sum_{m=1}^M\me\big[ |J_m(h_m-\bar h_m)|^2\big]
\notag\\
&=
\sum_{m=1}^M
\sum_{\ell_1,\cds,\ell_m=1}^d
\int_0^T
\int_0^{s_m}
\cds
\int_0^{s_2}
\Big|
h_m^{\ell_1,\cds,\ell_m}(s_1,\cds,s_m)
 - \bar h_m^{\ell_1,\cds,\ell_m}(s_1,\cds,s_m)
\Big|^2
\rd s_1\cds\rd s_m\,.
\label{w1227e4a}
\end{align}

We next estimate the integrand in \eqref{w1227e4a}. Fix $\ell_1,\cds,\ell_m\in\{1,\cds,d\}$ and $(s_1,\cds,s_m)\in I_{\Bj_m}$. By \eqref{w524e7},
\begin{align*}
&h_m^{\ell_1,\cds,\ell_m}(s_1,\cds,s_m)
-\bar h_m^{\ell_1,\cds,\ell_m}(s_1,\cds,s_m)\\
&=\frac{1}{|I_{\Bj_m}|} \int_{I_{\Bj_m}}
\big[h_m^{\ell_1,\cds,\ell_m}(s_1,\cds,s_m)
-h_m^{\ell_1,\cds,\ell_m}(r_1,\cds,r_m)\big]
\rd r_1\cds\rd r_m\,.
\end{align*}
By condition \eqref{lip-h-2}, for $(s_1,\cds,s_m), (r_1,\cds,r_m)\in I_{\Bj_m}$,
\beq
\big|h_m^{\ell_1,\cds,\ell_m}(s_1,\cds,s_m)-h_m^{\ell_1,\cds,\ell_m}(r_1,\cds,r_m)\big|
\leq L^{1/2}\sum_{j=1}^m |s_j-r_j|^{1/2}
\leq L^{1/2}m\t^{1/2}\,,
\eeq
which, along with H\"older's inequality, gives
\begin{align*}
&\Big|
h_m^{\ell_1,\cds,\ell_m}(s_1,\cds,s_m)
-
\bar h_m^{\ell_1,\cds,\ell_m}(s_1,\cds,s_m)
\Big|^2
\notag\\
&\leq
\frac{1}{|I_{\Bj_m}|}
\int_{I_{\Bj_m}}
\big|
h_m^{\ell_1,\cds,\ell_m}(s_1,\cds,s_m)
-
h_m^{\ell_1,\cds,\ell_m}(r_1,\cds,r_m)
\big|^2
\rd r_1\cds\rd r_m\notag\\
&\leq Lm^2\t\,.
\end{align*}
Substituting this estimate into \eqref{w1227e4a}, we obtain
\begin{align}
\me\big[|\G_My_T-\G_M^Ny_T|^2\big]
\leq L\tau \sum_{m=1}^\infty \frac{d^m m^2T^m}{m!}
\leq \cC\tau\,.
\label{w1227e4}
\end{align}
Here $\cC$ is independent of $M$ and $N$.

Finally, \eqref{w1102e2}, together with estimates \eqref{w1227e6} through \eqref{w1227e4}, leads to the assertion \eqref{w416e1}.
This completes the proof.
\end{proof}

\br{w506r1}
For $y_T\in\dbD^{m,2}(\dbR^n)$, the deterministic kernels in the Wiener chaos expansion \eqref{ex-y} can be represented, componentwise with respect to the Brownian coordinates, by expectations of the iterated Malliavin derivatives. More precisely, for any $k=1,2,\cds,m$ and any $\ell_1,\cds,\ell_k\in\{1,\cds,d\}$,
$h_k^{\ell_1,\cds,\ell_k}(s_1,\cds,s_k) =\me\big[D_{s_1,\cds,s_k}^{\ell_1,\cds,\ell_k}y_T\big]$
for a.e. $(s_1,\cds,s_k)\in[0,T]^k$, where
$h_k=\big\{h_k^{\ell_1,\cds,\ell_k}\big\}_{1\leq\ell_1,\cds,\ell_k\leq d}\in L^2_{S,d}(T^k;\dbR^n)$.
Thus, condition \eqref{lip-h-2} may equivalently be formulated as a $1/2$-H\"older continuity condition on the kernels
$h_k^{\ell_1,\cds,\ell_k}$, uniformly over
$(\ell_1,\cds,\ell_k)\in\{1,\cds,d\}^k$.

For the convergence rate with respect to a fixed partition $I_\t$, it is sufficient to impose the weaker gridwise condition
\beq
\bal
\big|h_k^{\ell_1,\cds,\ell_k}(s_1,\cds,s_k)-h_k^{\ell_1,\cds,\ell_k}\big(\nu(s_1),\cds,\nu(s_k)\big)\big|
\leq\cC k\tau^{1/2}\,,
\eal
\eeq
for all $\ell_1,\cds,\ell_k\in\{1,\cds,d\}$,
$0\leq s_1\leq\cds\leq s_k\leq T$, and
$k=1,2,\cds,m$.
See \eqref{w1227e4} in the proof of Theorem \ref{ftm-rate}.
\er

\br{w630r1}
Suppose that Assumption {\rm\bf(H)$_M$} holds with $q>8$ and $M\geq2$.
Then the adapted solution $\big(y(\cd),Y(\cd)\big)$ satisfies
\beq
\bal
\me\int_0^T
\bigg[
&\Big(
\sum_{\ell=1}^d
\int_0^T
|D_\th^\ell y(t)|^2
\rd\th
\Big)^2
+\Big(\sum_{\ell=1}^d
\int_0^T
|D_\th^\ell Y(t)|^2
\rd\th\Big)^2\\
&\q+
\sum_{\ell_1,\ell_2=1}^d
\int_0^T\int_0^T
\big(\big|
D_{\th_1,\th_2}^{\ell_1,\ell_2}y(t)
\big|^2+\big|
D_{\th_1,\th_2}^{\ell_1,\ell_2}Y(t)
\big|^2\big)
\rd\th_1\rd\th_2
\bigg]\rd t
<\infty\,.
\eal
\eeq
Here
$D_\th^\ell y(t)\in\dbR^n$,
$D_\th^\ell Y(t)\in\dbR^{n\times d}$,
$D_{\th_1,\th_2}^{\ell_1,\ell_2}y(t)\in\dbR^n$, and
$D_{\th_1,\th_2}^{\ell_1,\ell_2}Y(t)\in\dbR^{n\times d}$.

Since $f(\cd,\cd,\cd)$ has uniformly bounded first- and second-order derivatives with respect to its second and third variables, the Malliavin chain rule implies that $f\big(\cd,y(\cd),Y(\cd)\big)\in L^2\big(0,T;\dbD^{2,2}(\dbR^n)\big)$.
Indeed, for any $\ell\in\{1,\cds,d\}$,
\beq
D_\th^\ell
f\big(t,y(t),Y(t)\big)
=\pr_yf\big(t,y(t),Y(t)\big)D_\th^\ell y(t)
+\pr_Yf\big(t,y(t),Y(t)\big)
\big[D_\th^\ell Y(t)\big]\,,
\eeq
and, for any $\ell_1,\ell_2\in\{1,\cds,d\}$,
\begin{align*}
&D_{\th_1,\th_2}^{\ell_1,\ell_2}
f\big(t,y(t),Y(t)\big)\\
&=
\pr_yf\big(t,y(t),Y(t)\big)
D_{\th_1,\th_2}^{\ell_1,\ell_2}y(t)
+\pr_Yf\big(t,y(t),Y(t)\big)
\big[D_{\th_1,\th_2}^{\ell_1,\ell_2}Y(t)\big]\\
&\q+\pr_{yy}f\big(t,y(t),Y(t)\big)
\big[
D_{\th_1}^{\ell_1}y(t),
D_{\th_2}^{\ell_2}y(t)
\big]
+\pr_{yY}f\big(t,y(t),Y(t)\big)
\big[
D_{\th_1}^{\ell_1}y(t),
D_{\th_2}^{\ell_2}Y(t)
\big]\\
&\q+\pr_{Yy}f\big(t,y(t),Y(t)\big)
\big[
D_{\th_1}^{\ell_1}Y(t),
D_{\th_2}^{\ell_2}y(t)
\big]
+\pr_{YY}f\big(t,y(t),Y(t)\big)
\big[
D_{\th_1}^{\ell_1}Y(t),
D_{\th_2}^{\ell_2}Y(t)
\big]\,.
\end{align*}
Therefore, by the second-order Wiener chaos truncation estimate,
\beq
\bal
\me\int_0^T
\big|f\big(s,y(s),Y(s)\big)
-\G_Mf\big(s,y(s),Y(s)\big)
\big|^2\rd s
\leq \frac{\cC}{M(M+1)}\,.
\eal
\eeq
Consequently, the convergence rate in \eqref{w416e1} can be improved to
\beq
\sup_{0\leq t\leq T}\me\big[\big|y_{N,M}(t)-y\big(\nu(t)\big)\big|^2\big]
+\me\int_0^T|Y_{N,M-1}(t)-Y(t)|^2\rd t
\leq \cC\Big(\frac{1}{N}
+\frac{1}{M(M+1)}\Big)\,.
\eeq
\er

\br{w0420r2}
{\rm (1)} The representation \eqref{w23e10a} in Theorem \ref{ftm-rate} resembles the Euler scheme for the projected BSDE
\bel{bsde2}
\left\{
\bal
&\mathrm{d} y(t)=\G_Mf\big(t,y(t),Y(t)\big)\rd t+Y(t)\rd W(t)\,,\quad t\in [0,T]\,,\\
&y(T)=y_T\,.
\eal
\right.
\ee
See, e.g., \cite{Bender-Denk07, Bouchard-Touzi04, ZhangJF04} for Euler schemes for BSDEs.

Nevertheless, the finite transposition method should not be identified with a direct application of the Euler method to \eqref{bsde2}. In particular, the $L^2$-regularity estimate for its martingale integrand corresponding to \eqref{w709e2} does not follow directly from the usual regularity theory for BSDEs. Such an estimate is essential in the standard error analysis of Euler schemes.

By contrast, the finite transposition method is defined by the finite-dimensional variational identity \eqref{w23e2}. Its algebraic formulation \eqref{w23e1}, \eqref{w23e6}, and \eqref{w23e9} determines the coefficients of the numerical solution without requiring the direct evaluation of conditional expectations. The conditional-expectation formulas in \eqref{w23e10a} are derived from the finite transposition solution and serve only as representations of that solution; they do not define the numerical scheme itself.


%
%

{\rm (2)} When BSDE \eqref{bsde0} is linear, under Assumptions {\rm\bf (P)} and {\rm\bf (B)}, the representation \eqref{w23e10a} reduces to
\beq
\lt\{
\bal
&y_{N,M}(\cd)=\sum_{k=0}^{N-1}\chi_{[t_k,t_{k+1})}(\cd)\me\Big(y^\pi_T-\int_{t_k}^T  f\big(\nu(t),y_{N,M}(t),Y_{N,M-1}(t)\big) \rd t \,\Big|\,\cF_{t_k}\Big)\,,\\
&Y_{N,M-1}(\cd)=\sum_{k=0}^{N-1}\chi_{[t_k,t_{k+1})}(\cd)\me\Big( \frac{1}{\tau }\Big[ y^\pi_T
- \int_{t_{k+1}}^T f\big(\nu(t),y_{N,M}(t),Y_{N,M-1}(t)\big)\rd t\Big] (\D_{k+1}W)^\top \,\Big|\,\cF_{t_k}\Big)\,.
\eal
\rt.
\eeq
\er

\br{w708r1}
%
The representation \eqref{w23e1} should be understood as a coefficient-based, rather than a pathwise, representation.
More precisely, the unknowns in \eqref{w23e1} are the deterministic coefficients $\{a_{k,i}\}_{k,i}$ and $\{b_{k,i}\}_{k,i}$, whereas the basis functions $\{e_{k,i}\}_{k,i}$ are random variables depending on the Brownian increments.
Consequently, once these coefficients have been determined, the approximate solution along any realization $\omega$ is recovered from
\beq
y_{N,M}(t_k,\omega)
=\!\!\sum_{i=1}^{\dim(\cH^M(k))}\!\!a_{k,i}e_{k,i}(t_k,\omega)\,,
\q
Y_{N,M-1}(t_k,\omega)
=\!\!\sum_{i=1}^{\dim(\cH^{M\!-\!1}(k))}\!\!b_{k,i}e_{k,i}(t_k,\omega)\,.
\eeq
Thus, the finite transposition method determines a common set of deterministic coefficients, while the pathwise values are obtained simply by evaluating the prescribed basis functions along the realization under consideration.
\er


\subsection{Convergence of the finite transposition method}\label{sec-ftm-convergence}

We now establish the convergence of the finite transposition method under Assumption {\rm\bf(B)}, without imposing the additional regularity conditions required for the convergence-rate estimate in Theorem \ref{ftm-rate}.

\bt{ftm-convergence}
Suppose that Assumption {\rm\bf(B)} holds and that $y_T\in L^2_{\cF_T}(\Omega;\dbR^n)$. Let $\big(y(\cd),Y(\cd)\big)$ be the adapted solution to BSDE \eqref{bsde0}.
For each $M\in\dbN$ and $N\geq N_*$, let $\big(y_{N,M}(\cd),Y_{N,M-1}(\cd)\big)$ be its finite transposition solution with terminal approximation $y_T^\pi=\G_M^N y_T$.
Then
\bel{w506e1}
\lim_{R\to\infty}\sup_{N,M\geq R}\Psi_{N,M}=0\,,
\ee
where
\beq
\Psi_{N,M}
= \sup_{0\leq t\leq T}\me\big[\big|y_{N,M}(t)-y\big(\nu(t)\big)\big|^2\big]
+\me\int_0^T\big|Y_{N,M-1}(t)-Y(t)\big|^2\rd t\,.
\eeq
Moreover,
\bel{w506e1a}
\lim_{M\to\infty}\limsup_{N\to\infty}\Psi_{N,M}
=\lim_{N\to\infty}\limsup_{M\to\infty}\Psi_{N,M}
=0\,.
\ee
\et

\begin{proof}
Fix $\e>0$. The proof is based on a standard density and stability argument. We first approximate the terminal condition and the generator by smoother data which are independent of the numerical parameters $N$ and $M$. Suppose that $\big(y^\e(\cd),Y^\e(\cd)\big)$ and $\big(y_{N,M}^\e(\cd),Y_{N,M-1}^\e(\cd)\big)$ are the corresponding adapted solution and finite transposition solution, respectively.
We then estimate the three differences $\big(y(\cd),Y(\cd)\big)-\big(y^\e(\cd),Y^\e(\cd)\big)$, $\big(y^\e(\cd),Y^\e(\cd)\big)-\big(y_{N,M}^\e(\cd),Y_{N,M-1}^\e(\cd)\big)$, and $\big(y_{N,M}^\e(\cd),Y_{N,M-1}^\e(\cd)\big)-\big(y_{N,M}(\cd),Y_{N,M-1}(\cd)\big)$ in Steps 1--3, respectively, and combine these estimates in Step 4.

\ss
\no
{\bf Step 1.} We first construct smooth approximations of the terminal condition $y_T$ and the generator $f(\cd,\cd,\cd)$.
Recall that the Wiener chaos expansion
$y_T=\me[y_T]+\sum_{m=1}^{\infty}J_m(h_m)$,
where $h_m =\big\{h_m^{\ell_1,\cds,\ell_m}\big\}_{1\leq\ell_1,\cds,\ell_m\leq d}\in L^2_{S,d}(T^m;\dbR^n)$ and
$\sum_{m=1}^{\infty}\frac{1}{m!} |h_m|_{L^2_{S,d}(T^m;\dbR^n)}^2
=\me\big[\big|y_T-\me[y_T]\big|^2\big]
<\infty$.
Hence, there exists $K_\e\in\dbN$ such that
\beq
\sum_{m=K_\e+1}^{\infty}
\frac{1}{m!}|h_m|_{L^2_{S,d}(T^m;\dbR^n)}^2
\leq \frac{\e}{2}\,.
\eeq
For each $1\leq m\leq K_\e$, by the density of smooth symmetric kernels in $L^2_{S,d}(T^m;\dbR^n)$, we can choose
$h_m^\e =\big\{h_m^{\e,\ell_1,\cds,\ell_m}\big\}_{1\leq\ell_1,\cds,\ell_m\leq d}\in L^2_{S,d}(T^m;\dbR^n)$ whose components are smooth on $[0,T]^m$, such that
\beq
\sum_{m=1}^{K_\e}
\frac{1}{m!}
|h_m-h_m^\e|_{L^2_{S,d}(T^m;\dbR^n)}^2
\leq \frac{\e}{2}\,.
\eeq
Define
\bel{w520e11}
y_T^\e
= \me[y_T]+\sum_{m=1}^{K_\e}J_m(h_m^\e)\,.
\ee
By the orthogonality of different Wiener chaoses and the It\^o isometry,
\bel{w520e12}
\me\big[\big|y_T-y_T^\e\big|^2\big]
\leq \e\,.
\ee

Since $y_T^\e$ belongs to a finite sum of Wiener chaoses with smooth deterministic kernels, it follows that $y_T^\e\in\bigcap_{r\geq1}\dbD^{r,q}(\dbR^n)$ for any $q\geq 2$.
Moreover, for every $1\leq k\leq K_\e$ and $\ell_1,\cds,\ell_k\in\{1,\cds,d\}$,
\beq
\me\big[D_{s_1,\cds,s_k}^{\ell_1,\cds,\ell_k}y_T^\e\big]
=h_k^{\e,\ell_1,\cds,\ell_k}(s_1,\cds,s_k)\,,
\eeq
whereas the left-hand side vanishes for $k>K_\e$.
Since the functions $h_k^\e$ are smooth on the compact sets $[0,T]^k$, all the regularity conditions on the terminal value in Assumption {\rm\bf(H)$_m$} hold for every $m\in\dbN$, with constants depending on $\e$ but independent of $m$.

We next regularize the generator with respect to its second and third variables. Let $\rho\in C_c^\infty(\dbR^n\times \dbR^{n\times d})$ be a nonnegative mollifier satisfying $\int_{\dbR^n\times \dbR^{n\times d}}\rho(z,Z)\rd z\rd Z$ $=1$.
Set $\rho_\e(z,Z)=\frac{1}{\e^{n(d+1)}} \rho\big(\frac z\e,\frac Z\e\big)$, and define
\beq
f^\e(t,y,Y) = \int_{\dbR^{n}\times\dbR^{n\times d}} f(t,y-z,Y-Z)\rho_\e(z,Z)\rd z\rd Z\,.
\eeq
By Assumption {\rm\bf(B)}, $f^\e(\cd,\cd,\cd)$ satisfies Assumption {\rm\bf(B)} with a Lipschitz constant independent of $\e$, and
\bel{w506e3}
\sup_{\substack{
t\in[0,T],\,y\in\dbR^n\\
Y\in\dbR^{n\times d}
}} |f(t,y,Y)-f^\e(t,y,Y)|
\leq \cC\e\,.
\ee
Furthermore, $f^\e(\cd,\cd,\cd)$ remains $\frac12$-H\"older continuous with respect to $t$ and is smooth with respect to $(y,Y)$. Its first-order derivatives with respect to $(y,Y)$ are bounded uniformly in $\e$, while, for each fixed $\e>0$, its second-order derivatives are also bounded.
Consequently, for each fixed $\e>0$, the data $\big(y_T^\e,f^\e(\cd,\cd,\cd)\big)$ satisfy Assumption {\rm\bf(H)$_m$} for every $m\in\dbN$, with a common regularity constant which may depend on $\e$ but is independent of $m$.

Let $\big(y^\e(\cd),Y^\e(\cd)\big)$ be the adapted solution to the BSDE with terminal condition $y_T^\e$ and generator $f^\e$. The standard stability estimate for BSDEs, together with \eqref{w520e12} and \eqref{w506e3}, yields
\bel{w506e4}
\bal
&\me\big[\sup_{0\leq t\leq T}|y(t)-y^\e(t)|^2\big]
+\me\int_0^T|Y(t)-Y^\e(t)|^2\rd t\\
&\leq \cC\Big\{\me\big[\big|y_T-y_T^\e\big|^2\big]
+\me\int_0^T\big|f\big(t,y(t),Y(t)\big) -f^\e\big(t,y(t),Y(t)\big)\big|^2\rd t
\Big\}\\
&\leq \cC(\e+\e^2)\,.
\eal
\ee
Here and below, $\cC$ is independent of $\e$.

\ss
\no
{\bf Step 2. }
For each $M\in\dbN$ and $N\geq N_*$, let $\big(y_{N,M}^\e(\cd),Y_{N,M-1}^\e(\cd)\big)$ be the finite transposition solution associated with the regularized data $\big(y_T^\e,f^\e(\cd,\cd,\cd)\big)$ and the terminal approximation $y_T^{\pi,\e}=\G_M^N y_T^\e$. For fixed $\e>0$, the data $\big(y_T^\e,f^\e(\cd,\cd,\cd)\big)$ satisfy Assumption {\rm\bf(H)$_M$} for every $M\in\dbN$. Therefore, Theorem \ref{ftm-rate} gives
\bel{w520e10}
\bal
\sup_{0\leq t\leq T} \me\big[ \big|y_{N,M}^\e(t)\!-\!y^\e\big(\nu(t)\big)\big|^2\big]
\!+\!\me\!\int_0^T\! \big|Y_{N,M-1}^\e(t)\!-\!Y^\e(t)\big|^2\rd t
\leq \cC_\e \big(\frac1N\!+\!\frac1{M+1}\big)\,,
\eal
\ee
where $\cC_\e$ may depend on $\e$, but is independent of $N$ and $M$.

\ss
\no
{\bf Step 3.}
We next establish a stability estimate between the finite transposition solutions corresponding to the original and regularized data.
For simplicity, write
\beq
&\delta y_k =y_{N,M}(t_k)-y_{N,M}^\e(t_k)\,,
\qquad
\delta Y_k =Y_{N,M-1}(t_k)-Y_{N,M-1}^\e(t_k)\,,\\
\mbox{and}\q
&\delta f_k=\G_M f\big(t_k,y_{N,M}(t_k),Y_{N,M-1}(t_k)\big)-\G_M f^\e\big(t_k,y_{N,M}^\e(t_k),Y_{N,M-1}^\e(t_k)\big)\,.
\eeq
At the terminal node, we set
$\delta y_N=\G_M^Ny_T-\G_M^Ny_T^\e$.

The representation \eqref{w23e10a}, applied successively at $t_k$ and $t_{k+1}$, gives the local relations
\beq
&\delta y_k=\me\big(\delta y_{k+1}\,\big|\,\cF_{t_k}\big)-\tau\,\delta f_k\,,\\
\mbox{and}\q
&\delta Y_k=\frac1\tau\me\big(\delta y_{k+1}(\D_{k+1}W)^\top\,\big|\,\cF_{t_k}\big)\,.
\eeq
Since the orthogonal projection $\G_M$ is a contraction in $L^2_{\cF_T}(\Omega;\dbR^n)$, Assumption {\rm\bf(B)} and \eqref{w506e3} imply
\bel{w506e5}
\me\big[|\delta f_k|^2\big]
\leq
\cC\big(\me\big[|\delta y_k|^2\big]
+\me\big[|\delta Y_k|^2\big]
+\e^2\big)\,.
\ee

We now derive a one-step stability estimate.
$\me(\delta y_{k+1}\,|\,\cF_{t_k})$ and $\delta Y_k\D_{k+1}W$ are orthogonal in $L^2_{\cF_T}(\Omega;\dbR^n)$, and the latter is the orthogonal projection of $\delta y_{k+1}-\me(\delta y_{k+1}\,|\,\cF_{t_k})$ onto the space generated by $\D_{k+1}W$ with $\cF_{t_k}$-measurable coefficients. Hence
\beq
\me\big[|\delta y_{k+1}|^2\big]
\geq \me\big[ \big|\me(\delta y_{k+1}\,|\,\cF_{t_k})\big|^2\big]
+\tau\me\big[|\delta Y_k|^2\big]\,.
\eeq
Using $\me(\delta y_{k+1}\,|\,\cF_{t_k})=\delta y_k+\tau\delta f_k$, Young's inequality, and \eqref{w506e5}, we obtain, for sufficiently small $\tau$,
\bel{w506e6}
\me\big[|\delta y_k|^2\big]
+\frac{\tau}{2}\me\big[|\delta Y_k|^2\big]
\leq (1+\cC\tau) \me\big[|\delta y_{k+1}|^2\big]+\cC\tau\e^2\,.
\ee
%
Applying the discrete Gronwall inequality to \eqref{w506e6}, and then summing the resulting estimates over $k=0,1,\cds,N-1$, we arrive at
\begin{align*}
\sup_{0\leq t\leq T} \me\big[|y_{N,M}(t)-y_{N,M}^\e(t)|^2\big]
+\me\int_0^T|Y_{N,M-1}(t)-Y_{N,M-1}^\e(t)|^2\rd t
\leq \cC\big(\me\big[|\G_M^N(y_T-y_T^\e)|^2\big]+\e^2\big)\,.
\end{align*}
Since $\G_M^N$ is an $L^2$-orthogonal projection, it follows from \eqref{w520e12} that
\bel{w506e8}
\bal
\sup_{0\leq t\leq T}\me\big[|y_{N,M}(t)-y_{N,M}^\e(t)|^2\big]
+\me\int_0^T|Y_{N,M-1}(t)-Y_{N,M-1}^\e(t)|^2\rd t
\leq \cC(\e+\e^2)\,.
\eal
\ee

\ss
\no
{\bf Step 4.} We finally combine the preceding estimates. By the triangle inequality, \eqref{w506e4}, \eqref{w520e10}, and \eqref{w506e8}, we obtain
\bel{w506e9}
0\leq \Psi_{N,M}
\leq \cC(\e+\e^2)
+\cC_\e\big( \frac1N+\frac1{M+1}\big)\,.
\ee
Here $\cC$ is independent of $\e,N,M$, whereas $\cC_\e$ is
independent of $N$ and $M$.

For any fixed $\e>0$, \eqref{w506e9} yields
\beq
0\leq \limsup_{R\to\infty}
\sup_{N,M\geq R}\Psi_{N,M}
\leq
\cC(\e+\e^2)\,.
\eeq
Since $\e>0$ is arbitrary, letting $\e\to0$ gives
\beq
\lim_{R\to\infty}
\sup_{N,M\geq R}\Psi_{N,M}=0\,.
\eeq

The same argument yields \eqref{w506e1a}.
This completes the proof.
\end{proof}


\subsection{Picard iteration}\label{sec-picard}

The representation \eqref{w23e10a}, or \eqref{w23e6}--\eqref{w23e9}, shows that the finite transposition method defined by the variational identity \eqref{w9e3} is implicit.
Except in some special cases, such as linear BSDEs or BSDEs whose generators are independent of $y$, it is generally impractical to compute the finite transposition solution directly from \eqref{w9e3}.
We therefore introduce the following Picard iteration.

\ss\no
{\it For any $N,M\in \dbN$, set $\big(y^{0}_{N,M}(\cd),\,Y^0_{N,M-1}(\cd)\big)=(0,0)$.
For $p\geq1$, let $\big(y^p_{N,M}(\cd),\,Y^p_{N,M-1}(\cd)\big)\in \dbH_{N,M}^n\times \dbH_{N,M-1}^{n\times d}$ be uniquely determined by the following variational identity
\begin{align}\label{w602e3}
\me\big[\langle z(T),y_T^\pi\rangle\big]
=\me\!\int_{0}^T \! \Big[ \big\langle z(t), f\big(\nu(t),y^{p-1}_{N,M}(t),Y^{p}_{N,M-1}(t)\big)\big\rangle
\!+\! \big\langle u(t),y^{p}_{N,M}(t)\big\rangle
\!+\!\big\langle v(t), Y^{p}_{N,M-1}(t)\big\rangle \Big]\! \rd t\,,
\end{align}
where $u(\cd)\in \dbH_{N,M}^n,\, v(\cd)\in \dbH_{N,M-1}^{n\times d}$, $z(\cd)$ solves \eqref{w9e4}, and $y^\pi_T\in[\cH^M(N)]^n$.
Whenever the terminal value is involved, we use the convention $y^p_{N,M}(T)=y_T^\pi$.
}

\ss

Write
\bel{w602e4}
\bal
y^p_{N,M}(\cd)=\sum_{k=0}^{N-1}\sum_{i=1}^{\dim(\cH^M(k))}a^p_{k,i}e_{k,i}(\cd)\,,\q
Y^p_{N,M-1}(\cd)=\sum_{k=0}^{N-1}\sum_{i=1}^{\dim(\cH^{M-1}(k))}b^p_{k,i}e_{k,i}(\cd)\,,\q p\in \dbN\,.
\eal
\ee
Proceeding as in Section \ref{ftm1}, we obtain the following representation of the Picard iterates.

\bt{w602t1}
Suppose that Assumptions {\rm\bf (P)} and {\rm\bf (B)} hold. Then, for any $p,N,M\in\dbN$, the Picard iterate $\big(y^p_{N,M}(\cd),Y^p_{N,M-1}(\cd)\big)\in \dbH_{N,M}^{n}\times\dbH_{N,M-1}^{n\times d}$ obtained by \eqref{w602e3} admits the following coefficient representation:
\bel{w602e5}
\lt\{
\bal
&a^p_{k,i}=\me\big[\sqrt{\tau } h_{k,i}\, y^\pi_T\big]
   -\me\int_{t_k}^T \sqrt{\tau } h_{k,i}\, f\big(\nu(t),y^{p-1}_{N,M}(t),Y^{p}_{N,M-1}(t)\big) \rd t\,,\\
&b^p_{k,i}=\me\Big[\frac{h_{k,i}}{\sqrt\tau}
\Big(y_T^\pi-\int_{t_{k+1}}^T
 f\big(\nu(t),y^{p-1}_{N,M}(t),Y^p_{N,M-1}(t)\big)\rd t\Big)
(\D_{k+1}W)^\top\Big]\,.
\eal
\rt.
\ee
Equivalently,
{\small
\bel{w602e6}
\lt\{
\bal
&y^p_{N,M}(\cd)=\sum_{k=0}^{N-1}\chi_{[t_k,t_{k+1})}(\cd)\me\Big(y^\pi_T
       -\int_{t_k}^T  \G_M f\big(\nu(t),y^{p-1}_{N,M}(t),Y^{p}_{N,M-1}(t)\big) \rd t \,\Big|\,\cF_{t_k}\Big)\,,\\
&Y^p_{N,M-1}(\cd)=\sum_{k=0}^{N-1}\chi_{[t_k,t_{k+1})}(\cd)\me\Big( \frac{1}{\tau }\Big[ y^\pi_T
      - \int_{t_{k+1}}^T \G_M f\big(\nu(t),y^{p-1}_{N,M}(t),Y^{p}_{N,M-1}(t)\big)\rd t\Big]\\
 &\qq\qq\qq\times (\D_{k+1}W)^\top \,\Big|\,\cF_{t_k}\Big)\,.
\eal
\rt.
\ee
}
\et

By Remark \ref{w21r1} {\rm (2)}, the projected BSDE \eqref{w19a4} admits a unique adapted solution.
Moreover, the representations \eqref{w23e10a} and \eqref{w602e6} have the same conditional-expectation form as the implicit Euler approximation of BSDE \eqref{w19a4}, whereas $\big(y^p_{N,M}(\cd),\,Y^p_{N,M-1}(\cd)\big)$ is obtained by linearizing the dependence of the generator on $y$ through Picard iteration. Hence, under Assumption {\rm\bf (B)}, adapting the argument of \cite[Theorem 5]{Bender-Denk07} and using the $L^2$-contraction property of $\G_M$, we obtain
 \bel{w602e7}
 \bal
 \sup_{0\leq t\leq T} \me\big[|y_{N,M}(t)-y^p_{N,M}(t)|^2\big]+\me\int_0^T\!|Y_{N,M-1}(t)-Y^{p}_{N,M-1}(t)|^2\rd t
\leq \cC_1\big(1+\me\big[|y_T^\pi|^2\big]\big)\big(\frac 1 2+\frac{\cC_2}{N} \big)^p\,,
 \eal
 \ee
where $\cC_1,\,\cC_2$ are constants independent of $p,N,M$. Combining this estimate with Theorem \ref{ftm-rate}, we obtain the convergence rate of the Picard-iteration-based finite transposition method.

\bt{w602t2}
Let $p,M\in\dbN$ and $N\geq N_0$, where $N_0$ is a sufficiently large integer independent of $p$ and $M$. Suppose that Assumptions {\rm\bf(P)}, {\rm\bf(B)}, and {\rm\bf(H)$_M$} hold. Let $\big(y(\cd),Y(\cd)\big)$ be the adapted solution to BSDE \eqref{bsde0}, and let $\big(y^p_{N,M}(\cd),Y^p_{N,M-1}(\cd)\big)\in \dbH_{N,M}^{n}\times\dbH_{N,M-1}^{n\times d}$ be the Picard iterate defined by \eqref{w602e3} with $y^\pi_T=\G_M^Ny_T$.
Then there exist constants $\cC_1,\cC_2>0$, independent of $p,N,M$, such that
\beq
\bal
\sup_{0\leq t\leq T} \me\big[\big|y\big(\nu(t))-y^p_{N,M}(t)\big|^2\big]
   +\me\int_0^T\big|Y(t)-Y^{p}_{N,M-1}(t)\big|^2\rd t 
\leq \cC_1\Big[\frac 1 N+\frac{1}{M+1}
+ \big( \frac 1 2+\frac{\cC_2}{N} \big)^p\Big]\,.
\eal
\eeq
%
\et


\subsection{Finite transposition method for Markovian BSDEs}\label{Markov-BSDE}

In Theorem \ref{ftm-rate}, Assumption {\rm\bf(H)$_M$} is imposed to derive the convergence rate of the finite transposition method. We now replace this assumption with conditions tailored to Markovian BSDEs.
Moreover, instead of taking an abstract approximation $y_T^\pi\in[\cH^M(N)]^n$, we construct an explicit approximation of the terminal value $y(T)$.
%
Here, the terminal condition of BSDE
\eqref{bsde0} satisfies
\beq
y(T)=\varphi\big(x(T)\big)\,,
\eeq
where $x(\cd)$ solves the following SDE
\bel{fsde1}
\left\{
\bal
&\mathrm d x(t)=b\big(t,x(t)\big)\rd t+\si\big(t,x(t)\big)\rd W(t)\,,\q t\in [0,T]\,,\\
&x(0)=x_0\in\dbR^m\,.
\eal
\right.
\ee
To derive a convergence rate in the Markovian setting, we impose the following assumption:
\begin{enumerate}
\item [{\bf (F)}]
$ b: [0,T]\times\dbR^m\longrightarrow\dbR^m$, $\si:[0,T]\times\dbR^m\longrightarrow\dbR^{m\times d}$, and $\varphi: \dbR^m\longrightarrow\dbR^n$ are deterministic functions. Moreover, there exists a constant $L>0$ such that, for any $t_1,t_2,t\in[0,T]$ and $x_1,x_2\in\dbR^m$,
\begin{equation*}
\begin{array}{c}
\ds |b(t_1,x_1)-b(t_2, x_2)|+|\si(t_1,x_1)-\si(t_2, x_2)|\leq
     L(|t_1-t_2|^{1/2}+|x_1-x_2|)\,,\\
\ns\ds  |\varphi (x_1)-\varphi (x_2)|\leq
     L|x_1-x_2|\,,\\
\ns\ds |b(t, 0)|+|\si(t,0)|+|\varphi (0)|\leq L\,.
\end{array}
\end{equation*}

\end{enumerate}

We first recall two auxiliary results.

\begin{lemma}[{\cite[Proposition 1.2.4]{Nualart06}}]\label{w0312l2}
Let $\phi:\,\dbR^m\longrightarrow \dbR^r$ be a Lipschitz continuous function with $m,r\in\dbN$, and suppose that $\xi=(\xi_1,\cds,\xi_m)^\top\in \dbD^{1,2}(\dbR^m)$. Then $\phi(\xi)\in\dbD^{1,2}(\dbR^r)$.
\end{lemma}

\begin{lemma}[{\cite[Theorem 2.2.1]{Nualart06}}]\label{w0312l3}
Suppose that $b(\cd,\cd)$ and $\sigma(\cd,\cd)$ satisfy the globally Lipschitz and linear growth conditions.
Then the solution $x(\cd)$ to SDE \eqref{fsde1} satisfies $x(t)\in \dbD^{1,\infty}(\dbR^m)$ for any $t\in [0,T]$, and
\beq
\sup_{\th\in[0,T]}\sum_{\ell=1}^d
\me\Big[\sup_{t\in[\th,T]} \big|D_\th^\ell x(t)\big|^2\Big]
<\infty\,.
\eeq
\end{lemma}

We next construct an explicit approximation by using the Euler--Maruyama approximation of SDE \eqref{fsde1}. For details concerning the Euler--Maruyama method, see, e.g., \cite[Chapter 10]{Kloeden-Platen92}. Set
\bel{w9em1}
\lt\{
\bal
&x^\pi(t_0)=x_0\,,\\
&x^\pi(t_{k+1})\!=\!x^\pi(t_k)\!+\!b\big(t_k,x^\pi(t_k)\big)\tau \!+\!\si\big(t_k,x^\pi(t_k)\big)\D_{k+1}W\,,\,\, k=0,\cds,N\!-\!1\,,
\eal
\rt.
\ee
and define the step process $x^\pi(\cd)$ from the discrete values $\{x^\pi(t_k)\}_{k=0}^N$ by
\bel{w0420e1}
x^\pi(\cd)\deq x^\pi\big(\nu(\cd)\big)\,.
\ee

The following standard estimate gives the strong convergence rate of the Euler--Maruyama approximation; see, e.g., \cite[Theorem 10.2.2]{Kloeden-Platen92}.

\bl{fsde-convergence}
Suppose that Assumptions {\rm\bf(P)} and {\rm\bf(F)} hold, $x(\cd)$ is the adapted solution to SDE \eqref{fsde1}, and $x^\pi(\cd)$ is defined by \eqref{w9em1}--\eqref{w0420e1}.
Then there exists a positive constant $\cC$, depending only on $T$ and $L$, such that
\beq
\sup_{t\in [0,T]}\me\big[|x(t)-x^\pi(t)|^2\big]\leq \cC\tau \,.
\eeq
\el

Obviously, the numerical solution $x^\pi(\cd)\in L^2_{\dbF}(0,T;\dbR^m)$. Moreover, the recursive relation \eqref{w9em1} implies the following measurability property.

\bl{w9l1}
Under Assumptions {\rm\bf(P)} and {\rm\bf(F)}, it holds that
\bel{w9e1}
x^\pi(t_k)\in L^2_{\mathcal{G}_{t_k}}(\Omega;\dbR^m)\,, \qquad k=0,1,\cds,N\,.
\ee
\el

\br{w531r2}
Motivated by Lemma \ref{w9l1}, we choose the following terminal approximation
\beq
y^\pi_T=\G_M\varphi \big(x^\pi(T)\big)\,.
\eeq
Since $\varphi \big(x^\pi(T)\big)\in L^2_{\cG_T}(\Omega;\dbR^n)$, Remark \ref{w531r1} implies that
$y^\pi_T\in [\cH^M(N)]^n$.
\er

The following lemma provides the Malliavin regularity required in the subsequent convergence analysis under the Lipschitz condition in Assumption {\rm\bf(B)}.

\bl{malliavin-lipschitz}
Suppose that Assumption {\rm\bf(B)} holds and $y_T\in\dbD^{1,2}(\dbR^n)$. Let $\big(y(\cd),Y(\cd)\big)$ be the adapted solution to BSDE \eqref{bsde0}. Then
\beq
y(\cd)\in L^2\big(0,T;\dbD^{1,2}(\dbR^n)\big)
\q\mbox{and}\q
Y(\cd)\in L^2\big(0,T;\dbD^{1,2}(\dbR^{n\times d})\big)\,.
\eeq
Moreover,
\bel{malliavin-lipschitz-est}
\bal
&\sum_{\ell=1}^d\int_0^T\int_0^t \me\Big[ |D_\th^\ell y(t)|^2+|D_\th^\ell Y(t)|^2 \Big]\rd\th\rd t
\leq \cC \sum_{\ell=1}^d\int_0^T \me\big[|D_\th^\ell y_T|^2\big]\rd\th\,.
\eal
\ee
Furthermore, $f\big(\cd,y(\cd),Y(\cd)\big) \in L^2\big(0,T;\dbD^{1,2}(\dbR^n)\big)$, and
\bel{malliavin-generator-est}
\sum_{\ell=1}^d\int_0^T\int_0^t \me\Big[ \big|D_\th^\ell f\big(t,y(t),Y(t)\big)\big|^2\Big]\rd\th\rd t
\leq \cC \sum_{\ell=1}^d\int_0^T \me\big[|D_\th^\ell y_T|^2\big]\rd\th\,.
\ee
\el

\begin{proof}
Let $f^\e(\cd,\cd,\cd)$ be the regularized generator introduced in Step 1 of the proof of Theorem \ref{ftm-convergence}.
Recall that $f^\e$ is deterministic and smooth with respect to its second and third variables, its first-order derivatives are bounded uniformly in $\e$, and \eqref{w506e3} holds.

For $\e>0$, let $\big(y^\e(\cd), Y^\e(\cd)\big)$ be the adapted solution to
\bel{malliavin-bsde-eps}
\lt\{
\bal
&\mathrm d y^\e(t)
=f^\e\big(t, y^\e(t), Y^\e(t)\big)\rd t + Y^\e(t)\rd W(t)\,,
\q t\in[0,T]\,,\\
& y^\e(T)=y_T\,.
\eal
\rt.
\ee
Since $y_T\in\dbD^{1,2}(\dbR^n)$ and $f^\e(\cd,\cd,\cd)$ is deterministic and continuously differentiable with respect to its second and third variables, with bounded and continuous first derivatives, the standard Malliavin differentiability result for BSDEs gives
$ y^\e(\cd)\in
L^2\big(0,T;\dbD^{1,2}(\dbR^n)\big)$
and
$ Y^\e(\cd)\in
L^2\big(0,T;\dbD^{1,2}(\dbR^{n\times d})\big)$;
see, e.g.,
\cite[Theorem 5.1]{Mastrolia-Possamai-Reveillac17}.

For each $\ell=1,\cds,d$, adaptedness implies that
\beq
D_\th^\ell y^\e(t)=0\,,
\q
D_\th^\ell Y^\e(t)=0\,,
\eeq
for $\rd\th\otimes\rd t\otimes\rd\dbP$-a.e.
$(\th,t,\omega)$ with $t<\th$. For a.e. $\th\in[0,T]$, the
Malliavin chain rule then gives
\bel{malliavin-linearized}
\lt\{
\bal
&\mathrm d D_\th^\ell y^\e(t)
=\!\Big[{\mathcal A}^\e(t)D_\th^\ell y^\e(t)
\!+\!{\mathcal B}^\e(t)\big[D_\th^\ell Y^\e(t)\big]\Big]\rd t
\!+\!D_\th^\ell Y^\e(t)\rd W(t)\,, \,\,\, t\in[\th,T]\,,\\
&D_\th^\ell y^\e(T)
=D_\th^\ell y_T\,,
\eal
\rt.
\ee
where ${\mathcal A}^\e(t) \!=\!\pr_y\!f^\e\!\big(t,y^\e(t),Y^\e(t)\big)
\!\in\dbR^{n\times n}$,
while ${\mathcal B}^\e(t) \!=\!\pr_Y\!f^\e\!\big(t,y^\e(t),Y^\e(t)\big):
\dbR^{n\times d}\longrightarrow\dbR^n$.
By the uniform boundedness of the first-order derivatives of $f^\e(\cd,\cd,\cd)$,
\beq
\big|{\mathcal A}^\e(t)u\big|
+\big|{\mathcal B}^\e(t)[V]\big|
\leq\cC\big(|u|+|V|\big)
\eeq
for all $u\in\dbR^n$ and $V\in\dbR^{n\times d}$, for
$\rd t\otimes\rd\dbP$-a.e. $(t,\omega)$, where $\cC$ is independent
of $\e$.

Fix $\ell\in\{1,\cds,d\}$ and a.e. $\th\in[0,T]$. The standard $L^2$-estimate for linear BSDEs, together with the uniform boundedness of the first-order derivatives of $f^\e(\cd,\cd,\cd)$, yields
\beq
\me\big[|D_\th^\ell y^\e(t)|^2\big]
+\me\int_t^T|D_\th^\ell Y^\e(s)|^2\rd s
\leq \cC\me\big[|D_\th^\ell y_T|^2\big]\,,
\q
\th\leq t\leq T\,,
\eeq
where $\cC$ is independent of $\e,\th$, and $\ell$.
Subsequently,
\bel{malliavin-uniform-yY}
\bal
\sup_{\e>0}\bigg\{\sup_{0\leq t\leq T}\sum_{\ell=1}^d\int_0^t\me\big[|D_\th^\ell y^\e(t)|^2\big]\rd\th
+\sum_{\ell=1}^d\int_0^T\int_0^t
\me\big[|D_\th^\ell Y^\e(t)|^2\big]
\rd\th\rd t\bigg\}
\leq \cC\sum_{\ell=1}^d\int_0^T\me\big[|D_\th^\ell y_T|^2\big]\rd\th\,.
\eal
\ee


We next pass to the limit as $\e\to0$.
Since \eqref{bsde0} and \eqref{malliavin-bsde-eps} have the same terminal condition, the standard stability estimate for BSDEs and \eqref{w506e3} yield
\bel{malliavin-solution-conv}
\me\Big[ \sup_{0\leq t\leq T}| y^\e(t)-y(t)|^2\Big]
+\me\int_0^T | Y^\e(t)-Y(t)|^2\rd t
\leq \cC\e^2\,.
\ee

We first consider $ y^\e(\cd)$. Set $
{\mathcal K}\deq L^2(0,T;\dbR^n)$
and let ${\mathcal K}^d$ be the $d$-fold product of ${\mathcal K}$, endowed with its natural product norm.
By the closedness of the Malliavin derivative defined in Section \ref{fts} and Fubini's theorem, the operator $
D^{(1)}: L^2\big(0,T;\dbD^{1,2}(\dbR^n)\big)
\subset L^2(\Omega;{\mathcal K})
\longrightarrow L^2\big(\Omega\times(0,T);{\mathcal K}^d\big)$
is a closed linear operator. Its graph is therefore a closed linear subspace of
$ L^2(\Omega;{\mathcal K}) \times L^2\big(\Omega\times(0,T);{\mathcal K}^d\big)$ and hence is weakly closed.

By \eqref{malliavin-uniform-yY}, $\{D^{(1)}y^\e\}_{\e>0}$ is bounded in $L^2\big(\Omega\times(0,T);{\mathcal K}^d\big)$.
Therefore, there exist a sequence $\e_j\downarrow0$ and $G\in L^2\big(\Omega\times(0,T);{\mathcal K}^d\big)$ such that
\beq
D^{(1)} y^{\e_j}
\rightharpoonup G
\quad\mbox{weakly in }
L^2\big(\Omega\times(0,T);{\mathcal K}^d\big)\,.
\eeq
On the other hand, \eqref{malliavin-solution-conv} gives
\beq
 y^{\e_j}(\cd)\longrightarrow y(\cd)
\quad\mbox{strongly in }
L^2(\Omega;{\mathcal K})\,.
\eeq
Consequently,
\beq
\big(y^{\e_j},D^{(1)}y^{\e_j}\big)
\rightharpoonup (y,G)
\quad\mbox{weakly in } L^2(\Omega;{\mathcal K})
\times
L^2\big(\Omega\times(0,T);{\mathcal K}^d\big)\,.
\eeq
Since the graph of $D^{(1)}$ is weakly closed, we conclude that $y(\cd)\in L^2\big(0,T;\dbD^{1,2}(\dbR^n)\big)$ and $D^{(1)}y=G$.
Moreover, by the weak lower semicontinuity
of the norm and \eqref{malliavin-uniform-yY},
\bel{malliavin-y-est}
\sum_{\ell=1}^d\int_0^T\int_0^t
\me\big[|D_\th^\ell y(t)|^2\big]\rd\th\rd t
\leq
\cC
\sum_{\ell=1}^d\int_0^T
\me\big[|D_\th^\ell y_T|^2\big]\rd\th\,.
\ee

Applying the same argument to $Y^\e(\cd)$, with $L^2(0,T;\dbR^{n\times d})$ in place of ${\mathcal K}$, and using \eqref{malliavin-uniform-yY} and \eqref{malliavin-solution-conv}, we obtain
$Y(\cd)\in L^2\big(0,T;\dbD^{1,2}(\dbR^{n\times d})\big)$
and
\bel{malliavin-Y-est}
\sum_{\ell=1}^d\int_0^T\int_0^t
\me\big[|D_\th^\ell Y(t)|^2\big]\rd\th\rd t
\leq
\cC
\sum_{\ell=1}^d\int_0^T
\me\big[|D_\th^\ell y_T|^2\big]\rd\th\,.
\ee
Combining \eqref{malliavin-y-est} and
\eqref{malliavin-Y-est} proves
\eqref{malliavin-lipschitz-est}.

It remains to consider the generator. Under Assumption {\rm\bf(B)}, for any $t\in[0,T]$, the mapping $(y,Y)\longmapsto f(t,y,Y)$ from $\dbR^n\times\dbR^{n\times d}$ into $\dbR^n$ is Lipschitz continuous.
Therefore, by the Lipschitz chain rule for Malliavin derivatives (see, e.g., \cite[Proposition 1.2.4]{Nualart06}), $f\big(t,y(t),Y(t)\big)\in\dbD^{1,2}(\dbR^n)$ for a.e. $t\in[0,T]$, and
\beq
\big|D_\th^\ell f\big(t,y(t),Y(t)\big)\big|^2
\leq \cC\big[ |D_\th^\ell y(t)|^2
+|D_\th^\ell Y(t)|^2\big]
\eeq
for $\ell=1,\cds,d$ and $\rd\th\otimes\rd t\otimes\rd\dbP$-a.e. $(\th,t,\omega)$. Hence,
\begin{align*}
&\sum_{\ell=1}^d\int_0^T\int_0^t
\me\Big[
\big|D_\th^\ell
f\big(t,y(t),Y(t)\big)\big|^2
\Big]\rd\th\rd t\\
&\leq \cC \sum_{\ell=1}^d\int_0^T\int_0^t
\me\Big[ |D_\th^\ell y(t)|^2
+|D_\th^\ell Y(t)|^2 \Big]\rd\th\rd t\\
&\leq \cC \sum_{\ell=1}^d\int_0^T
\me\big[|D_\th^\ell y_T|^2\big]\rd\th\,,
\end{align*}
which proves \eqref{malliavin-generator-est}.
\end{proof}

The following theorem provides an error estimate for the terminal approximation and gives the corresponding convergence rate of the finite transposition method.

\bt{w531t2}
Let Assumptions {\rm\bf(P)}, {\rm\bf(B)} and {\rm\bf(F)} hold, and $y^\pi_T=\G_M\varphi \big(x^\pi(T)\big)$. Then
\bel{w531e5}
\me\big[|y(T)-y^\pi_T|^2\big]\leq \cC\big(\frac 1 N+\frac 1{M+1}\big)\,.
\ee
Furthermore,
the finite transposition solution $\big(y_{N,M}(\cd),Y_{N,M-1}(\cd)\big)$ satisfies
\bel{w104e1}
\bal
\sup_{0\leq t\leq T} \me\big[\big|y_{N,M}(t)\!-\!y\big(\nu(t)\big)\big|^2\big]\!+\!\me\!\int_0^T\!|Y_{N,M-1}(t)\!-\!Y(t)|^2\rd t
\leq \cC\big(\frac 1 N\!+\!\frac 1 {M+1}\big)\,,
\eal
\ee
where $\big(y(\cd),Y(\cd)\big)$ is the adapted solution to BSDE \eqref{bsde0}, and $\cC$ is independent of $N$ and $M$.

\et

\begin{proof}
Lemma \ref{w0312l3} and Assumption {\rm\bf (F)} imply that $x(T)\in \dbD^{1,2}(\dbR^m)$. Then Lemma \ref{w0312l2} yields $y(T)\in \dbD^{1,2}(\dbR^n)$. Therefore, utilizing Lemma \ref{fsde-convergence}, we have
\begin{align*}
\me\big[|y(T)\!-\!y^\pi_T|^2\big]
&\leq  2\Big(\me\big[\big|\varphi \big(x(T)\big)\!-\!\G_M\varphi \big(x(T)\big)\big|^2\big]\!+\!\me\big[\big|\G_M\varphi \big(x(T)\big)\!-\!\G_M\varphi \big(x^\pi(T)\big)\big|^2\big]\Big)\\
&\leq 2\frac 1{M+1} \sum_{\ell=1}^d \me\int_0^T \big|D_\th^\ell\varphi(x(T))\big|^2\rd\th
+2\me\big[\big|\varphi \big(x(T)\big)-\varphi \big(x^\pi(T)\big)\big|^2\big]\\
&\leq \cC\big(\frac 1 N +\frac 1{M+1}\big)\,,
\end{align*}
which proves \eqref{w531e5}.

Under Assumptions {\rm\bf(P)}, {\rm\bf(B)} and {\rm\bf(F)}, the
$L^2$-regularity estimate \eqref{w709e2} holds; see Remark
\ref{w1102r1}. Therefore, the argument in Step 2 of the proof of
Theorem \ref{ftm-rate} yields \eqref{w1102e2}.

By Lemmas \ref{w0312l2} and \ref{w0312l3}, we have $y(T)=\varphi\big(x(T)\big)\in\dbD^{1,2}(\dbR^n)$. Hence, Lemma \ref{malliavin-lipschitz} gives the solution
$\big(y(\cd),Y(\cd)\big) \in L^2\big(0,T;\dbD^{1,2}(\dbR^n)\big)
\times L^2\big(0,T;\dbD^{1,2}(\dbR^{n\times d})\big)$,
and $f\big(\cd,y(\cd),Y(\cd)\big)\in L^2\big(0,T;\dbD^{1,2}(\dbR^n)\big)$.
Since $f\big(t,y(t),Y(t)\big)$ is $\cF_t$-measurable, $D_\th^\ell f\big(t,y(t),Y(t)\big)=0$ for $\th>t$ and $\ell=1,\cds,d$.
Therefore, combining \eqref{w1227e6a} with \eqref{malliavin-generator-est}, we obtain
\beq
\me\int_0^T
\big|f\big(t,y(t),Y(t)\big)
-\G_M f\big(t,y(t),Y(t)\big)\big|^2\rd t
\leq \frac{\cC}{M+1}\,.
\eeq
Finally, \eqref{w104e1} follows by combining the above estimate, \eqref{w1102e2}, and \eqref{w531e5}.
This completes the proof.
\end{proof}


\subsection{Implementation by the Monte Carlo method}\label{MC}

The coefficient representation \eqref{w602e5} of the Picard iterate $\big(y^p_{N,M}(\cd),Y^p_{N,M-1}(\cd)\big)$ involves expectations. In this subsection, we approximate these expectations by Monte Carlo sampling.
By Theorem \ref{w602t1}, the pair of processes $\big(y^p_{N,M}(\cd),Y^p_{N,M-1}(\cd)\big)$ admits the following representation: for any $t\in [t_k, t_{k+1})$ with $k=0,1,\cds,N-1$,
\bel{w603e1}
\lt\{
\bal
&y^p_{N,M}(t)=\me\Big(y^{p}_{N,M}(t_{k+1})
       -\int_{t_k}^{t_{k+1}}  \G_M f\big(t_k,y^{p-1}_{N,M}(t_k),Y^{p}_{N,M-1}(t_k)\big) \rd t \,\Big|\,\cF_{t_k}\Big)\,,\\
&Y^p_{N,M-1}(t)=\me\Big( \frac{1}{\tau } y^{p}_{N,M}(t_{k+1})(\D_{k+1}W)^\top\, \Big|\,\cF_{t_k}\Big)\,.
\eal
\rt.
\ee
Equivalently, the corresponding coefficients satisfy
\bel{w603e2}
\lt\{
\bal
&a^p_{k,i}=\me\big[ \sqrt{\tau } h_{k,i}\, y^{p}_{N,M}(t_{k+1})\big]
   \!-\!\me\!\int_{t_k}^{t_{k+1}}\! \sqrt{\tau } h_{k,i}\,  f\big(t_k,y^{p-1}_{N,M}(t_k),Y^{p}_{N,M-1}(t_k)\big) \rd t\,,\\
 & \qq\qq\qq\qq\qq\qq\qq\qq\qq\qq i=1,2,\cds, \dim\big(\cH^M(k)\big)\,, \\
&b^p_{k,i}=\me\Big[\frac{h_{k,i}}{\sqrt{\tau }} y^{p}_{N,M}(t_{k+1})(\D_{k+1}W)^\top\Big]\,,
\qq\qq i=1,2,\cds, \dim\big(\cH^{M-1}(k)\big)\,, \\
\eal
\rt.
\ee
where
\bel{w603e3}
\bal
y^p_{N,M}(t)=\sum_{i=1}^{\dim(\cH^M(k))}\frac{a^p_{k,i}}{\sqrt{\tau }} h_{k,i}\,, \q
Y^p_{N,M-1}(t)=\sum_{i=1}^{\dim(\cH^{M-1}(k))}\frac{b^p_{k,i}}{\sqrt{\tau }} h_{k,i}\,.
\eal
\ee
Fix a sample size $K\in\dbN$. We now construct coefficient-wise Monte Carlo approximations of the Picard iterates $y^p_{N,M}(\cd)$ and $Y^p_{N,M-1}(\cd)$. For each $k=0,1,\cds, N-1$, generate independent and identically distributed samples
\beq
&\big( h_{k,i}\, y^{p}_{N,M}(t_{k+1})\big)_j,\,
\big( h_{k,i}\,  f\big(t_k,y^{p-1}_{N,M}(t_k),Y^{p}_{N,M-1}(t_k)\big) \big)_j\,,\\
&\qq\qq\qq\qq\qq j=1,2,\cds,K,\, i=1,2,\cds,\dim\big(\cH^M(k)\big)\,,
\eeq
and
\beq
\big(  h_{k,i}\, y^{p}_{N,M}(t_{k+1})(\D_{k+1}W)^\top\big)_j,\,
\qq j=1,2,\cds,K,\, i=1,2,\cds,\dim\big(\cH^{M-1}(k)\big)\,.
\eeq
For each fixed $j$, all random variables carrying the subscript $j$ are evaluated along the same realization, while the resulting random vectors and matrices are independent and identically distributed with respect to $j$. The Monte Carlo samples are also independent of the realization used to evaluate the basis functions in \eqref{w603e4}. All random variables are defined on a suitably enlarged probability space, and $\me$ denotes expectation with respect to all sources of randomness.
The coefficient-wise Monte Carlo approximation of
$\big(y^p_{N,M}(t),Y^p_{N,M-1}(t)\big)$ is then defined by
\bel{w603e4}
\bal
y^p_{N,M,K}(t)=\sum_{i=1}^{\dim(\cH^M(k))}\frac{a^p_{k,i,K}}{\sqrt{\tau }} h_{k,i}\,, \q
Y^p_{N,M-1,K}(t)=\sum_{i=1}^{\dim(\cH^{M-1}(k))}\frac{b^p_{k,i,K}}{\sqrt{\tau }} h_{k,i}\,,
\eal
\ee
where $a^p_{k,i,K}$ and $b^p_{k,i,K}$ are the sample-average estimators of $a^p_{k,i}$ and $b^p_{k,i}$, respectively:
\bel{w603e5}
\lt\{
\bal
&a^p_{k,i,K}=\frac {\sqrt{\tau }} K\sum_{j=1}^K\Big[\big(  h_{k,i}\, y^{p}_{N,M}(t_{k+1})\big)_j
   -\t \big( h_{k,i}\, f\big(t_k,y^{p-1}_{N,M}(t_k),Y^{p}_{N,M-1}(t_k)\big) \big)_j \Big]\,,\\
&b^p_{k,i,K}=\frac{1}{K\sqrt{\tau }}\sum_{j=1}^K\big(  h_{k,i}\, y^{p}_{N,M}(t_{k+1})(\D_{k+1}W)^\top\big)_j
 \,.
\eal
\rt.
\ee

\br{w604r1}
For any fixed $p,N,M\in\dbN$, and for each $k=0,1,\cdots,N-1$ and $i=1,2,\cdots,\dim\big(\cH^M(k)\big)$, we have
\beq
 h_{k,i}\, y^{p}_{N,M}(t_{k+1})\,,\,\,
\int_{t_k}^{t_{k+1}} h_{k,i}\,  f\big(t_k, y^{p-1}_{N,M}(t_k),Y^{p}_{N,M-1}(t_k)\big) \rd t \in L^1_{\cF_T}(\Omega;\dbR^n)\,.
\eeq
Moreover, for $i=1,2,\cds, \dim\big(\cH^{M-1}(k)\big)$,
\beq
 h_{k,i}\, y^{p}_{N,M}(t_{k+1})( \D_{k+1}W)^\top \in L^1_{\cF_T}(\Omega;\dbR^{n\times d})\,.
\eeq
Therefore, the strong law of large numbers yields
\bel{w604e1}
a^p_{k,i,K}\rightarrow a^p_{k,i}\,,\q b^p_{k,i,K}\rightarrow b^p_{k,i} \q \as
\ee
as $K\to\infty$. Consequently, for any fixed $t\in[0,T)$,
\bel{w604e2}
y^p_{N,M,K}(t)\rightarrow y^p_{N,M}(t)\,,\q
Y^p_{N,M-1,K}(t)\rightarrow Y^p_{N,M-1}(t) \q \as
\ee
\er

We next estimate the error between the Picard iterate and its coefficient-wise Monte Carlo approximation.

\bl{wang524a1}
For any $x\in\dbR$ and $m,r\in\dbN_0$, it holds that
\bel{wang524a2}
m!r!H_m(x)H_r(x)=\sum_{k=0}^{m\wedge r} k! C_m^k C_r^k (m+r-2k)!H_{m+r-2k}(x)\,.
\ee
Moreover, if $X$ is a standard normal random variable, then
\beq
\bal
(r!)^2\me\big[H_r^4(X)\big]
=\sum_{k=0}^{r} C_{2k}^{k}\big(C_r^k\big)^2\,.
\eal
\eeq
\el

\begin{proof}
The assertion \eqref{wang524a2} follows from \cite[p.~195]{Erdelyi-Magnus-Oberhettinger-Tricomi53}. Setting $m=r$ in \eqref{wang524a2}, we obtain
\beq
H_r^2(x)=\sum_{k=0}^r \frac{k!}{(r!)^2} (C_r^k )^2 (2r-2k)! H_{2r-2k}(x)
=\sum_{k=0}^r \frac{(2r-2k)!}{k!\big[(r-k)!\big]^2} H_{2r-2k}(x)\,.
\eeq
Then, applying Lemma \ref{w512l1}, we conclude that
\beq
(r!)^2\me[ H_r^4(X)]
= \sum_{l=0}^r\frac{(2r-2l)! (r!)^2}{(l!)^2\big[(r-l)!\big]^4}
=\sum_{l=0}^{r} C_{2r-2l}^{r-l} (C_r^l )^2=
\sum_{k=0}^r
C_{2k}^{k}(C_r^{k})^2\,,
\eeq
where, in the last equality, we have used the change of index $k=r-l$ and the identity $C_r^{r-k}=C_r^k$.
This completes the proof.
\end{proof}

The following lemma gives an $L^4$ estimate for the projection onto the Wiener chaoses up to order $M$.

\bl{w512l3}
For any $M,r\in\dbN$ and $X\in L^2_{\cF_T}(\Omega;\dbR^r)$, it holds that
\bel{w512e0}
\me\big[|\G_M X|^4\big]
\leq \Big( \frac{3^{M+1}-1}{2} \Big)^2\big(\me\big[|X|^2\big]\big)^2\,.
\ee
\el

\begin{proof}
Let $X=\sum_{m=0}^{\infty}J_m(h_m)$ be the Wiener chaos expansion of $X$, where $h_m\in L^2_{S,d}(T^m;\dbR^r)$; see \eqref{ex-y}. Write
\beq
J_m(h_m)
=\big(J_m^1,\cds,J_m^r\big)^\top,
\eeq
where each $J_m^i$ belongs to the scalar Wiener chaos of order $m$.
By the hypercontractivity of Wiener chaoses,
$\big|J_m^i\big|_{L_{\cF_T}^4(\Omega;\dbR)}
\leq 3^{m/2}\big|J_m^i\big|_{L_{\cF_T}^2(\Omega;\dbR)}$ for $i=1,\cds,r$; see, e.g., \cite[Theorem 1.4.1]{Nualart06}.
Hence, by the Cauchy--Schwarz inequality,
\begin{align*}
\me\big[|J_m(h_m)|^4\big]
\leq\Big(\!\sum_{i=1}^r\!
\big(\me[|J_m^i|^4]\big)^{1/2}\Big)^2
\leq3^{2m}
\Big(\!\sum_{i=1}^r\!\me[|J_m^i|^2]\Big)^2\!
=3^{2m}\big(\me[|J_m(h_m)|^2]\big)^2\,.
\end{align*}
Therefore,
\beq
\bal
|\G_MX|_{L_{\cF_T}^4(\Omega;\dbR^r)}
&\leq
\sum_{m=0}^M\big|J_m\big(h_m(\cd)\big)\big|_{L_{\cF_T}^4(\Omega;\dbR^r)}
\leq \Big( \sum_{m=0}^M3^m \Big)^{1/2}
\Big( \sum_{m=0}^M \big|J_m\big(h_m(\cd)\big)\big|_{L_{\cF_T}^2(\Omega;\dbR^r)}^2\Big)^{1/2}\,.
\eal
\eeq
Applying the orthogonality of the Wiener chaoses and then raising both sides to the fourth power yields \eqref{w512e0}.
\end{proof}

We next derive a fourth-moment estimate for the Picard iterates.
\bl{w512l2}
Let $M,N\in\dbN$ satisfy $N\geq N_{+}\deq \max\big\{N_*,\,\lceil4\cC_2\rceil\big\}$, where $N_*$ is given in \eqref{w-Nstar}, and $\cC_2$ is the constant appearing in \eqref{w602e7}. Suppose that $y_T^\pi\in[\cH^M(N)]^n$ and that $f(\cd,\cd,\cd)$ satisfies Assumption {\rm\bf(B)}. Then, for any $p\in\dbN$,
\beq
\max_{0\leq k\leq N}\me\big[\big|y_{N,M}^p(t_k)\big|^4\big]
+\sum_{k=0}^{N-1}\t^2\me\big[\big|Y_{N,M-1}^p(t_k)\big|^4\big]\leq \cC 9^M\big(1+\me\big[|y_T^\pi|^2\big]\big)^2\,,
\eeq
where $\cC$ is independent of $p,N,M$.
\el

\begin{proof}
Let $\big(y_{N,M}(\cd),Y_{N,M-1}(\cd)\big)$ be the finite transposition solution corresponding to $y_T^\pi$.
By applying Steps 2(a) and 2(b) in the proof of Theorem \ref{ftm-rate}, and arguing as in \eqref{w0310e3}, we obtain
\begin{align*}
&\me\big[\big|y_{N,M}(t_k)\big|^2\big]+\t \me\big[\big|Y_{N,M-1}(t_k)\big|^2\big]\\
&\leq \me\big[\big|y_0(t_k)\big|^2\big]+\me\Big[\int_{t_k}^{t_{k+1}}\big|Y_0(t)\big|^2\rd t\Big]\\
&\leq (1+\e)\me\big[\big|y_{N,M}(t_{k+1})\big|^2\big]
+3L^2\big(1+\e^{-1}\big)\t^2\Big(1 +\me\big[\big|y_{N,M}(t_k)\big|^2\big]+\me\big[\big|Y_{N,M-1}(t_k)\big|^2\big]\Big)\,,
\end{align*}
where the pair $\big(y_0(\cd),Y_0(\cd)\big)$ solves \eqref{w0310e1}.
By taking
$\e=12L^2\t$, the term involving $\me\big[\big|Y_{N,M-1}(t_k)\big|^2\big]$ can be absorbed into the left-hand side. Consequently, the discrete Gronwall inequality yields
\bel{w728e1}
\max_{0\leq k\leq N}\me\big[\big|y_{N,M}(t_k)\big|^2\big]
+\sum_{k=0}^{N-1}\t\me\big[\big|Y_{N,M-1}(t_k)\big|^2\big]
\leq \cC\big(1+\me\big[|y_T^\pi|^2\big]\big)\,.
\ee
Combining \eqref{w728e1} with the Picard-iteration estimate \eqref{w602e7}, and noting that $\frac12+\frac{\cC_2}{N}\leq\frac34$ for $N\geq 4\cC_2$, we have
\bel{w718e3}
\max_{0\leq k\leq N}\me\big[\big|y_{N,M}^p(t_k)\big|^2\big]+\sum_{k=0}^{N-1}\t \me\big[\big|Y_{N,M-1}^p(t_k)\big|^2\big]
\leq \cC\big(1+\me\big[|y_T^\pi|^2\big]\big)\,.
\ee

For every $k=0,1,\cds,N$, we have $y_{N,M}^p(t_k)\in[\cH^M(k)]^n$.
Thus, Lemma \ref{w512l3} gives
\bel{w512e0i}
\me\big[\big|y_{N,M}^p(t_k)\big|^4\big]
\leq
\Big( \frac{3^{M+1}-1}{2} \Big)^2
\big(\me\big[\big|y_{N,M}^p(t_k)\big|^2\big]\big)^2\,.
\ee
Similarly, since $Y_{N,M-1}^p(t_k)\in[\cH^{M-1}(k)]^{n\times d}$, we have
\begin{align}
\sum_{k=0}^{N-1}\tau^2\me\big[|Y_{N,M-1}^p(t_k)|^4\big]
&\leq\Big( \frac{3^{M}-1}{2} \Big)^2 \sum_{k=0}^{N-1} \tau^2\big( \me\big[|Y_{N,M-1}^p(t_k)|^2 \big]\big)^2\notag \\
&\leq\Big( \frac{3^{M}-1}{2} \Big)^2 \Big(\sum_{k=0}^{N-1} \tau \me\big[|Y_{N,M-1}^p(t_k)|^2 \big]\Big)^2\,.\label{w718e2}
\end{align}
Combining \eqref{w718e3}, \eqref{w512e0i}, and \eqref{w718e2} proves the assertion.
\end{proof}

\bt{w604t1}
Let $p,M,N,K\in\dbN$ with $N\geq N_+$, where $N_+$ is the constant specified in Lemma \ref{w512l2}. Suppose that $y_T^\pi\in[\cH^M(N)]^n$ and that $f(\cd,\cd,\cd)$ satisfies Assumption {\rm\bf(B)}. Then the following estimate holds
\bel{w604e3}
\bal
&\sup_{t\in[0,T)}\me\big[\big| y^p_{N,M}(t)-y^p_{N,M,K}(t) \big |^2\big]
+\me\int_0^T\big|Y^p_{N,M-1}(t)-Y^p_{N,M-1,K}(t)  \big |^2\rd t \\
&\qq\qq\leq  \frac{\cC\big(1+\me\big[|y_T^\pi|^2\big]\big)}{K} C_{dN+M-1}^M 9^{M}\,,
\eal
\ee
where $\cC$ is independent of $p,N,M$, and $K$.
\et

\begin{proof}
We adopt the Wiener chaos expansion to estimate the error of the Monte Carlo approximation. Set $f_k= f\big(t_k,y^{p-1}_{N,M}(t_k),Y^{p}_{N,M-1}(t_k)\big)$ with $k=0,1,\cds,N-1$.
Recall from Remark \ref{w9r1} that $\{h_{k,i}\}_{i=1}^{\dim(\cH^M(k))}$ is an orthonormal basis of $\cH^M(k)$, where
$
h_{k,i}
=\sqrt{\a^{k,i}!}\prod_{j=1}^k\prod_{\ell=1}^d
H_{\a_j^{k,i,\ell}}\big(\dbW(g_{j,\ell})\big)$.
By \eqref{w603e2}--\eqref{w603e5} and Remark
\ref{w9r1}, for any $t\in [t_k,t_{k+1})$, we have
\begin{align}
&\me\big[\big| y^p_{N,M}(t)-y^p_{N,M,K}(t) \big |^2\big]
=\frac{1}{\t}\sum_{i=1}^{\dim(\cH^M(k))}
  \me\big[|a^p_{k,i}-a^p_{k,i,K}|^2\big]
\notag\\
&\qq\leq\frac{2}{K} \sum_{i=1}^{\dim(\cH^M(k))}
\Big\{\me\big[h_{k,i}^2 |y^p_{N,M}(t_{k+1})|^2\big]
+\t^2 \me\big[h_{k,i}^2|f_k|^2\big] \Big\}\,.
\label{wang524a7}
\end{align}
Indeed, we have used the fact that, for any square-integrable random vector $Z$ and its i.i.d.\ copies $Z_1,\cds,Z_K$,
\beq
\me\Big[\Big| \frac1K\sum_{j=1}^K Z_j-\me[Z]\Big|^2\Big]
=\frac1K \me\big[ |Z-\me[Z]|^2\big]
\leq \frac1K\me[|Z|^2]\,.
\eeq

For $k=0$, \eqref{wang524a7} and \eqref{w718e3} yield $\me\big[\big| y^p_{N,M}(t)-y^p_{N,M,K}(t) \big |^2\big]\leq \frac {\cC(1+\me[|y_T^\pi|^2]) } {K}$.
We next consider $k\geq 1$. By H\"older's inequality,
\begin{align*}
\sum_{i=1}^{\dim(\cH^M(k))}
\me\big[h_{k,i}^2|y^p_{N,M}(t_{k+1})|^2\big]
\leq \Big(\me\big[|y^p_{N,M}(t_{k+1})|^4\big]\Big)^{1/2}
\sum_{i=1}^{\dim(\cH^M(k))}\big(\me[h_{k,i}^4]\big)^{1/2}\,.
\end{align*}
For each $i$, Lemma \ref{wang524a1} gives
\begin{align}
&\me[h_{k,i}^4]
=\prod_{j=1}^k\prod_{\ell=1}^d
\big(\a_j^{k,i,\ell}!\big)^2\me\big[H_{\a_j^{k,i,\ell}}^4\big(\dbW(g_{j,\ell})\big)\big]
=\prod_{j=1}^k\prod_{\ell=1}^d
\sum_{r=0}^{\a_j^{k,i,\ell}}C_{2r}^{r}
\big(C_{\a_j^{k,i,\ell}}^{r}\big)^2
\notag\\
&\q\leq \prod_{j=1}^k\prod_{\ell=1}^d
\sum_{r=0}^{\a_j^{k,i,\ell}}
2^{2r}\big(C_{\a_j^{k,i,\ell}}^{r}\big)^2
\leq \prod_{j=1}^k\prod_{\ell=1}^d
\Big(\sum_{r=0}^{\a_j^{k,i,\ell}}2^r C_{\a_j^{k,i,\ell}}^{r}\Big)^2
=9^{|\a^{k,i}|}\,.
\label{wang524a9}
\end{align}
Moreover, $\#\big\{ \a\in\Lambda_d(k):|\a|=m\big\}=C_{kd+m-1}^m$; see Remark \ref{w9r1} (2).
Hence, by Lemma \ref{w512l2},
\begin{align}
\sum_{i=1}^{\dim(\cH^M(k))}
\me\big[ h_{k,i}^2|y^p_{N,M}(t_{k+1})|^2\big]
&\leq \Big(\me\big[|y^p_{N,M}(t_{k+1})|^4\big] \Big)^{1/2}
\sum_{m=0}^M C_{kd+m-1}^m3^m\notag\\
& \leq \cC3^M
\big(1+\me[|y_T^\pi|^2]\big)
\sum_{m=0}^M C_{kd+m-1}^m3^m.
\label{wang524a8}
\end{align}

Similarly, by Assumption {\rm\bf(B)}, H\"older's inequality,
\eqref{wang524a9}, and Lemma \ref{w512l2}, we obtain
\begin{align}
\t^2 \sum_{i=1}^{\dim(\cH^M(k))}
\me\big[h_{k,i}^2|f_k|^2\big]
&\leq \cC\t^2
\Big(1+\me\big[|y^{p-1}_{N,M}(t_k)|^4\big]
+\me\big[|Y^p_{N,M-1}(t_k)|^4\big]\Big)^{1/2}
\sum_{m=0}^M C_{kd+m-1}^m3^m
\notag\\
&\leq \cC3^M
\big(1+\me[|y_T^\pi|^2]\big)
\sum_{m=0}^M C_{kd+m-1}^m3^m.
\label{wang529a8}
\end{align}
Here we have used the estimate for the $Y$-component in
Lemma \ref{w512l2} and absorbed the fixed factor $\t\leq T$ into
$\cC$.

Combining \eqref{wang524a7}, \eqref{wang524a8}, and \eqref{wang529a8}, and observing that $C_{kd+m-1}^m\leq C_{kd+M-1}^M$ for $m=0,1,\cds,M$,
we obtain, for $t\in[t_k,t_{k+1})$,
\bel{wang524b1}
\bal
&\me\big[\big|y^p_{N,M}(t)-y^p_{N,M,K}(t)\big|^2\big]
\leq \frac{\cC3^M\big(1\!+\!\me[|y_T^\pi|^2]\big)}{K}C_{kd+M-1}^M
\sum_{m=0}^M3^m
\leq\frac{\cC \big(1\!+\!\me[|y_T^\pi|^2]\big)}{K}C_{dN+M-1}^M9^M\,.
\eal
\ee

We next estimate the $Y$-component. Since $b^p_{k,i},b^p_{k,i,K}\in\dbR^{n\times d}$, the same orthogonality
argument gives
\begin{align}
&\me\big[\big|Y^p_{N,M-1}(t)-Y^p_{N,M-1,K}(t)\big|^2\big]
=\frac{1}{\t}\sum_{i=1}^{\dim(\cH^{M-1}(k))}
\me\big[|b^p_{k,i}-b^p_{k,i,K}|^2\big]
\notag\\
&\qq\q\leq \frac{1}{K\t} \sum_{i=1}^{\dim(\cH^{M-1}(k))}
\me\Big[h_{k,i}^2 |y^p_{N,M}(t_{k+1})|^2
\Big| \frac{\D_{k+1}W}{\sqrt{\t}}\Big|^2 \Big]\,.
\label{w604e5}
\end{align}
Here and below, $|\cdot|$ denotes the Frobenius norm for
$\dbR^{n\times d}$-valued quantities; in particular,
\beq
\big| y^p_{N,M}(t_{k+1})(\D_{k+1}W)^\top\big|^2
=|y^p_{N,M}(t_{k+1})|^2|\D_{k+1}W|^2\,.
\eeq

Since $\frac{\D_{k+1}W}{\sqrt{\t}}$ is a standard Gaussian random vector in $\dbR^d$, it follows that
$
\me\big[\big|\frac{\D_{k+1}W}{\sqrt{\t}}\big|^4\big]
=d(d+2)$.
Moreover, $\frac{\D_{k+1}W}{\sqrt{\t}}$ is independent of $h_{k,i}$.
Therefore, by H\"older's inequality,
\begin{align*}
\me\Big[ h_{k,i}^2 |y^p_{N,M}(t_{k+1})|^2
\Big|\frac{\D_{k+1}W}{\sqrt{\t}}\Big|^2\Big]
&\leq \Big(\me\big[|y^p_{N,M}(t_{k+1})|^4\big]\Big)^{1/2}
\Big(\me\Big[h_{k,i}^4\Big|\frac{\D_{k+1}W}{\sqrt{\t}}\Big|^4\Big]
\Big)^{1/2}
\notag\\
&=\sqrt{d(d+2)}
\Big(\me\big[ |y^p_{N,M}(t_{k+1})|^4\big]\Big)^{1/2}
\big(\me[h_{k,i}^4]\big)^{1/2}\,.
\end{align*}

For $k\geq1$, summing over $i=1,\cds,\dim(\cH^{M-1}(k))$, using \eqref{wang524a9} and Lemma \ref{w512l2}, gives
\begin{align}
\sum_{i=1}^{\dim(\cH^{M-1}(k))}
\me\Big[h_{k,i}^2 |y^p_{N,M}(t_{k+1})|^2
\Big|\frac{\D_{k+1}W}{\sqrt{\t}}\Big|^2\Big]
\leq \cC3^M \big(1+\me[|y_T^\pi|^2]\big)
\sum_{m=0}^{M-1} C_{kd+m-1}^m3^m\,.
\label{w604e6}
\end{align}
The same estimate, with the last sum replaced by $1$, holds for $k=0$.

Combining \eqref{w604e5} and \eqref{w604e6}, we arrive at
\begin{align}
\me\int_0^T\big| Y^p_{N,M-1}(t)-Y^p_{N,M-1,K}(t) \big|^2\rd t
\leq \frac{ \cC3^M\big(1+\me[|y_T^\pi|^2]\big)}{K}
\Big[1+\sum_{m=0}^{M-1} 3^m\sum_{k=1}^{N-1}C_{kd+m-1}^m
\Big]\,.
\label{w604e7}
\end{align}
For $m=0,1,\cds,M-1$, it holds that
$\sum_{k=1}^{N-1} C_{kd+m-1}^m
\leq \sum_{r=1}^{d(N-1)} C_{r+m-1}^m
=C_{d(N-1)+m}^{m+1}$.
Furthermore, $C_{d(N-1)+m}^{m+1}\leq C_{dN+M-1}^{M}$.
Consequently,
\begin{align}
\me\int_0^T \big| Y^p_{N,M-1}(t)-Y^p_{N,M-1,K}(t)\big|^2 \rd t
&\leq \frac{ \cC3^M \big(1+\me[|y_T^\pi|^2]\big)}{K} C_{dN+M-1}^{M}
\Big(1+\sum_{m=0}^{M-1}3^m \Big)
\notag\\
&\leq \frac{ \cC\big(1+\me[|y_T^\pi|^2]\big) }{K}
C_{dN+M-1}^{M}9^M\,.
\label{w604e8}
\end{align}
Combining \eqref{wang524b1} and \eqref{w604e8} proves
\eqref{w604e3}.
This completes the proof.
\end{proof}

\section{Algorithm and numerical examples}\label{examples}
In this section, we first present the finite transposition algorithm for BSDEs and then conduct numerical experiments illustrating its performance.

\begin{algorithm}\label{alg1}
Choose $N$, the number of uniform time steps on $[0,T]$; $M$, the order of the Wiener chaos expansion; and $P$, the number of Picard iterations.
Let $\{h_{k,i}\}_{i=1}^{\dim(\cH^M(k))}$ be the basis from Remark \ref{w9r1}.

\begin{enumerate}
\item [{\rm(1)}]
%
Choose an approximation $y_T^\pi\in[\cH^M(N)]^n$ of $y(T)$.
Relying on the orthonormal basis $\{h_{N,i}\}_{i=1}^{\dim(\cH^M(N))}$ introduced in Remark \ref{w9r1}, write
\beq
y^\pi_T=\sum_{i=1}^{\dim(\cH^M(N))}
\frac{a_{N,i}}{\sqrt{\tau }}h_{N,i}\,,
\eeq
where $a_{N,i}\in\dbR^n$ are the corresponding coefficients.

\item [{\rm(2)}] Set $y^0_{N,M}(\cd)=0$,  $p:=1$, $a^p_{N,i}=a_{N,i}$, and  $k:=N-1$.
Compute the coefficients $a^p_{k,i}\in\dbR^n$ and $b^p_{k,i}\in\dbR^{n\times d}$ via \eqref{w602e3} as follows.

\begin{enumerate}
\item [{\rm(2-1)}] Compute $b^p_{k,i}$ for $i=1,2,\cds,\dim\big(\cH^{M-1}(k)\big)$
by
\beq
b^p_{k,i}
=\me\Big[ \frac{h_{k,i}}{\sqrt{\tau}} y^p_{N,M}(t_{k+1})(\D_{k+1}W)^\top\Big]
\in\dbR^{n\times d}\,.
\eeq
Then
\beq
Y^p_{N,M-1}(t_k)=\sum_{i=1}^{\dim(\cH^{M-1}(k))}\frac{b^p_{k,i}}{\sqrt{\tau }}h_{k,i}\,.
\eeq

\item [{\rm(2-2)}] Compute $a^p_{k,i}$ for $i=1,2,\cds, \dim\big(\cH^{M}(k)\big)$
by
\beq
a^p_{k,i}=\me\big[ \sqrt{\tau } h_{k,i}\, y^p_{N,M}(t_{k+1})\big]
-\me\big[ \sqrt{\tau } h_{k,i}\,  f\big(t_k,y^{p-1}_{N,M}(t_k),Y^p_{N,M-1}(t_k)\big)\big] \tau \in\dbR^n\,.
\eeq
Then
\beq
y^p_{N,M}(t_k)=\sum_{i=1}^{\dim(\cH^M(k))}\frac{a^p_{k,i}}{\sqrt{\tau }}h_{k,i}\,.
\eeq

\item [{\rm(2-3)}] When $k=0$, terminate the backward iteration. Otherwise, set $k:=k-1$ and return to Step {\rm(2-1)}.

\end{enumerate}

\item [{\rm(3)}] If $p=P$, terminate the algorithm. Otherwise, set $p:=p+1$, $a^p_{N,i}:=a_{N,i}$ for $i=1,\cds,\dim\big(\cH^M(N)\big)$, and $k:=N-1$, and return to Step {\rm(2-1)}.

\end{enumerate}

\end{algorithm}

Below, we present numerical results for option-hedging problems to illustrate the performance of the finite transposition method developed above.
The first three examples concern scalar BSDEs driven by a one-dimensional Brownian motion, corresponding to $n=d=1$, while the last example concerns a scalar BSDE driven by a five-dimensional Brownian motion, corresponding to $n=1$ and $d=5$.
For the linear BSDEs considered below, no Picard iteration is needed, and we denote the coefficient-wise Monte Carlo approximation directly by $\big(y_{N,M,K}(\cd),Y_{N,M-1,K}(\cd)\big)$. For nonlinear BSDEs, the superscript $p$ is retained to indicate the Picard iteration index.
The computations were performed on a laptop equipped with a 2.7 GHz Intel Core i7 processor and 16 GB of RAM.

In Examples \ref{w624ex1}--\ref{w729ex1}, $x(\cd)$ denotes the one-dimensional stock-price process in the standard Black--Scholes model, satisfying
\bel{w624e1}
\lt\{
\bal
& \mathrm d x(t)=\mu_S x(t)\rd t+\si x(t)\rd W(t)\,,\q t\in [0,T]\,,\\
& x(0)=x_0\,.
\eal
\rt.
\ee
The option-value process $y(\cd)$ is characterized by the following BSDE:
\bel{w624e3}
\lt\{
\bal
&\mathrm d y(t)=\Big[ ry(t)+\frac{\mu_S-r}{\si}Y(t)-(R-r)\big(y(t)-\frac{Y(t)}{\si} \big)_{-}\Big]\rd t\\
&\qq\qq\qq+Y(t)\rd W(t)\,, \q t\in [0,T]\,,\\
&y(T)=\varphi\big(x(T)\big)\,,
\eal
\rt.
\ee
where $\varphi:\dbR\longrightarrow\dbR$ denotes the payoff function, which will be specified in each example below.
We use the Euler method to approximate the solution of SDE \eqref{w624e1} and the finite transposition method to approximate the solution of BSDE \eqref{w624e3}.

\bex{w624ex1}
Let the terminal condition be $y(T)= |x(T)-K_0 |$, and choose the parameters $x_0=100,\, \si=0.2,\, \mu_S=0.05,\,R=r=0.01,\, K_0=100,\, T=2$.
For these parameters, the corresponding Black--Scholes valuation gives the reference values $\big(y(0), Y(0)\big)=(22.32, 3.36)$.
Table \ref{tab1} reports the convergence of $\big(y_{N,M,K}(0),Y_{N,M-1,K}(0)\big)$ as the number $N$ of time steps increases, with $M=2$ and $K=10^5$.
The results indicate that the finite transposition method performs well even with the low Wiener chaos order $M=2$.



%
%
%

\begin{table}[!ht]
\normalfont
\caption{Computed values of $\big(y_{N,M,K}(0),Y_{N,M-1,K}(0)\big)$, the corresponding relative errors, and CPU times for different time partitions, with $M=2$ and $K=10^5$.}
\label{tab1}
\vspace{6pt}
\centering
\setlength{\tabcolsep}{4pt}
\begin{tabular}{@{}lrrrrr@{}}
\hline
$N$
& \multicolumn{1}{c}{$y_{N,M,K}(0)$}
& \multicolumn{1}{c}{Rel. error ($y$)}
& \multicolumn{1}{c}{$Y_{N,M-1,K}(0)$}
& \multicolumn{1}{c}{Relative error ($Y$)}
& \multicolumn{1}{c@{}}{CPU time}\\
\hline
$2^3$ & $22.2487$ & $0.32\%$ & $3.7576$ & $11.83\%$ & $0.069$ s\\
$2^4$ & $22.2929$ & $0.12\%$ & $3.6164$ & $7.63\%$ & $0.169$ s\\
$2^5$ & $22.3265$ & $0.03\%$ & $3.3354$ & $0.73\%$ & $0.497$ s\\
\hline
\end{tabular}
\end{table}

%
%
%

\ex

\bex{w624ex2} {\rm (European call option)}
Let
$y(T)=\max\{x(T)-K_0,0\}$, and take the parameters $x_0=100,\, \si=0.2,\, \mu_S=0.03,\,R=r=0.01,\, K_0=80,\, T=1$. For these parameters, direct evaluation of the Black--Scholes formula gives the reference value $y(0)=21.86$.
In Table \ref{tab2}, we take $ M=2,\, N=2^4$, and report the computed values of $y_{N,M,K}(0)$ as the sample size $K$ increases.

%
%
%

\begin{table}[!ht]
\normalfont
\caption{Computed values of $y_{N,M,K}(0)$, the corresponding relative errors, and CPU times for different numbers of Monte Carlo samples, with $M=2$ and $N=2^4$.}
\label{tab2}
\vspace{6pt}
\centering
\begin{tabular}{@{}lrrr@{}}
\hline
&
\multicolumn{1}{c}{$y_{N,M,K}(0)$}
&
\multicolumn{1}{c}{Relative error}
&
\multicolumn{1}{c@{}}{CPU time}
\\
\hline
$K=10^3$
& $22.1060$
& $1.13\%$
& $0.020$ s
\\
$K=10^4$
& $21.7551$
& $0.48\%$
& $0.033$ s
\\
$K=10^5$
& $21.9010$
& $0.19\%$
& $0.165$ s
\\
\hline
\end{tabular}
\end{table}

%
%
%

\ex

\bex{w729ex1}
In this example, we demonstrate the effectiveness of the finite transposition method for nonlinear BSDEs. Take the terminal condition $y(T)=|x(T)-K_0 |$, and choose the parameters $x_0=100,\, \si=0.2,\, \mu_S=0.05,\,R=0.06,\, r=0.01,\, K_0=100,\, T=2$.
Figure \ref{fig:picard-iteration} shows the evolution of $\big(y^p_{N,M,K}(0),Y^p_{N,M-1,K}(0)\big)$ as the Picard iteration index $p$ increases, with $N=2^5$, $M=2$, and $K=10^5$.
With the iteration number set to $p=6$, the computation takes approximately $10.6$ seconds.
\begin{figure}[tbph]
\centering
\includegraphics[width=1\textwidth,height=0.4\textheight]{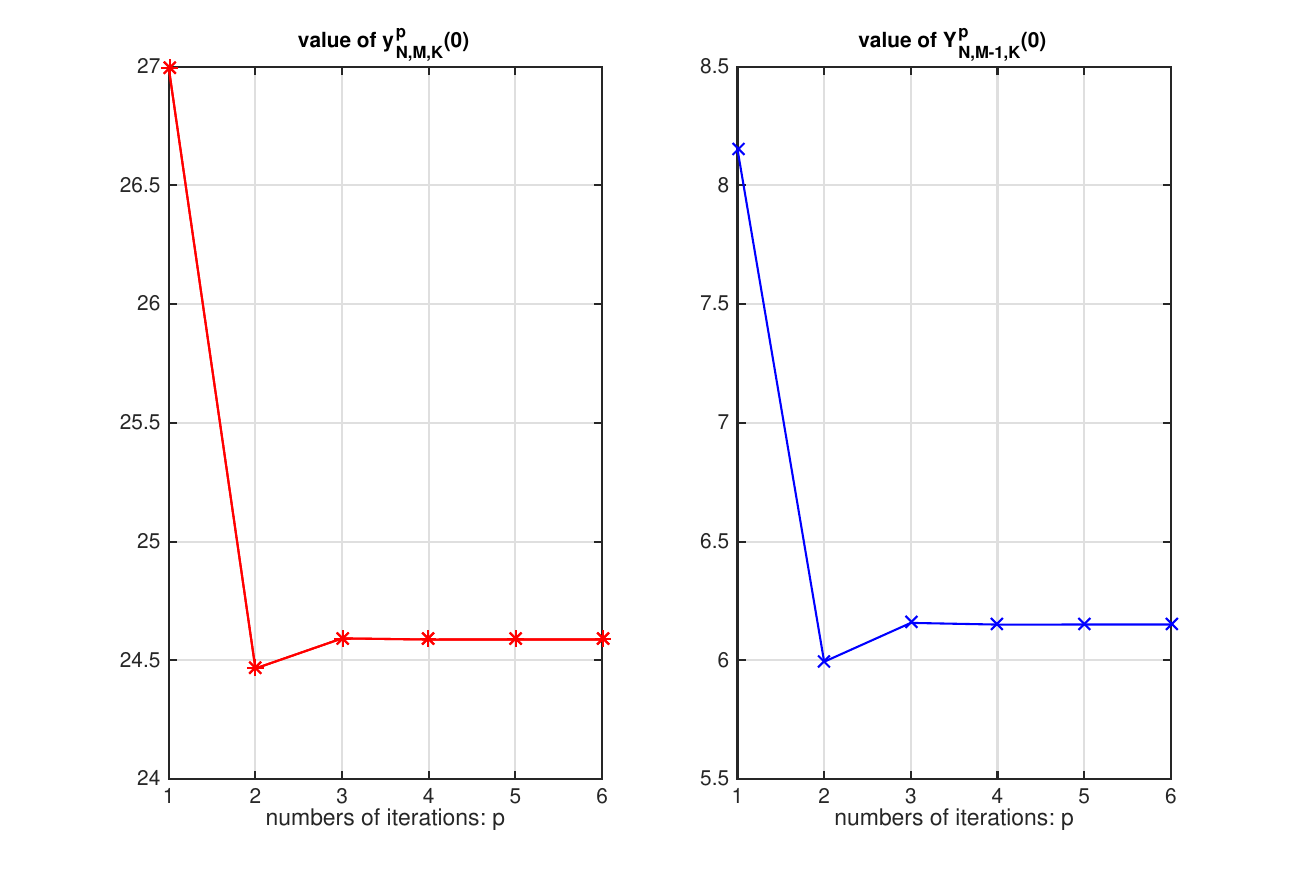}
\caption{Evolution of $\big(y^p_{N,M,K}(0),Y^p_{N,M-1,K}(0)\big)$ as the Picard iteration index $p$ increases, with $N=2^5$, $M=2$, and $K=10^5$.}
\label{fig:picard-iteration}
\end{figure}

\ex

We finally consider an option-pricing problem driven by a five-dimensional Brownian motion to illustrate the performance of the finite transposition method in the multidimensional setting.
In this example, the BSDE is scalar, i.e., $n=1$, while the Brownian dimension is $d=5$.

\bex{w423ex1} {\rm (Basket option)}
Let $W(\cd)=\big(W^1(\cd),W^2(\cd),\cds,W^5(\cd)\big)^\top$ be a five-dimensional Brownian motion.
The price process $x(\cd)=\big(x_1(\cd),x_2(\cd), \cds, x_5(\cd)\big)^\top$ satisfies
\bel{w423e1}
\lt\{
\bal
& \mathrm{d}x(t)=0.11 x(t)\rd t+0.5 \operatorname{diag}\big(x(t)\big)A \rd W(t)\,,\q t\in [0,T]\,,\\
& x(0)=(100,100,\cds,100)^\top\in\dbR^5\,,
\eal
\rt.
\ee
where $A$ is the Cholesky factor of the matrix $(\gamma_{i,j})\in \dbR^{5\times 5}$, with $\gamma_{i,j}=\frac 1 2\delta_{i,j}+\frac 1 2 $. The option-value process $y(\cd)$ satisfies the following BSDE
\bel{w423e2}
\lt\{
\bal
&\mathrm d y(t)=\big[ 0.1 y(t)\!+\!0.02 Y(t)A^{-1}(1,1,1,1,1)^\top \big]\rd t\!+\!Y(t)\rd W(t)\,, \,\, t\in [0,T]\,,\\
&y(T)=\max\big\{ 0.05x_1(T)+0.15x_2(T)+0.2x_3(T)\\
&\qq\qq\qq+0.25x_4(T)+0.35x_5(T)-90, 0\big\}\,.
\eal
\rt.
\ee
For $T=1$ and $T=3$, the corresponding reference prices are 
$25.3757$ and 
$42.7673$, respectively; see \cite{Ju02}.
In Tables \ref{tab3} and \ref{tab4}, we take $ M=2,\, N=2^4$, and report the computed values of $y_{N,M,K}(0)$ for increasing sample sizes $K$.
%
\begin{table}[!ht]
\normalfont
\caption{Computed values of $y_{N,M,K}(0)$, the corresponding relative errors, and CPU times for different numbers of Monte Carlo samples, with $T=1$.}
\label{tab3}
\vspace{6pt}
\centering
\begin{tabular}{@{}lrrr@{}}
\hline
&
\multicolumn{1}{c}{$y_{N,M,K}(0)$}
&
\multicolumn{1}{c}{Relative error}
&
\multicolumn{1}{c@{}}{CPU time}
\\
\hline
$K=10^3$
& $26.3731$
& $3.93\%$
& $0.081$ s
\\
$K=10^4$
& $25.1437$
& $0.91\%$
& $0.326$ s
\\
$K=10^5$
& $25.1523$
& $0.88\%$
& $2.606$ s
\\
\hline
\end{tabular}
\end{table}
\begin{table}[!ht]
\normalfont
\caption{Computed values of $y_{N,M,K}(0)$, the corresponding relative errors, and CPU times for different numbers of Monte Carlo samples, with $T=3$.}
\label{tab4}
\vspace{6pt}
\centering
\begin{tabular}{@{}lrrr@{}}
\hline
&
\multicolumn{1}{c}{$y_{N,M,K}(0)$}
&
\multicolumn{1}{c}{Relative error}
&
\multicolumn{1}{c@{}}{CPU time}
\\
\hline
$K=10^3$
& $44.6774$
& $4.47\%$
& $0.081$ s
\\
$K=10^4$
& $43.8647$
& $2.57\%$
& $0.307$ s
\\
$K=10^5$
& $43.0759$
& $0.72\%$
& $2.558$ s
\\
\hline
\end{tabular}
\end{table}
\ex


\end{document}